\documentclass{amsart}
\usepackage{amssymb}
\usepackage{amscd}
\usepackage{verbatim}
\usepackage{epsfig}
\begin{document}
\newcommand\nin{\notin}
\newcommand\identity{\operatorname{id}}
\newcommand\id{\operatorname{id}}
\newcommand\Id{\operatorname{Id}}
\newcommand\Real{\mathbb{R}}
\newcommand\pos{\Real^+}
\newcommand\Rnp{\Real\setminus\{0\}}
\newcommand\nzero{\setminus\{0\}}
\newcommand\Cx{\mathbb{C}}
\newcommand\Cxp{\Cx^+}
\newcommand\Cxm{\Cx^-}
\newcommand\Nat{\mathbb{N}}
\newcommand\halfNat{{\frac{1}{2}}\mathbb{N}}
\newcommand\intgr{\mathbb{Z}}
\newcommand\im{\operatorname{Im}}
\newcommand\re{\operatorname{Re}}
\newcommand\sign{\operatorname{sign}}
\newcommand\codim{\operatorname{codim}}
\newcommand\End{\operatorname{End}}
\newcommand\Ker{\operatorname{Ker}}
\newcommand\Hom{\operatorname{Hom}}
\newcommand\ideal{{\mathcal I}}
\newcommand\Span{\operatorname{span}}
\newcommand\Range{\operatorname{Ran}}
\newcommand\graph{\operatorname{graph}}
\newcommand\slim{\operatornamewithlimits{s-lim}}
\newcommand\diag{\operatorname{diag}}
\newcommand\Rn{\Real^n}
\newcommand\Rm{\Real^m}
\newcommand\RN{\Real^N}
\newcommand\RtN{\Real^{2N}}
\newcommand\RM{\Real^M}
\newcommand\sphere{\mathbb{S}}
\newcommand\Sn{\sphere^{n-1}}
\newcommand\Sm{\sphere^{m-1}}
\newcommand\Snp{\sphere^n_+}
\newcommand\Smp{\sphere^m_+}
\newcommand\SN{\sphere^{N-1}}
\newcommand\SNp{\sphere^N_+}
\newcommand\circlep{\sphere^1_+}
\newcommand\Phom{P_{h}}
\newcommand\Shom{S_{h}}
\newcommand\distance{\operatorname{dist}}
\newcommand\cl{\operatorname{cl}}
\newcommand\interior{\operatorname{int}}
\newcommand\Fa{\operatorname{Fa}}
\newcommand\ff{\operatorname{ff}}
\newcommand\mf{\operatorname{mf}}
\newcommand\cf{\operatorname{cf}}
\newcommand\scf{\operatorname{sf}}
\newcommand\lf{\operatorname{lb}}
\newcommand\rf{\operatorname{rb}}
\newcommand\indfam{{\mathcal K}}
\newcommand\calX{{\mathcal X}}
\newcommand\calK{{\mathcal K}}
\newcommand\calF{{\mathcal F}}
\newcommand\calO{{\mathcal O}}
\newcommand\calC{{\mathcal C}}
\newcommand\calCL{{\mathcal C}_{\text L}}
\newcommand\calCR{{\mathcal C}_{\text R}}
\newcommand\Cinf{{\mathcal C}^{\infty}}
\newcommand\dist{{\mathcal C}^{-\infty}}
\newcommand\dCinf{\dot\Cinf}
\newcommand\ddist{\dot\dist}
\newcommand\Cj{{\mathcal C}^j}
\newcommand\Linf{L^{\infty}}
\newcommand\bcon{{\mathcal A}}
\newcommand\bconc{{\mathcal A}_{\text{phg}}}
\newcommand\Sch{{\mathcal S}}
\newcommand\temp{\Sch^{\prime}}
\newcommand\Diff{\operatorname{Diff}}
\newcommand\Diffb{\operatorname{Diff}_{\text{b}}}
\newcommand\Diffc{\operatorname{Diff}_{\text{c}}}
\newcommand\Diffsc{\operatorname{Diff}_{\text{sc}}}
\newcommand\DiffI{\operatorname{Diff}_{\text{I}}}
\newcommand\DiffIq{\operatorname{Diff}_{\text{I},q}}
\newcommand\supp{\operatorname{supp}}
\newcommand\ssupp{\operatorname{sing\ supp}}
\newcommand\csupp{\operatorname{cone\ supp}}
\newcommand\esupp{\operatorname{ess\ supp}}
\newcommand\Fr{{\mathcal F}}
\newcommand\Frinv{\Fr^{-1}}
\newcommand\bop{{\mathcal B}}
\newcommand\spec{\operatorname{spec}}
\newcommand\pspec{\spec_{pp}}
\newcommand\cspec{\spec_{c}}
\newcommand\FIO{{\mathcal I}}
\newcommand\SP{\operatorname{SP}}
\newcommand\Symc{S_c}
\newcommand\Symca{S_c^{\alpha}}
\newcommand\Symczero{S_c^{0,...,0}}
\newcommand\sci{{}^{\text{sc}}}
\newcommand\sct{\sci T^*}
\newcommand\scct{\sci\bar{T}^*}
\newcommand\Csc{C_{\text{sc}}}
\newcommand\SNpscd{(\SNp)^2_{\text{sc}}}
\newcommand\scdiag{\Delta_{\text{sc}}}
\newcommand\projscl{\pi^L_{\text{sc}}}
\newcommand\projscr{\pi^R_{\text{sc}}}
\newcommand\scHL{\sci H^{2,0}_{|\zeta|^2-\lambda^2}}
\newcommand\scHrg{\sci H^{2,0}_{\sqrt{g}}}
\newcommand\Hsc{H_{\text{sc}}}
\newcommand\WF{\operatorname{WF}}
\newcommand\WFp{\operatorname{WF^{\prime}}}
\newcommand\WFsc{\operatorname{WF}_{\text{sc}}}
\newcommand\WFscp{\operatorname{WF_{sc}^{\prime}}}
\newcommand\WFC{\operatorname{WF}_C}
\newcommand\WFCi{\operatorname{WF}_{C_i}}
\newcommand\elliptic{\operatorname{ell}}
\newcommand\Psop{\operatorname{\Psi}}
\newcommand\Psiscrs{\operatorname{\Psi_{sc}^{-2,\infty}}}
\newcommand\Psiscr{\operatorname{\Psi_{sc}^{-2,0}}}
\newcommand\Psiscrm{\operatorname{\Psi_{sc}^{0,2}}}
\newcommand\PsiscHam{\operatorname{\Psi_{sc}^{2,0}}}
\newcommand\Psisci{\operatorname{\Psi_{sc}^{*,*}}}
\newcommand\Psiscid{\operatorname{\Psi_{sc}^{0,0}}}
\newcommand\Psiscis{\operatorname{\Psi_{sc}^{0,\infty}}}
\newcommand\Psiscsi{\operatorname{\Psi_{sc}^{-\infty,0}}}
\newcommand\Psiscs{\operatorname{\Psi_{sc}^{-\infty,\infty}}}
\newcommand\Psiscalg{\operatorname{\Psi_{sc}^{\infty,-\infty}}}
\newcommand\nullHam{{\mathcal N}}
\newcommand\charD{\Sigma_{\Delta-\lambda^2}}
\newcommand\charLap{\Sigma_{\Delta-\lambda}}
\newcommand\Snl{\Sn_{\lambda}}
\newcommand\SNl{\SN_{\lambda}}
\newcommand\gammat{\tilde\gamma}
\newcommand\gammasc{\gamma}
\newcommand\taut{\tilde\tau}
\newcommand\PlV{P_V(\lambda)}
\newcommand\PlVc{P_V^{\flat}(\lambda)}
\newcommand\Pl{P_0(\lambda)}
\newcommand\SVl{S_V(\lambda)}
\newcommand\Sjr{S_j^{\reduced}}
\newcommand\Rkp{{\mathcal R}^k_+}
\newcommand\Rkm{{\mathcal R}^k_-}
\newcommand\Rkpm{{\mathcal R}^k_{\pm}}
\newcommand\Phys{{\mathcal P}}
\newcommand\Pc{\overline{\mathcal P}}
\newcommand\pip{\pi^{\perp}}
\newcommand\xit{\tilde\xi}
\newcommand\zetat{\tilde\zeta}
\newcommand\etat{\tilde\eta}
\newcommand\sigmat{\tilde\sigma}
\newcommand\sigmahat{\hat\sigma}
\newcommand\thetat{\tilde\theta}
\newcommand\psit{\tilde\psi}
\newcommand\phit{\tilde\phi}
\newcommand\chit{\tilde\chi}
\newcommand\rhot{\tilde\rho}
\newcommand\xib{\bar\xi}
\newcommand\zetab{\bar\zeta}
\newcommand\thetab{\bar\theta}
\newcommand\etab{\bar\eta}
\newcommand\iotal{\iota_{\lambda}}
\newcommand\rhoat{\rhot_{\alpha_1}}
\newcommand\Lambdat{\tilde\Lambda}
\newcommand\poles{\Lambda'}
\newcommand\rpoles{\Lambda_p}
\newcommand\thresholds{\Lambda}
\newcommand\Vt{\tilde V}
\newcommand\half{{\frac{1}{2}}}
\newcommand\sigmah{\sigma^{1/2}}
\newcommand\bX{\partial X}
\newcommand\Deltabt{\tilde\Delta_0}
\newcommand\strip{\Omega_T}
\newcommand\Kf{K^{\flat}}
\newcommand\Gs{G^{\sharp}}
\newcommand\Gt{\tilde G}
\newcommand\Gb{\bar G}
\newcommand\Osc{\sci\Omega}
\newcommand\Osch{\sci\Omega^{\half}}
\newcommand\Oscmh{\sci\Omega^{-\half}}
\newcommand\Isc{I_{sc}}
\newcommand\Qzl{Q^0_{-\lambda}}
\newcommand\Lie{{\mathcal L}}
\newcommand\bl{{\text b}}
\newcommand\scl{{\text{sc}}}
\newcommand\sccl{{\text{scc}}}
\newcommand\Scl{{\text{3sc}}}
\newcommand\ScLl{{\text{Sc,L}}}
\newcommand\ScRl{{\text{Sc,R}}}
\newcommand\Sccl{{\text{Scc}}}
\newcommand\sfl{{\operatorname{s\Phi}}}
\newcommand\sfi{{}^\sfl}
\newcommand\sus{{\text{sus}}}
\newcommand\Osfh{\sfi\Omega^{\half}}
\newcommand\Isf{I_\sfl}
\newcommand\XXb{X^2_\bl}
\newcommand\XXsc{X^2_\scl}
\newcommand\XXSc{X^2_\Scl}
\newcommand\XXScL{X^2_\ScLl}
\newcommand\XXScR{X^2_\ScRl}
\newcommand\MMsc{M^2_\scl}
\newcommand\Deltab{\Delta_\bl}
\newcommand\Deltasc{\Delta_\scl}
\newcommand\DeltaSc{\Delta_\Scl}
\newcommand\DeltaScL{\Delta_\ScLl}
\newcommand\DeltaScR{\Delta_\ScRl}
\newcommand\prs{\sigma}
\newcommand\Nsc{N_\scl}
\newcommand\Nscp{N_{\scl,p}}
\newcommand\Nff{N_{\ff}}
\newcommand\Nffz{N_{\ff,0}}
\newcommand\Nffzp{N_{\ff,0,p}}
\newcommand\Nffl{N_{\ff,l}}
\newcommand\Nffml{N_{\ff,-l}}
\newcommand\Nmf{N_{\mf}}
\newcommand\Nmfz{N_{\mf,0}}
\newcommand\Nmfl{N_{\mf,l}}
\newcommand\Nmfml{N_{\mf,-l}}
\newcommand\ffb{\operatorname{bf}}
\newcommand\Ffb{\operatorname{bf'}}
\newcommand\ffsc{\operatorname{sf}}
\newcommand\ffSc{\operatorname{sf_C}}
\newcommand\Ffsc{\operatorname{sf'}}
\newcommand\rff{\rho_{\ff}}
\newcommand\rmf{\rho_{\mf}}
\newcommand\rffb{\rho_{\ffb}}
\newcommand\rffsc{\rho_{\ffsc}}
\newcommand\rFfsc{\rho_{\Ffsc}}
\newcommand\rffSc{\rho_{\ffSc}}
\newcommand\rinf{\rho_{\infty}}
\newcommand\CL{C_L}
\newcommand\CR{C_R}
\newcommand\betab{\beta_\bl}
\newcommand\betasc{\beta_\scl}
\newcommand\betaSc{\beta_\Scl}
\newcommand\BetaSc{\bar\beta_\Scl}
\newcommand\betaScL{\beta_\ScLl}
\newcommand\betaScR{\beta_\ScRl}
\newcommand\sfT{{}^\sfl T^*}
\newcommand\sfN{{}^\sfl N^*}
\newcommand\ScT{{}^\Scl T^*}
\newcommand\SccT{{}^\Scl \bar T^*}
\newcommand\ScS{{}^\Scl S^*}
\newcommand\scS{{}^\scl S^*}
\newcommand\Tb{{}^\bl T}
\newcommand\Tsc{{}^\scl T}
\newcommand\Tsf{{}^\sfl T}
\newcommand\TSc{{}^\Scl T}
\newcommand\CSc{C_\Scl}
\newcommand\Lambdasc{{}^\scl\Lambda}
\newcommand\XXXb{X^3_\bl}
\newcommand\XXXsc{X^3_\scl}
\newcommand\XXXSc{X^3_\Scl}
\newcommand\XXXScO{X^3_{\Scl,O}}
\newcommand\XXXScF{X^3_{\Scl,F}}
\newcommand\XXXScS{X^3_{\Scl,S}}
\newcommand\XXXScC{X^3_{\Scl,C}}
\newcommand\KDsc{\operatorname{KD^{\half}_\scl}}
\newcommand\SDsc{\operatorname{SD^{\half}_\scl}}
\newcommand\SDsf{\operatorname{SD^{\half}_\sfl}}
\newcommand\KDSc{\operatorname{KD^{\half}_\Scl}}
\newcommand\KDScEF{\operatorname{KD^{E,F}_\Scl}}
\newcommand\Oh{\operatorname{\Omega^{\half}}}
\newcommand\WFSc{\WF_\Scl}
\newcommand\WFScmf{\WF_{\Scl,\mf}}
\newcommand\WFScff{\WF_{\Scl,\ff}}
\newcommand\WFScs{\WF_{\Scl,\prs}}
\newcommand\WFScp{\WF'_\Scl}
\newcommand\WFScmfp{\WF'_{\Scl,\mf}}
\newcommand\WFScffp{\WF'_{\Scl,\ff}}
\newcommand\WFScsp{\WF'_{\Scl,\prs}}
\newcommand\Diffscc{\Diff_\sccl}
\newcommand\DiffSc{\Diff_\Scl}
\newcommand\Diffsf{\Diff_\sfl}
\newcommand\DiffScc{\Diff_\Sccl}
\newcommand\DiffscI{\Diff_{\scl,\text{I}}}
\newcommand\VscI{\Vf_{\scl,\text{I}}}
\newcommand\DiffsV{\operatorname{Diff}_{\sus(V)}}
\newcommand\DiffsVsc{\operatorname{Diff}_{\sus(V),\scl}}
\newcommand\DiffsVCsc{\operatorname{Diff}_{\sus(V)-C,\scl}}   
\newcommand\Psisc{\Psop_\scl}
\newcommand\Psiscc{\Psop_\sccl}
\newcommand\PsiSc{\Psop_\Scl}
\newcommand\PsiScc{\Psop_\Sccl}
\newcommand\PsiSccml{\Psop^{m,l}_\Sccl}
\newcommand\PsiScxx{\Psop^{*,*}_\Scl}
\newcommand\PsiScml{\Psop^{m,l}_\Scl}
\newcommand\PsiScmz{\Psop^{m,0}_\Scl}
\newcommand\PsiScmmz{\Psop^{-m,0}_\Scl}
\newcommand\PsiSckz{\Psop^{k,0}_\Scl}
\newcommand\PsiScmmml{\Psop^{-m,-l}_\Scl}
\newcommand\Psiscmkk{\Psop^{-k,k}_\scl}
\newcommand\Psiscmmmkk{\Psop^{-m-k,k}_\scl}
\newcommand\Psiscmoo{\Psop^{-1,1}_\scl}
\newcommand\Psiscmz{\Psop^{m,0}_\scl}
\newcommand\Psiscmmz{\Psop^{-m,0}_\scl}
\newcommand\PsiSckmkl{\Psop^{km,kl}_\Scl}
\newcommand\PsiScmplp{\Psop^{m',l'}_\Scl}
\newcommand\PsiScmmpllp{\Psop^{m+m',l+l'}_\Scl}
\newcommand\Psiscml{\Psop^{m,l}_\scl}
\newcommand\PsiScid{\Psop^{0,0}_\Scl}
\newcommand\PsiSczo{\Psop^{0,1}_\Scl}
\newcommand\PsiScmii{\Psop^{-\infty,\infty}_\Scl}
\newcommand\PsiScmiz{\Psop^{-\infty,0}_\Scl}
\newcommand\PsiScmoo{\Psop^{-1,1}_\Scl}
\newcommand\PsisCid{\Psop^{0,0}_{\scl-C}}
\newcommand\PsisC{\Psop_{\scl-C}}
\newcommand\Psiinf{\Psop_{\infty}}
\newcommand\Psiinfid{\Psop_{\infty}^0}
\newcommand\PsiFinf{\Psop_{\infty-\Fr}}
\newcommand\PsisVscml{\Psop^{m,l}_{\sus(V),\scl}}
\newcommand\PsisVsc{\Psop_{\sus(V),\scl}}
\newcommand\PsisVpsc{\Psop_{\sus(V_p),\scl}}
\newcommand\PsisVCSc{\Psop_{\sus(V)-C,\scl}}
\newcommand\SFinf{S_{\infty-\Fr}}
\newcommand\YsVC{Y^2_{\sus(V)-C,\scl}}
\newcommand\ffYsc{\ffsc_{\sus(V)}}
\newcommand\SXC{S(X;C)}
\newcommand\Ios{I_{\text{os}}}
\newcommand\pbL{\pi^2_{\bl,{\text L}}}
\newcommand\pbR{\pi^2_{\bl,{\text R}}}
\newcommand\pscL{\pi^2_{\scl,{\text L}}}
\newcommand\pscR{\pi^2_{\scl,{\text R}}}
\newcommand\PbO{\pi^3_{\bl,{\text O}}}
\newcommand\PscO{\pi^3_{\scl,{\text O}}}
\newcommand\PScO{\pi^3_{\Scl,{\text O}}}
\newcommand\PScF{\pi^3_{\Scl,{\text F}}}
\newcommand\PScC{\pi^3_{\Scl,{\text C}}}
\newcommand\PScS{\pi^3_{\Scl,{\text S}}}
\newcommand\pScL{\pi^2_{\Scl,{\text L}}}
\newcommand\pScR{\pi^2_{\Scl,{\text R}}}
\newcommand\CLF{\CL^F}
\newcommand\CLO{\CL^O}
\newcommand\CLS{\CL^S}
\newcommand\CLC{\CL^C}
\newcommand\DeltaYb{\Delta_{\bl,Y}}
\newcommand\DeltaYsc{\Delta_{\sus-\scl}}
\newcommand\Vf{{\mathcal V}}
\newcommand\Vb{{\mathcal V}_{\bl}}
\newcommand\Vsc{{\mathcal V}_{\scl}}
\newcommand\VSc{{\mathcal V}_{\Scl}}
\newcommand\Vsf{{\mathcal V}_{\sfl}}
\newcommand\VfI{\Vf_{\text{I}}}
\newcommand\VfIq{\Vf_{\text{I},q}}
\newcommand\scH{{}^\scl H}
\newcommand\scHg{\scH_g}
\newcommand\xh{\hat x}
\newcommand\xb{\bar x}
\newcommand\Yh{\hat Y}
\newcommand\Yb{\bar Y}
\newcommand\hb{\bar h}
\newcommand\xih{\hat\xi}
\newcommand\etah{\hat\eta}
\newcommand\muh{\hat\mu}
\newcommand\mut{\tilde\mu}
\newcommand\mub{\bar\mu}
\newcommand\mubh{\widehat{\bar\mu}}
\newcommand\yb{\bar y}
\newcommand\phib{\bar \phi}
\newcommand\ub{\bar u}
\newcommand\Qb{\bar Q}
\newcommand\Wbp{{\bar W}^\perp}
\newcommand\Wp{W^\perp}
\newcommand\Kt{\tilde K}
\newcommand\Wt{\tilde W}
\newcommand\Ut{\tilde U}
\newcommand\xt{\tilde x}
\newcommand\yt{\tilde y}
\newcommand\ft{\tilde f}
\newcommand\fs{f^{\sharp}}
\newcommand\at{\tilde a}
\newcommand\htil{\tilde h}
\newcommand\gt{\tilde g}
\newcommand\Ht{\tilde H}
\newcommand\Mt{\tilde M}
\newcommand\St{\tilde S}
\newcommand\Pt{\tilde P}
\newcommand\Rt{\tilde R}
\newcommand\qt{\tilde q}
\newcommand\Qt{\tilde Q}
\newcommand\Xb{\bar X}
\newcommand\lambdat{\tilde\lambda}
\newcommand\epst{\tilde\epsilon}
\newcommand\At{\tilde A}
\newcommand\Ah{\hat A}
\newcommand\Bh{\hat B}
\newcommand\Gh{\hat G}
\newcommand\Hh{\hat H}
\newcommand\Qh{\hat Q}
\newcommand\Ph{\hat P}
\newcommand\Nh{\hat N}
\newcommand\Sh{\hat S}
\newcommand\Gcal{{\mathcal G}}
\newcommand\GcalC{{\mathcal G}_C}
\newcommand\Jcal{{\mathcal J}}
\newcommand\JcalC{{\mathcal J}_C}
\newcommand\Miff{\ \text{iff}\ }
\newcommand\Mif{\ \text{if}\ }
\newcommand\Mand{\ \text{and}\ }
\newcommand\Mor{\ \text{or}\ }
\newcommand\Mst{\ \text{s.t.}\ }
\setcounter{secnumdepth}{3}
\newtheorem{lemma}{Lemma}[section]
\newtheorem{prop}[lemma]{Proposition}
\newtheorem{thm}[lemma]{Theorem}
\newtheorem{cor}[lemma]{Corollary}
\newtheorem{result}[lemma]{Result}
\newtheorem*{thm*}{Main Results}
\numberwithin{equation}{section}
\theoremstyle{remark}
\newtheorem{rem}[lemma]{Remark}
\theoremstyle{definition}
\newtheorem{Def}[lemma]{Definition}
\def\signature#1#2{\par\noindent#1\dotfill\null\\*
{\raggedleft #2\par}}

\newcommand\RR{\mathbb{R}}
\newcommand\CC{\mathbb{C}}
\newcommand\Lap{\Delta}
\newcommand\pa{\partial}
\newcommand\ep{\epsilon}
\renewcommand\l{\lambda}
\newcommand\symbint{\sigma_{\operatorname{int}}}
\newcommand\symbb{\sigma_{\pa}}
\newcommand\Dprime{{\mathcal D}'}
\newcommand\td{\tilde}
\def\ang#1#2{\big\langle #1 , #2 \big\rangle}
\newcommand\ol{\overline}
\newcommand\hz{\hat z}
\newcommand\Nu{{\mathcal V}}
\renewcommand\sp{\operatorname{sp}}
\newcommand\lb{\operatorname{lb}}
\newcommand\rb{\operatorname{rb}}
\newcommand\bfc{\operatorname{bf}}
\newcommand\ul{\underline}
\newcommand\Ls{L^{\sharp}}
\newcommand\Lt{\tilde L}
\newcommand\Imag{\operatorname{Im}}

\title[Spectral projections and resolvent]
{The spectral projections and the resolvent for scattering metrics}
\author{Andrew Hassell}
\address{Centre for Mathematics and its Applications, Australian National University,
Canberra ACT 0200 Australia}
\email{hassell@maths.anu.edu.au}
\author[Andras Vasy]{Andr\'as Vasy}
\address{Department of Mathematics, University of California, Berkeley,
CA 94720, U.S.A.}
\email{andras@math.berkeley.edu}
\date{\today}
\begin{abstract}

In this paper we consider a compact manifold with boundary $X$ equipped with a
scattering metric $g$ as defined by Melrose \cite{RBMSpec}. That is, $g$ is a
Riemannian metric in the interior of $X$ that can be brought to the form
$g=x^{-4}\,dx^2+x^{-2}h'$ near the boundary, where $x$ is a boundary defining
function and $h'$ is a smooth symmetric 2-cotensor which
restricts to a metric $h$ on $\pa X$. Let $H=\Delta+V$ where $V\in
x^2\Cinf(X)$ is real, so
$V$ is a `short-range' perturbation of $\Delta$.
Melrose and Zworski started a detailed analysis of various operators
associated to $H$
in \cite{RBMZw} and showed that the scattering
matrix of $H$ is a Fourier integral operator associated to the geodesic flow
of $h$ on $\bX$ at distance $\pi$ and that the kernel of the Poisson operator
is a Legendre distribution on $X\times\bX$ associated to an intersecting
pair with conic points. In this paper we describe the kernel of the
spectral projections and the resolvent,
$R(\sigma \pm i0)$, on the positive real axis. We define a class of
Legendre distributions on certain types of manifolds with
corners, and show that the kernel of the spectral
projection is a Legendre distribution associated to a conic
pair on the b-stretched product $\XXb$ (the blowup of $X^2$ about the corner,
$(\pa X)^2$). The structure of the resolvent is only slightly more complicated.

As applications of our results we show that there are `distorted Fourier transforms'
for $H$, ie, unitary operators which intertwine $H$ with a multiplication operator
and determine the scattering matrix; and give a scattering wavefront set estimate
for the resolvent $R(\sigma \pm i0)$ applied to a distribution $f$. 
\end{abstract}
\maketitle

\section{Introduction}

In this paper we study the basic analytic operators associated to short range
Schr\"odinger operators $H$ on a manifold with boundary, $X$, with scattering metric.
The analytic operators of interest are the Poisson operator, scattering matrix,
spectral projections and resolvent family. These will be defined and described in
detail later in the introduction. The first two operators were analyzed rather
completely by Melrose and Zworski \cite{RBMZw} and we will use their analysis 
extensively in this paper. 

A scattering metric $g$ on a manifold with boundary $X^n$ is a smooth 
Riemannian metric in the
interior of $X$ which can be brought to the form
\begin{equation}
g=\frac{dx^2}{x^4}+\frac{h'}{x^2}
\label{scmetric}
\end{equation}
near $\bX$ for some choice of a boundary defining function
$x$ and for some smooth symmetric
2-cotensor $h'$ on $X$ which restricts to a metric
$h$ on $\bX$. We note that this fixes $x$ up to $x^2\Cinf(X)$.
To understand the geometry that such a metric endows, suppose that for some choice
of product coordinates $(x,y)$ near the boundary, $g$ takes the warped product form
$$
\frac{dx^2}{x^4}+\frac{h(y)}{x^2} \quad \text{for } x \leq x_0.
$$
Then changing variable to $r = 1/x$, this looks like $dr^2 + r^2 h(y)$ for $r \geq r_0$,
which is the `large end of a cone'. A general scattering metric may be considered a
`short range' perturbation of this geometry. A particularly important example is $\RR^n$, 
with the flat metric, which after radial compactification and after introduction of 
the boundary defining function $x = 1/r$ takes the form \eqref{scmetric} with $h$ the
standard metric on the sphere. 

Let $\Delta$ be the (positive) Laplacian of $g$, and let $V \in x^2 \Cinf(X)$ be a real
valued function. We will study operators of the form $H = \Delta + V$ which we consider
a `short range Schr\"odinger operator' by analogy with the Euclidean situation. The
case $V \equiv 0$ is already of considerable interest. Since the interior of $X$
is a complete Riemannian manifold, $H$ is essentially self-adjoint on $L^2(X)$
(\cite{Taylor}, chapter 8). 

One of the themes of this paper is that the operator $\Lap + V$ on $\RR^n$ is very
typical of a general $H$ as above. Melrose and Zworski \cite{RBMZw} showed that
the scattering-microlocal structure 
of the Poisson operator and scattering matrix 
can be described purely in terms of geodesic flow on the
boundary $\pa X$. This is equally true of the spectral projections and resolvent,
and looked at from this point of view the
case of $\RR^n$ is a perfect guide to the general situation. 

Our results are phrased in terms of Legendrian distributions introduced by Melrose
and Zworski. We need to extend their definitions to the case of manifolds with
corners of codimension 2, since the Schwartz kernels of the spectral projections and
resolvent are defined on $X^2$ which has codimension 2 corners. One
of our motivations for writing this paper is to understand Legendrian distributions
on manifolds with corners, in the belief that similar techniques
may be used to do analogous constructions for $N$ body Schr\"odinger
operators.
For the three body problem, Legendrian techniques have already proved useful
--- see \cite{Hassell:Plane}.

\

The structure of generalized eigenfunctions of a scattering metric has been described
by Melrose \cite{RBMSpec}. He showed 
that for any $\lambda\in\Real\nzero$ and $a\in\Cinf(\bX)$ there is a
unique $u\in\dist(X)$
which satisfies $(H-\lambda^2)u=0$ and which is of the form
\begin{equation}\label{eq:smooth-exp-1}
u=e^{-i\lambda/x}x^{(n-1)/2}v_-+e^{i\lambda/x}x^{
(n-1)/2}v_+,\ v_\pm\in\Cinf(X),\ v_-|_{\bX}=a.
\end{equation}
Of course, such a representation has been known for $\Lap + V$ on $\RR^n$ for a
long time. 
Note that changing the sign of $\lambda$ reverses the role of $v_+$ and
$v_-$, and correspondingly the role of the `boundary data' $v_+|_{\bX}$
and $v_-|_{\bX}$.
Thus, if we wish to take $\lambda>0$, this statement says that
to get a generalized eigenfunction $u$ of $H$ of this form,
we can either prescribe $v_+|_{\bX}$ or $v_-|_{\bX}$ freely, but the
prescription of either of these determines $u$ uniquely. In particular,
prescribing $v_+|_{\bX}$ determines $v_-|_{\bX}$ and vice versa.

The Poisson operator is defined as the map from boundary data to generalized
eigenfunctions of $H$ \cite{RBMZw}. Thus, with the notation of \eqref{eq:smooth-exp-1},
for $\lambda\in\Real\nzero$, $P(\lambda):\Cinf(\bX)\to\dist(X)$ is given by
$$
P(\lambda)a=u.
$$
Moreover, the scattering matrix $S(\lambda):\Cinf(\bX)
\to\Cinf(\bX)$ is the map
\begin{equation}
S(\lambda)a=v_+|_{\bX}.
\end{equation}
Thus, $S(\lambda)v_-|_{\bX}=v_+|_{\bX}$, i.e.\ the scattering matrix links the
parametrizations of generalized eigenfunctions $u$
of $H$ by the boundary data $v_-|_{\bX}$ and $v_+|_{\bX}$.
We remark that the normalization of $P(\lambda)$ and $S(\lambda)$
here is the opposite of the one used in \cite{RBMSpec} (i.e.\ the roles
of $S(\lambda)$ and $S(-\lambda)$ are interchanged); it matches that
of \cite{RBMZw}.

The kernels of these operators are distributions on $X \times \pa X$ and
$\pa X \times \pa X$ respectively. 
Both kernels were analyzed in \cite{RBMZw} using scattering microlocal
analysis (described in the next section). The authors showed that the 
kernel of the
Poisson operator is a conic Legendrian pair which can be described in terms of
geodesic flow on $\pa X$ and that the kernel of the
scattering matrix is a Fourier integral operator whose canonical relation is 
geodesic flow at time $\pi$ on $\pa X$. 

Melrose also showed that the resolvent family $R(\sigma) = (H - \sigma)^{-1}$ is
in the scattering calculus (a class of pseudodifferential operators on $X$
defined in \cite{RBMSpec}) for $\im \sigma \neq 0$. We are interested in this paper in 
understanding the behaviour of the resolvent kernel as $\im \sigma \to 0$ when
$\re\sigma \neq 0$. This is related to the spectral measure of $H$ via
Stone's theorem
\begin{equation}
\text{s-}\lim_{\ep \to 0} \frac1{2\pi i} \int_a^b R(\sigma + i\ep) -
R(\sigma - i\ep) \, d\sigma = \frac{1}{2} \Big( E_{(a,b)} + E_{[a,b]} \Big).
\label{Stone}
\end{equation}
We will show in Lemma~\ref{Poisson=sp} that in fact
\begin{equation}\label{eq:intro-55}
R(\lambda^2+i0)-R(\lambda^2-i0)=\frac{i}{2\lambda}P(\lambda)P(\lambda)^*
\end{equation}
as operators $\dCinf(X)\to\dist(X)$. Modulo a few technicalities 
(see Lemma~\ref{continuity})
this shows that the spectral measure $dE(\sigma)$
of $H$ is differentiable as a map from
$\dCinf(X)$ to $\dist(X)$, and we may write
$$
(dE)(\l^2) \equiv 2\l \sp(\l) \, d\l = \frac1{2\pi} P(\lambda)P(\lambda)^* d\l, \quad
\l > 0
$$
for some operator $\sp(\l) : \dCinf(X) \to \dist(X)$ which we call the
generalized spectral projection, or (somewhat incorrectly) just the spectral
projection at energy $\l^2 > 0$. Thus $\sp(\l) =
(4\pi \l)^{-1}P(\lambda)P(\lambda)^*$. 
Our approach in analyzing the kernel of the resolvent is to understand the
microlocal nature of the composition $P(\lambda)P(\lambda)^*$ and then
to construct the resolvent using a regularization of the formal expression
\begin{equation}
R(\l^2 \pm i0) = 
\int_{-\infty}^\infty (\sigma - \l^2 \mp i0)^{-1} dE(\sigma) . 
\end{equation}

\begin{figure}\centering
\epsfig{file=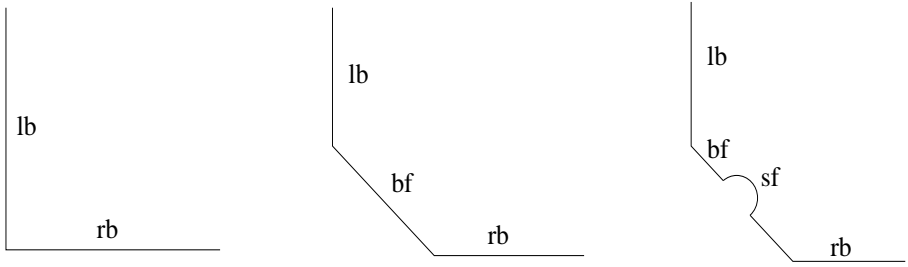,width=12cm,height=3.5cm}
\caption{Blowing up to produce the double scattering space}
\label{fig:XXb}
\end{figure}

Our results on the spectral projection and the resolvent are stated
precisely in Theorems~\ref{thm:sp-proj} and
Theorem~\ref{thm:resolvent}. Here we will give an informal description. Let
$\XXb$ denote the blow-up of
$X^2=X\times X$ about its corner $(\bX)^2$ (see Section~\ref{sec:sc-calculus}
for a detailed description of the space, and Figure~\ref{fig:XXb} for
a picture). Then we show 

\begin{thm*} (1) The kernel of the spectral projection $\sp(\l)$ is a
Legendre distribution associated to an intersecting pair $(L(\l), \Ls(\l))$
with conic points. Here $L(\l)$ and $\Ls(\l)$ are two Legendrians contained
in the boundary of a `stretched cotangent bundle' over $\XXb$ which can be defined
purely geometrically (see equations \eqref{eq:sp-1} and
\eqref{eq:L-sharp}). 

(2) The kernel of the boundary value $R(\lambda^2+i0)$ of the resolvent of $H$ on the
positive real axis is of the form $R_1+R_2+R_3$, where $R_1$ is in
the scattering calculus (and thus qualitatively similar to the resolvent off
the real axis), $R_2$ is an intersecting Legendre distribution, supported
near $\partial\diag_\bl$ and $R_3$ is a Legendre distribution associated to 
$(L(\l), \Ls(\l))$. 
\end{thm*}

\

The paper is organized as follows. In section~\ref{sec:sc-calculus} we review
elements of the scattering calculus, including Legendrian distributions.
Sections~\ref{sec:corner-geo} and \ref{sec:corner-Leg-dist} are devoted to
the definition and parametrization of Legendre distributions on manifolds
with codimension two corners with a `stratification' of the boundary,
in a sense made precise there. In section~\ref{sec:Poisson} we recall Melrose
and Zworski's results on the Poisson operator and prove \eqref{eq:intro-55}. 
Sections~\ref{sec:sp} and \ref{sec:resolvent} contain the statements and proof
of the main results, Theorems~\ref{thm:sp-proj} and \ref{thm:resolvent}. 
This is preceded by section~\ref{sec:example} which discusses the Euclidean
spectral projection and resolvent in Legendrian terms, and turns out to be a
useful guide to the general situation. 

The final two sections contain applications of the main results. In 
section~\ref{sec:distorted-FT}, we show that
the Poisson operator can be viewed as a `distorted Fourier transform' for $H$. 
Let us define the operators 
$P^*_\pm :\Cinf_c(X) \to \dist(\pa X \times \RR_+)$ by
the formula
\begin{equation}
(P^*_\pm u)(y,\l) = (2\pi)^{-1/2} (P(\pm \l)^* u)(y).
\label{P^*pm}
\end{equation}
Then $P^*_\pm$ extends to an isometry from $H_{ac}(H)$, the absolutely continuous
spectral subspace of $H$, to $L^2(\pa X \times \RR_+)$, with adjoint 
\begin{equation}
(P_\pm f)(z) = (2\pi)^{-1/2} \int_0^\infty P(\pm\l)(f(\cdot, \l)) \, d\lambda .
\label{Ppm}
\end{equation}
If the operator $S$ is defined on $\Cinf(\pa X \times \RR_+)$ by
\begin{equation}
(Sf)(y,\l) = (S(\l) f(\cdot, \l))(y),
\end{equation}
then we find that $S = P_- P^*_+$, analogous to a standard formula in scattering
theory for the scattering matrix in terms of the distorted Fourier transforms. 

In section~\ref{sec:wavefront-reln}, we derive a bound for the scattering wavefront
set of $R(\l^2 \pm i0) f$ in terms of the scattering wavefront set of $f$. This
is not a new result, since Melrose proved it in \cite{RBMSpec} using positive
commutator estimates. However, we wish to emphasize that it follows from
a routine calculation, once one 
understands the Legendrian structure of the resolvent.

In future work, we plan to outline a symbol calculus for 
these types of distributions which can
be used to construct the resolvent kernel `directly', that is, without going
via the spectral projections. We also plan to use our results here to get 
wavefront set bounds on the Schr\"odinger kernel $e^{itH}$, and the wave
kernel $e^{it\sqrt{H}}$, applied to suitable distributions $f$. 

\

{\it Notation and conventions. } On a compact manifold with boundary, $X$, we
use $\dCinf(X)$ to denote the class of smooth functions, all of whose derivatives
vanish at the boundary, with the usual topology, and $\dist(X)$ to denote its
topological dual. On the radial compactification of $\RR^n$ these correspond to the
space of Schwartz functions and tempered distributions, respectively. The Laplacian
$\Delta$ is taken to be positive. The space $L^2(X)$ is taken with respect to the
Riemannian density induced by the scattering metric $g$, thus has the form
$a \, dx dy/x^{n+1}$ near the boundary, where $a$ is smooth. 

The outgoing resolvent $R(\l^2 + i0)$, $\l > 0$, is such that application to a
$\dCinf(X)$ function yields a function of the form 
$$
R(\l^2 + i0) = x^{(n-1)/2} e^{i\l/x} a(x,y) , \quad a \in \Cinf(X).
$$
This corresponds to $\tau$ being {\it negative} in the scattering wavefront set.

If $M$ is a manifold with corners and $S$ is a p-submanifold (that is, the inclusion
of $S$ in $M$ is smoothly modelled on $\RR^{n',k'} \subset \RR^{n,k}$), then
$[M;S]$ denotes the blowup of $S$ inside $M$ (see \cite{RBMCalcCon}). The new
boundary hypersurface thus created is called the `front face'. 

{\it Acknowledgements. } We wish to thank Richard Melrose and Rafe Mazzeo for suggesting
the problem, and Richard Melrose in particular for many very helpful conversations. 
A.H. also thanks the MIT mathematics department
for its hospitality during several visits, and the Australian Research Council for
financial support. Finally, we thank the referee for helpful comments. 

\section{Scattering Calculus and Legendrian distributions}\label{sec:sc-calculus}

In this section we describe very briefly elements of the scattering calculus. Further 
information can be found in \cite{RBMSpec}.

A scattering metric $g$ on a manifold with boundary $X$ determines a class of `sc-vector
fields', denoted $\Vsc(X)$, namely those $\Cinf$ vector fields which are of finite 
length with respect to $g$. Locally this consists of all $\Cinf$ vector fields
in the interior of $X$, and near the boundary, is the $\Cinf(X)$-span of the
vector fields $x^2 \pa_x$ and $x \pa_{y_i}$, as can be seen from \eqref{scmetric}.
The set $\Vsc(X)$ is in fact a Lie Algebra. From it we may define $\Diffsc^k(X)$, the
set of differential operators which can be expressed as a finite sum of products of at
most $k$ sc-vector fields with $\Cinf(X)$-coefficients. 

In \cite{RBMSpec} 
Melrose constructed a natural `microlocalization' $\Psisc(X)$ of $\Diffsc(X)$,
ie an algebra of operators which extends the differential operators in a similar way to
that in which pseudodifferential operators on a closed manifold extend the differential
operators. The operators are described by specifying their distributional kernels. 
Since pseudodifferential operators on a closed manifold have kernels conormal to the
diagonal, elements of $\Psisc(X)$ should be conormal to the diagonal in the interior
of $X^2$; what remains to be specified is their behaviour at the boundary. For this
purpose, $X^2$ has rather singular geometry, since the diagonal intersects the
boundary in a non-normal way. This difficulty is overcome by blowing up $X^2$ to
resolve the singularities. There are two blowups to be performed in order to
describe the kernels of $\Psisc(X)$. First we form the b-double space
$$
\XXb = [X^2; (\pa X)^2].
$$
The diagonal lifts to a product-type submanifold $\diag_\bl$. We next blow up the boundary
of this submanifold, obtaining the scattering double space $\XXsc$:
$$
\XXsc = [\XXb; \pa \diag_\bl].
$$
The lift of the diagonal to $\XXsc$ is denoted $\diag_\scl$. The boundary hypersurfaces
are denoted lb, rb, bf and sc which arise from the boundaries $\pa X \times X$,
$X \times \pa X$, the first blowup and the second blowup, respectively (see figure).

It is not hard to check that $\Diffsc(X)$ corresponds precisely to all distributions
on $\XXsc$ conormal to, and supported on, $\diag_\scl$, which are smooth up to the
boundary (using the Riemannian density to define integration against functions). 
Melrose defines $\Psisc^{k,0}(X; \Osch)$ to be the class of kernels on $\XXsc$ of the
form $\kappa |d\mu|^{1/2}$ where $\kappa$ is conormal of order $k$ to $\diag_\scl$,
smooth up to sc, and rapidly decreasing at lb, rb and bf, and $|d\mu|$ is the 
product of Riemannian densities from the left and right factors of $X^2$, lifted to
$\XXsc$. These kernels act naturally on half-densities on $X$. One can show that
they map $\dCinf$ half-densities to $\dCinf$ half-densities, and thus, by duality,
also map $\dist$ half-densities to $\dist$ half-densities. The space
$\Psisc^{k,l}(X)$ is defined to be $x^l \Psisc^{k,0}(X)$.  
See \cite{RBMSpec} for more 
details, or \cite{Hassell-Vasy:Symbolic} 
for a more detailed description where $X$ is the compactification of $\RR^n$.

We next describe the symbols of these operators. First, we need to consider a rescaled
version of the cotangent bundle which matches the degeneration of the sc-vector fields
at the boundary. Since $\Vsc(X)$ can always be expressed locally as the 
$\Cinf(X)$-span of $n$ local sections, there is a vector
bundle $\Tsc X$ whose smooth sections are exactly the sc-vector
fields: $\Vsc(X)=\Cinf(X;\Tsc X)$. The dual bundle of $\Tsc X$ is denoted
by $\sct X$, is called the scattering cotangent bundle, and a basis of
its smooth sections is given by $\frac{dx}{x^2}$, $\frac{dy_j}{x}$, $j=1,
\ldots, N-1$, in the local coordinates described above.
We can write a covector $v\in\sct_p M$ as $v=\tau\frac{dx}{x^2}
+\mu\cdot\frac{dy}{x}$, and thereby we have coordinates $(x,y,\tau,\mu)$
on $\sct X$ near $\partial X$. Notice that $\tau = x^2 \xi$ and $\mu = x \eta$ in terms
of the dual cotangent variables $(\xi, \eta)$ to $(x, y)$, so we will sometimes refer
to $(\tau, \mu)$ as rescaled cotangent variables. 
The half-density bundle induced by $\sct X$ is denoted by $\Osch X$.

We may alternatively think of sc-vector fields as smoothly varying linear forms on
$\sct X$. Since $\Vsc(X)$ is a Lie Algebra, there is a well-defined notion of the
$k$th degree part of an operator $P \in \Diffsc^k(X)$; its symbol is the corresponding
homogeneous polynomial of degree $k$ on the fibres of $\sct X$. This is the
interior principal symbol of $P$, denoted $\symbint(P)$. The symbol map extends
to $\symbint : \Psisc^{k,0}(X) \to S^k(\sct X)/S^{k-1}(\sct X)$, which is multiplicative,
$$
\symbint(AB) = \symbint(A) \cdot \symbint(B),
$$
and such that there is an exact sequence:
$$
0 \longrightarrow \Psisc^{k-1,0}(X) \longrightarrow \Psisc^{k,0}(X) \longrightarrow 
S^k(\sct X) / S^{k-1}(\sct X) \longrightarrow 0.
$$

There is also a boundary principal symbol map which carries more information than just
the restriction of $\symbint$ to the scattering cotangent bundle over the boundary. 
Since $[\Vsc(X), \Vsc(X)] \subset x \Vsc(X)$, there is also a well-defined notion of
the order $x^0$ part of $P \in \Diffsc^k(X)$ at $\pa X$. The {\it full} symbol of this
order zero part is a (in general nonhomogeneous) polynomial of degree $k$ on
$\sct_{\pa X} X$, denoted
$\symbb(P)$, whose order $k$ part agrees with $\symbint(P)$ restricted
to the boundary. This symbol map extends to $\symbb : \Psisc^{k,0}(X) \to
S^k(\sct_{\pa X} X)$ which is also multiplicative,
$$
\symbb(AB) = \symbb(A) \cdot \symbb(B),
$$
and gives another exact sequence
$$
0 \longrightarrow \Psisc^{k,1}(X) \longrightarrow \Psisc^{k,0}(X) \longrightarrow 
S^k(\sct_{\pa X} X) \longrightarrow 0.
$$

Surprisingly at first sight, perhaps, this symbol map is totally well defined, i.e.,
there is no ambiguity modulo symbols of order $k-1$ (expressed by the exact sequence
above). This is because we already defined the symbol map as only taking values
at $x=0$; the analogous procedure for the interior symbol would be to divide by $|\xi|^k$
and restrict to the sphere at infinity, whereupon (at least for classical operators)
the interior symbol would also become completely well defined.

The combinination of the two principal symbols $\symbint(A)$ and
$\symbb(A)$, called the joint symbol $j_\scl(A)$ in \cite{RBMSpec},
thus determines $A\in\Psisc^{k,0}(X)$ modulo $\Psisc^{k-1,1}(X)$.
In particular, if $A\in\Psisc^{0,0}(X)$ and both $\symbint(A)=0$ and
$\symbb(A)=0$ then $A$ is compact on $L^2(X)$.
In order to talk about both symbols of $A$ at the same time, it is
useful to introduce
$$
\Csc X=\scS X \cup \sct_{\bX} X
$$
where $\scS X$ is the scattering cosphere bundle, i.e.\ the quotient of
$\sct X\setminus 0$
by the natural $\Real^+$-action. Thus, $j_\scl(A)$ is a function
on $\Csc X$; it is given by $\symbb(A)$ on $\sct_{\bX}(X)$ and the
rescaled $\symbint(A)$ on $\scS X$ (we divide out by $|\xi|^k$,
so for $k\neq 0$, $j_\scl(A)$ over $\scS X$
is invariantly defined as a section of a
trivial line bundle).

\

The bundle $\sct X$, restricted to $\pa X$, plays a fundamental role in the 
scattering calculus since scattering microlocal analysis takes place on it. We now briefly
explain these ideas. For notational simplicity let us denote the restriction of
$\sct X$ to the boundary by $K$.
Given local coordinates $(x, y)$ near the boundary
of $\pa X$, there are coordinates $(y, \tau, \mu)$ induced on $K$ as above.

We have already seen that the joint symbol $j_\scl(A)$ of $A \in \Psisc^{k,0}(X)$
is a function on $\Csc X$. We will say that $A$ is elliptic at $q \in \Csc X$ if 
$j_\scl(A)(q) \neq 0$ and define the characteristic variety of $A$, $\Sigma(A)$, to
be the set of non-elliptic points of $\Csc X$. 
This notion of ellipticity allows us to define the
`scattering wavefront set' $\WFsc(u)$ 
of a distribution $u \in \dist(X)$. Namely,
we define it to be the subset of $\Csc X$ whose complement is
\begin{equation}\begin{gathered}
\WFsc(u)^\complement = \{ q \in \Csc X \mid \ \exists \ A \in \Psisc^{0,0}(X) \text{ such
that } \\
A \text{ is elliptic at } q \text{ and } Au \in \dCinf(X) \}.
\end{gathered}
\end{equation}
For example, the scattering wavefront set of $e^{i/x}$ is contained in
$\{ \tau = -1, \mu = 0 \}$. This is easy to see since the function $e^{i/x}$ is
annihilated by $-i x^2 \pa x - 1$ and $-i \pa_{y_i}$ whose boundary symbols are
$-\tau - 1$ and $\mu_i$, respectively.
Note that the restriction of $\WFsc(u)$ to $\scS_{X^o} X$ is just the
usual wave front set $\WF(u)$.

The scattering wavefront set of $u$ is empty if and only if
$u \in \dCinf(X)$; so its part in $K$
may be regarded as a measure of nontrivial behaviour of the distribution at the
boundary, where `trivial' means `all derivatives vanish at the boundary'. 
This is useful in solving equations such as $(\Lap - \l^2) u = 0$. The 
interior symbol of $\Lap - \l^2$ is elliptic, hence for any such $u$,
$\WFsc(u)\subset K$; in particular $u$ must be $\Cinf$
in the interior. On the other hand,
the boundary symbol, $\tau^2 + h(y, \mu) - \l^2$ is not
elliptic when $\l > 0$ and $u$ can then be expected to have nontrivial scattering
wavefront set in $K$. For this reason, we will concentrate on the geometry
and analysis associated $K$. We remark that the corresponding statements on
$\scS_{X^o}X$ are the subject of standard microlocal analysis
(see e.g.\ \cite{FIO1}) and
are very similar in nature; the two are connected by (a localized version
of) the Fourier transform, see \cite{RBMSpec, RBMZw}.

The manifold $K$ carries a natural contact structure. To see this, recall
that $\sct X$ is 
canonically isomorphic to the usual cotangent bundle over the interior, 
so the symplectic
form may be pulled back to $\sct X$, though it blows up as $x^{-2}$ at the boundary
then. Contracting with the vector $x^2 \pa_x$, we obtain a smooth one-form 
whose restriction to $K$ we denote $\chi$.
In local coordinates $(y, \tau, \mu)$, it takes the form
\begin{equation}
\chi \equiv i_{x^2 \pa_x} (\omega) = d\tau + \mu \cdot dy.
\label{chi}\end{equation}
Thus, $\chi$ is a contact form, i.e.\ it is nondegenerate. Changing the boundary defining
function to $\td x = a x$, the form $\td \chi = a \chi$ changes by a nonzero multiple,
so the contact structure on $K$ is totally well-defined. 

Using the contact structure we define the Hamilton vector field of 
$A \in \Psisc^{k,0}(X)$ to be the Hamilton vector field of its boundary
symbol, and the bicharacteristics of $A$ to be the integral curves of its Hamilton
vector field. The bicharacteristics turn out to play a similar role in the scattering
calculus that they do in standard microlocal analysis. For example, Melrose 
\cite{RBMSpec} showed that 

\begin{thm} Suppose $A \in \Psisc^{k,0}(X)$ has real boundary symbol. Then for
$u \in \Cinf(X^o) \cap \dist(X)$,
we have $\WFsc(Au) \subset \WFsc(u)$, $\WFsc(u) \setminus \WFsc(Au)
\subset \Sigma(A)$, and $\WFsc(u) \setminus \WFsc(Au)$ is a union of maximally
extended bicharacteristics of $A$ inside $\Sigma(A) \setminus \WFsc(Au)$.
\end{thm}

\

Next, we discuss Legendrian distributions. A Legendre
submanifold of $K$ is a smooth submanifold $G$ of
dimension $n-1$ such that $\chi$ vanishes on $G$. 
The simplest examples of Legendre submanifolds are given by projectable
Legendrians, that is, those Legendrians $G$ for which the projection to the base,
$\pa X$, is a diffeomorphism, and therefore coordinates $y$ on the base
lift to smooth coordinates on $G$. It is not hard to see, using the explicit form
\eqref{chi} for the contact form, that $G$ is then given by the graph of a differential
$\phi(y)/x$ where $\tau = -\phi(y)$ on $G$. Conversely, any function $\phi$ gives rise
to a projectable Legendrian. We say that $\phi$ is a phase function which
parametrizes $G$. A Legendrian distribution associated to a projectable Legendrian
$G$ is by definition a function of the form 
$$
u = x^m e^{i\phi(y)/x} a(x, y),
$$
where $\phi$ parametrizes $G$ and $a \in \Cinf(X)$. It is quite easy to show that
$\WFsc(u) \subset G$. 

In the general case, where $G$ is not necessarily projectable, we say that a
function $\phi(y, w)/x$, $w \in \RR^k$, is a non-degenerate phase function parametrizing 
$G$, 
locally near $q = (y_0, \tau_0, \mu_0) \in G$, if for some $w_0$
$$
q = (y_0, d_{(x,y)}(\phi(y_0, w_0)/x)) \ \text{ and } \ d_w \phi(y_0, w_0) = 0,
$$
$\phi$ satisfies the non-degeneracy hypothesis 
\begin{equation}
d_{(y,w)}\big( \frac{\pa \phi}{\pa w_i}  \big) \text{ are linearly independent at }
(y_0, w_0),
\end{equation}
and locally near $q$, $G$ is given by
\begin{equation}
G=\{(y,d_{(x,y)}(\phi/x)):\ (y,w)\in C_\phi\} \text{ where } 
C_\phi=\{(y,w):\ d_w\phi(y,w)=0\}.
\label{Legdist}
\end{equation}
Notice
that non-degeneracy implies that $C_\phi$ is a smooth $(n-1)$-dimensional submanifold
of $(y, w)$-space and that the correspondence in \eqref{Legdist} is $\Cinf$. 

A scattering half-density valued
Legendre distribution $u\in\Isc^m(X,G;\Osch)$
associated to $G$ is a distribution in $\dist(X;
\Osch)$ such that 
\begin{equation}
u=u_0+\sum_{j=1}^J u_j\cdot\nu_j,
\end{equation}
where $u_0\in\dCinf(X;\Osch)$, $\nu_j\in\Cinf(X;\Osch)$, and $u_j$ are
given in local coordinates by integrals
\begin{equation}\label{eq:Leg-7}
u_j(x,y)=x^{m-\frac{k}{2}+\frac{n}{4}}
\int_{\Real^k} e^{i\phi_j(y,w)/x}\,a_j(x,y,w)\,dw,
\end{equation}
where $\phi_j$ are non-degenerate phase functions parametrizing $G$,
$a_j\in\Cinf_c([0,\epsilon)_x\times U_y\times U'_w)$. Then $u$ is
smooth in the interior of $X$, and as in the projectable case, 
the scattering wavefront of $w$ is contained
in $G$. Note that the definition of Legendre
distributions is such that for $A\in\Diffsc^k(X)$, $u\in\Isc^m(X,G;\Osch)$
we have $Au\in\Isc^m(X,G;\Osch)$, i.e.\ Legendre distributions are preserved
under the application of scattering differential operators. More
generally, this statement remains true for scattering pseudo-differential
operators $A\in\Psisc^{k,0}(X)$; see \cite{RBMZw}.

\

We also need to discuss two different types of distributions associated to the
union of two Legendre submanifolds. The simpler type is the class of intersecting
Legendrian distributions defined in \cite{Hassell:Plane} and modelled directly on
the intersecting Lagrangian distributions of Melrose and Uhlmann 
\cite{Rbm-Uhl:Intersecting}. A pair
$(G_0, G_1)$ is an intersecting pair of Legendre submanifolds if 
$G_0$ is a smooth Legendre submanifold and $G_1$ is a smooth Legendre submanifold
with boundary, meeting $G_0$ transversally and such that $\pa G_1 = G_0 \cap G_1$. 
A nondegenerate parametrization near $q \in \pa G_1$ is a phase function
$$
\phi(y, w, s)/x = (\phi_0(y, w) + s \phi_1(y, w, s))/x, \quad w \in \RR^k,
$$
such that
$$
q = (y_0, d_{(x,y)} \phi(y_0, w_0, 0)/x), \quad d_{w,s} \phi(y, w, 0) = 0
$$
and such that $\phi_0$ is a nondegenerate parametrization of $G_0$ near $q$,
$\phi$ parametrizes $G_1$ for $s > 0$ with 
$$
ds, \ d \frac{\pa \phi}{\pa v_i} \ \text{ and } \ d \frac{\pa \phi}{\pa s} \ 
\text{ linearly independent  at } q.
$$
A scattering half-density valued
Legendre distribution $u\in\Isc^m(X,(G_0, G_1);\Osch)$
associated to the pair $(G_0,G_1)$ is a distribution in $\dist(X;
\Osch)$ such that 
\begin{equation}
u=u_0+\sum_{j=1}^J u_j\cdot\nu_j,
\end{equation}
where $u_0\in \Isc^{m+1/2}(X,G_0;\Osch)+ I^{m}_{sc,c}(X,G_1;\Osch)$, 
$\nu_j\in\Cinf(X;\Osch)$, and $u_j$ are
given by integrals
\begin{equation}\label{eq:Leg-int}
u_j(x,y)=x^{m-\frac{k+1}{2}+\frac{n}{4}}
\int_{\Real^k} \int_0^\infty e^{i\phi_j(y,w,s)/x}\,a_j(x,y,w,s)\,ds\,dw,
\end{equation}
where the $\phi_j$ are non-degenerate phase functions parametrizing $(G_0,G_1)$, and the
$a_j\in\Cinf_c([0,\epsilon)_x\times U_y\times U'_w)\times [0,\infty)_s$. The subscript
$c$ in $I^{m}_{sc,c}(X,G_1;\Osch)$ above indicates that the microlocal support is
in the interior of $G_1$. The order of $u$ on $G_0$ is always equal to that
on $G_1$ plus one half (to see this, integrate by parts once in $s$). As before, $u$ is
smooth in the interior of $X$, and
the scattering wavefront of $u$ is contained
in $G_0 \cup G_1$.

Following \cite{RBMZw},
we also need to discuss certain pairs of Legendre submanifolds with
conic singularities and distributions corresponding to them.
So suppose first that $\td x$ is any boundary defining function and $\Gs$ is a smooth
Legendre submanifold of $\sct_{\partial X} X$ which is given
by the graph of the differential $\l \,d(1/\td x)$.
Suppose also that $G$
is a Legendre submanifold of $K$ and that its
closure satisfies
\begin{equation}
\cl(G)\setminus G\subset\Gs,
\end{equation}
and $\cl(G)\setminus G$ is the site of an at most conic singularity of
$\cl(G)$. That is, $\cl(G)$ is the image, under the blow-down map,
in $\sct_{\partial X} X$
of a smooth (closed
embedded) manifold with boundary, $\Gh$, of
\begin{equation}
[\sct_{\partial X}X;
\Span\{d(1/\td x)\}]
\end{equation}
which intersects the front face of the blow-up transversally.
This means in local coordinates generated by $(\td x, y)$,
(writing $\mu=r\muh=|\mu|\muh$) that
\begin{equation}
\Gh=\{(y,\tau,r,\muh):\ r\geq 0,\ |\muh|=1,
\ \tau=T(y,r,\muh),\ g_j(y,r,\muh)=0\},
\label{Ghat}
\end{equation}
\begin{equation}
\cl(G)=\{(y,\tau,\mu):\ \tau=T(y,|\mu|, \muh),\ g_j(y,|\mu|, \muh)=0\},
\ \muh=\mu/|\mu|
\end{equation}
where the $n$ functions $g_j$ are such that $d_{(y,\muh)}g_j$, $j=1,\ldots,
n$, are independent at the base point $(y_0,\muh_0,0)$. In particular,
$r=|\mu|$ has a non-vanishing differential on $\Gh$. Then the pair
$\Gt=(G,\Gs)$ is called an intersecting pair of Legendre submanifolds with
conic points.

\begin{rem}
To avoid confusion, let us emphasize at this point that we do not assume, as
was done in \cite{RBMZw}, that $\xt$ is the 
boundary defining function distinguished by a choice of scattering metric. That assumption
was relevant to the analysis of \cite{RBMZw}, but 
not to a description of the geometry of intersecting
pairs of Legendre submanifolds with conic points, which is independent of any choice of
metric structure. In this paper, the boundary defining function which 
describes $G^\sharp$ will {\it not} always be that which is given by the 
scattering metric.
\end{rem}

\begin{rem}\label{altblowup} 
It is equivalent, in specifying the meaning of the conic singularity
of $G$, to require that $\cl G$ is the image of a smooth manifold with boundary
in the slightly different blown up space
$$
[ \sct_{\pa X} X; \l d(1/\xt) ] ,
$$
i.e., we only need to blow up the graph of $\l d(1/\xt)$ (which is precisely $\Gs$), rather
than its linear span. 
\end{rem}

For later use we will write down the parametrization of a conic pair with
respect to a given boundary defining function $x$, not necessarily equal to $\td x$. 
Then $\Gs$ is parametrized by a function $\phi_0(y)$.  
Blowing up the span of $\Gs$, we have coordinates
$y, \tau, r = |\mu - d_y \phi_0 |$ and $\nu = \widehat{ \mu - d_y \phi_0} \in S^{n-2}$. 
We say that the function 
$\phi/x$, with $\phi \in\Cinf(U_y\times[0,s_0)_s\times
U'_u)$ gives a non-degenerate 
parametrization of the pair $\Gt$ near a singular point $q\in
\pa \hat G$ if $U\times[0, s_0)\times U'$ is a neighborhood
of $q' = (y_0,0,w_0)$ in $\partial X\times[0,\infty)\times \Real^{k}$, $\phi$ is
of the form
\begin{equation}
\phi(y,s,w)=\phi_0(y)+s\psi(y,s,w),\quad \psi\in\Cinf(U\times[0,\infty)\times
U'),
\end{equation}
$$
q = (y_0, -\phi(y_0), 0, \widehat{ d_y \psi(q')})
$$
in the coordinates $(y, \tau, r, \nu)$,
the differentials
\begin{equation}
d_{(y,w)}\psi \Mand d_{(y,w)} \frac{\pa \psi}{\partial {w_j}}
\quad j=1,\ldots,k,
\end{equation}
are independent at $(y_0, 0, w_0)$, and locally near $q$,
\begin{gather}
\hat G=\{(y,-\phi, s|d_y \psi|, \widehat{ d_y \psi}) \mid \ (y,s,w)\in C_\phi\},
\quad \text{where} \\
C_\phi=\{(y,s,w):\ \partial_s\phi(y,s,w)=0,\ \partial_w\psi(y,s,w)=0\}.
\end{gather}
It follows that $C_\phi$
is diffeomorphic to the desingularized submanifold $\Gh$.  

A Legendre distribution $u\in\Isc^{m,p}(X,\Gt;\Osch)$
associated to such an intersecting pair $\Gt$ is a distribution
$u\in\dist(X;\Osch)$ such that 
\begin{equation}
u=u_0+u^\sharp+\sum_{j=1}^J u_j\cdot\nu_j,
\end{equation}
where $u_0\in\Isc^m(X,G;\Osch)$, $u^\sharp\in\Isc^p(X,\Gs;\Osch)$,
$\nu_j\in\Cinf(X;\Osch)$ and the $u_j$ are given by integrals
\begin{equation}\label{eq:rev-con-16}
u_j(x,y)=\int\limits_{[0,\infty)\times\Real^k} e^{i\phi_j(y,s,w)/x}
a_j(x,y,x/s,s,w) \left(\frac{x}{s}\right)^{m+\frac{n}{4}-
\frac{k+1}{2}}\,s^{p+\frac{n}{4}-1}\,ds\,dw
\end{equation}
with $a_j\in\Cinf_c([0,\epsilon)_x\times U_y\times[0,\infty)_{x/s}
\times[0,\infty)_s\times U'_w)$ and $\phi_j$ phase functions 
parametrizing $\Gt$. Note that $k+1$ appears in place of $k$ in
the exponent of $x$ (present as $x/s$) since there are $k+1$ parameters,
$u$ and $s$.

\section{Legendre geometry on manifolds with codimension $2$ corners}
\label{sec:corner-geo}
In this section we will describe a class of Legendre submanifolds associated to
manifolds $M$ with codimension 2 corners which are endowed with certain extra structure that
we will discuss in detail. The aim of this section together with the next is to
describe a class of distributions that will include, in particular, the spectral
projections $\sp(\l)$ on $\XXb$. Since there are several cases of interest, we
will describe the situation in reasonable generality, but before doing so let us
consider the setting of primary interest, $M = \XXb$, to see what the ingredients of
the setup should be. 

Recall that $\XXb = [X^2; (\pa X)^2]$ is the blowup of the double space $X^2$ about
the corner. The boundary hypersurfaces are labelled lb, rb and bf according as they
arise from $(\pa X) \times X$, $X \times (\pa X)$ or $(\pa X)^2$. It has two
codimension two corners, $\lb \cap \bfc$ and $\rb \cap \bfc$; note that $\lb \cap \rb
= \emptyset$. Let $\betab:\XXb\to X^2$ be the blow-down map and 
$\pbL$ and $\pbR$ be the left and right stretched projections $\XXb\to X$,
\begin{equation}
\pbL=\pi_L\circ\betab:\XXb\to X,\ \pbR=\pi_R\circ\betab:\XXb\to X.
\end{equation}
If $z$ or $(x, y)$ are local coordinates on $X$, in the interior or near the boundary,
respectively, then we will denote by $z'$ or $(x', y')$ the lift of these coordinates
to $\XXb$ from the left factor (ie, pulled back
by ${\pbL} ^*$) and by $z''$ or $(x'', y'')$ the
lift from the right factor. Write $\sigma = x'/x''$. Then 
$(x'',\sigma,y',y'')$ are coordinates
on $M$ near $\betab^{-1}(p',p'')\setminus\rf\subset\ffb$ where
$\sigma=\frac{x'}{x''}$.
We have a similar result, with the role of $x'$ and $x''$
interchanged, near $\betab^{-1}(p',p'')\setminus\lf$. Away from $\ffb$,
$\betab$ is a diffeomorphism, so near the interior of $\lf$, $\rf$, or in
the interior of $\XXb$, we can use the product coordinates from $X^2$
directly.

We need to identify a Lie Algebra of vector fields, $\Nu$, on $\XXb$ that plays the
role of $\Vsc$ for the scattering calculus. The structure algebra $\Nu$ will
generate a rescaled cotangent bundle which will play the role of $K$ in the previous
section. It is not hard to guess what $\Nu$ should be. Recall from the introduction
that the spectral projection $\sp(\l) = (2\pi)^{-1} P(\l) P(\l)^*$; heuristically,
then, $(2\pi) \sp(\l)(z, z') = \int P(z', y) \ol{ P(z'',y) } dy$ and as a function
of $z$,  Melrose and Zworski showed that $P(\l)$ is a conic Legendrian pair, so
associated with $\Vsc(X)$ acting on the $z$ variable. Thus, we may confidently predict
that the structure algebra is generated by lifting $\Vsc(X)$ to $\XXb$ from the left and
the right:
\begin{equation}\label{eq:str-b-2}
\Nu = 
(\pbL)^*\Vsc(X)\oplus(\pbR)^*\Vsc(X).
\end{equation}

Let us consider the qualitative properties of this Lie Algebra. Near the interior
of bf, using coordinates $(x'', \sigma, y', y'')$, we have the vector fields
\begin{equation}
{x''} ^2 \pa_{x''} - x'' \sigma \pa_\sigma , \quad x'' \pa_{y''_i}, \quad
x'' \sigma \pa_{y'_j}, \quad x'' \sigma^2 \pa_{\sigma}.
\label{vfields}\end{equation}
Away from lb and rb, $\sigma$ is a smooth nonvanishing function, so this
generates the scattering Lie algebra on the non-compact
manifold with boundary $\XXb\setminus(\lf\cup\rf)$, i.e.\ they
generate $\Vsc(\XXb\setminus(\lf\cup\rf))$
as a $\Cinf(\XXb\setminus(\lf\cup\rf))$-module. Near the interior of
lb, using coordinates $(x', y', z'')$, we have the vector fields 
$$
{x'} ^2 \pa_{x'}, \quad x' \pa_{y'_i}, \quad \pa_{z''_j}.
$$
This is the fibred-cusp algebra developed by Mazzeo and Melrose 
\cite{Mazzeo-Melrose:Fibred}: the vector
fields restricted to lb are tangent to the fibers of a fibration,
in this case the projection
from lb $= (\pa X) \times X \to \pa X$, and in addition, for a distinguished
boundary defining function, here $x''$, $V x'' \in {x''} ^2 \Cinf (\XXb)$
if $V\in\Nu$. A similar
property holds for the right boundary rb. Near one of the corners, say $\lb \cap \bfc$,
we have the vector fields \eqref{vfields}. NB: at $\sigma = 0$, the vector
fields ${x''}^2 \pa_{x''}$ and $x'' \sigma \pa_\sigma$ are not separately in
this Lie algebra, only their difference is. This is important, say, in
checking \eqref{eq:stability}. 

The structure algebra $\Nu$ may be described globally in terms of a suitably chosen
total boundary defining function $x$ (that is, a product of boundary defining functions
for each hypersurface). We define $x$ by the Euclidean-like formula
\begin{equation}
x^{-2} = {x'}^{-2} + {x''}^{-2}.
\label{eq:XXb-bd-def-fn}
\end{equation}
Note that here, $x'$, say is an arbitrary defining function for $X$, lifted to
$\XXb$ via the left projection, but then $x''$ is the {\it same } function lifted
via the right projection. We have fibrations $\phi_{\lb}$ and $\phi_{\rb}$
on lb and rb, the projections down to $\pa X$, and the trivial fibration
$\phi_{\bfc} = \id : \bfc \to \bfc$. Then
\begin{equation}
\Nu = \{ V \in \Vb(\XXb) \mid V x \in x^2 \Cinf(\XXb), 
 V \text{ tangent to } \phi_{\lb}, \, \phi_{\rb}, \, \phi_{\bfc} \}.
\end{equation}

\

Now we consider a generalization of the structure just described. 
Let $M$ be a manifold with corners of
codimension $2$, $N=\dim M$, and let 
$M_1(M)$ denotes the set of its boundary hypersurfaces. Suppose
that for each boundary hypersurface $H$ of $M$ we have a fibration
\begin{equation}
\phi_H:H\to Z_H.
\end{equation}
We also impose the following compatibility conditions between the fibrations.
We assume that there is one boundary hypersurface of $M$ (possibly 
disconnected, but which we consider as one, and regard as a single element of
$M_1(M)$) which we denote by $\mf$ (`main face'), such that
\begin{equation}
\phi_{\mf}=\id,
\end{equation}
and that for all other $H,H'\in M_1(M)$, not equal to $\mf$, then either
$H = H'$ or
$H\cap H'=\emptyset$, and if $H \neq \mf \in M_1(M)$ 
then the base of the fibration $Z_H$ is a compact manifold without boundary,
$H\cap \mf
\neq\emptyset$ and the fibers of $\phi_H$ intersect $\mf$ transversally. 
Additionally, we suppose that a distinguished total boundary defining function
$x$ is given. 

In terms of local coordinates the conditions given above
mean that there are local coordinates
near $p\in\mf\cap H$ of the form $(\xt,x'',y',y'')$
(with $y'=(y'_1,\ldots,y'_l)$, $y''=(y''_1,\ldots,y''_m)$) such that
$x''=0$ defines $\mf$, $\xt \equiv  x/x''=0$ defines $H$, and the fibers of $\phi_H$
are locally given by
\begin{equation}
\{(0,x'',y',y''):\ y'=\text{const}\}.
\end{equation}
Near $p\in\interior(\mf)$ we simply have coordinates of the form
$(x,\yb)$, $\yb=(\yb_1,\ldots,\yb_{N-1})$ (the fibration being trivial there),
while near $p\in\interior(H)$, $H\in M_1(X)$, $H\neq\mf$,
we have coordinates $(x,y',z'')$ (with $y'=(y'_1,\ldots,y'_l)$,
$z''=(z''_1,\ldots,z''_{m+1})$) and the fibers of $\phi_H$ are
locally given by
\begin{equation}
\{(0,y',z''):\ y'=\text{const}\}.
\end{equation}

We let $\Phi$ be the collection of the boundary fibration-maps:
\begin{equation}
\Phi=\{\phi_H:\ H\in M_1(M)\}
\end{equation}
which we think of as a `stratification' of the boundary of $M$.
It is also convenient to introduce the non-compact manifold with boundary
$\Mt$ given by removing all boundary hypersurfaces but $\mf$ from $M$:
\begin{equation}
\Mt=M\setminus\bigcup\{H\in M_1(M):\ H\neq\mf\}.
\end{equation}

\begin{Def} The fibred-scattering structure on $M$, with respect to the 
stratification $\Phi$ and total boundary defining function $x$, is
the Lie algebra $\Vsf(M)$ of smooth vector fields satisfying
\begin{equation}\label{eq:Vsf-def}
V\in\Vsf(M)\Miff V\in\Vb(M),\ Vx\in x^2\Cinf(M),
\ V\ \text{is tangent to the fibres of}
\ \Phi.
\end{equation}
\end{Def}
It is straightforward to check that $\Vsf(M)$ is actually well-defined
if $x$ is defined up to multiplication by positive functions which
are constant on the fibers of $\Phi$.
Note that \eqref{eq:Vsf-def} gives
local conditions and near the interior of $\mf$, 
hence in $\Mt$, where
$\Phi$ is trivial, they are equivalent to requiring that $V\in\Vsc(\Mt)$.
Moreover, we remark that
in the interior of each boundary hypersurface
$H$, the structure is just a fibred cusp structure (though on
a non-compact manifold) which is described in \cite{Mazzeo-Melrose:Fibred}.

We write $\Diffsf(M)$ for the differential operator algebra generated
by $\Vsf(M)$. One can check in local coordinates (as was done explicitly for
$\XXb$ above) that locally these are precisely the $\Cinf(M)$-combinations of
$N$ independent vector fields, and therefore such vector fields are 
the space of all smooth sections of a vector
bundle, $\Tsf M$. We denote the dual bundle by $\sfT M$. 
It is spanned by one-forms of the form $d(f/x)$ where $f\in\Cinf(M)$ is
constant on the fibers of $\Phi$ (which is of course no restriction
in the interior of $\mf$). Since we are more interested in $\sfT M$, we will write
down explicit basis in local coordinates for $\sfT M$; the corresponding bases for
$\Tsf M$ are easily found by duality. Near $p\in\mf\cap H$, $H\in M_1(M)$,
we choose the local coordinates
$(x'',\xt,y',y'')$ as discussed above. 
Then a basis of $\sfT M$ is given by
\begin{equation}\label{eq:fib-b-2}
\frac{dx}{x^2},\ \frac{d\xt}{x},\ \frac{dy'_j}{x},\ \frac{dy''_k}{x''}
\ (j=1,\ldots,l,\ k=1,\ldots,m).
\end{equation}
Alternatively, a basis is given by
\begin{equation}
d\left(\frac1{x}\right),\ d\left(\frac1{x''}\right),
\ \frac{dy'_j}{x},\ \frac{dy''_k}{x''}
\ (j=1,\ldots,l,\ k=1,\ldots,m).
\end{equation}
Corresponding to these bases, we introduce coordinates
on $\sfT M$. Thus for $p\in\ffb$ we write $q\in\sfT_p M$
as
\begin{equation}\label{eq:corner-9-a}
q =\tau\frac{dx}{x^2}+\eta\frac{d\xt}{x}+\mu'\cdot\frac{dy'}{x}
+\mu''\cdot\frac{dy''}{x''},
\end{equation}
giving coordinates $(x'', \xt, y'_j, y''_k; \tau, \eta, \mu'_j, \mu''_k)$ on
$\sfT M$, or alternatively, 
\begin{equation}\label{eq:corner-9-b}
q =\td\tau\frac{dx}{x^2}+ \tau'' \frac{d x''}{(x'')^2}+\mu'\cdot\frac{dy'}{x}
+\mu''\cdot\frac{dy''}{x''},
\end{equation}
giving coordinates $(x'', \xt, y'_j, y''_k; \td \tau, \tau'', \mu'_j, \mu''_k)$. The
relation between the two is given by $\tau = \td \tau + \xt \tau''$, $\eta = -\tau''$.
In the interior of $H$ where we have coordinates $(x,y',z'')$, $z''$
being coordinates on an appropriate coordinate patch in the interior of $X$
(considered as the right factor in $X^2$), a basis is given by
\begin{equation}
\frac{dx}{x^2},\ \frac{dy'_j}{x}\ (j=1,\ldots,l),\ dz''_k
\ (k=1,\ldots,m+1).
\end{equation}

If $M=\XXb$ with $x$ given by \eqref{eq:XXb-bd-def-fn} and
the boundary fibration by $\Phi = (\phi_{\lb}, \phi_{\rb}, \phi_{\bfc})$, 
we can replace
$x$ in all these equations by $x'$ since near $\lf$ we
have $x=x'(1+\sigma^2)^{-1/2}$ where $\sigma$ can be used as
the coordinate $\xt$ in our general notation. Therefore,
\begin{equation}\label{eq:str-b-7}
\sfT \XXb=(\pbL)^*\sct X\oplus(\pbR)^*\sct X
\end{equation}
and correspondingly for the structure bundle $\Tsf \XXb$
we have \eqref{eq:str-b-2} as desired.

\

We have already discussed how the b-double space $\XXb$ fits into this
framework. Another example
is given by the following setup. Let $Y$ be a manifold with boundary and let
$C_j$, $j=1,\ldots,s$, be disjoint closed embedded submanifolds of
$\partial Y$. Let $x$ be a boundary defining function of $Y$.
Let $M=[Y;\cup_j C_j]$ be the blow up of $Y$ along the $C_j$. The pull-back of
$x$ to $M$ by the blow-down map $\beta:M\to Y$ gives a total boundary
defining function, and $\beta$ restricted to the boundary hypersurfaces of
$M$ gives the fibrations. Thus, the lift of $\partial Y$ to $M$ is the main
face, $\mf$, with the trivial fibration, and if $H_j=\beta^{-1}(C_j)$ is
the front face corresponding to the blow-up of $C_j$, then $\phi_{H_j}:
H_j\to C_j$. These data define a stratification on $\pa M$ as
discussed above.
Such a setting is a geometric
generalization of the three-body problem and it has been discussed in
\cite{Vasy:Propagation-2}, see also \cite{Vasy:Propagation},
where the structure bundle was denoted by $\TSc M$. We note that
$\sfT M$ can be
identified as $\beta^*\sct Y$; this is particularly easy to see by
considering that $\sfT M$ is spanned by $d(f/x)$ with $f$
constant on the fibers of $\beta$. We remark that
any other boundary defining function of $Y$ is of the form $ax$ where
$a\in\Cinf(Y)$, $a>0$, so $\beta^*a\in\Cinf(M)$ is constant on the
fibers of $\phi_H$, i.e.\ on those of $\beta$, so this structure does not
depend on the choice of the boundary defining function of $Y$ as long as
we use its pull-back as total bondary defining function on $M$.

\

Returning to the general setting, for each $H\in M_1(M)$, $H\neq\mf$,
there is a natural subbundle of $\sfT_{H} M$ spanned by
one-forms of the form $d(f/x)$ where $f$ {\em vanishes} on $H$.
Restricted to each fiber $F$ of $H$, this is isomorphic to the
scattering
cotangent bundle of that fiber, and we denote this subbundle by
$\sct (F;H)$. Indeed, if $f$ vanishes on $H$ then $f=\xt\ft$,
$\ft\in\Cinf(M)$, so modulo terms vanishing at $H$, $d(f/x)=d(\ft/x'')$ is
equal to $-\ft\,dx''/(x'')^2 + d_{y''}\ft /x''$, so locally
a basis of $\sct (F;H)$ is given by $dx''/(x'')^2$,
$dy''_j / x''$, $j=1,\ldots,m-1$.
The quotient bundle, $\sfT_{H}M/\sct (F;H)$, is
the pull back of a bundle $\sfN Z_H$ from the base $Z_H$ of the
fibration $\phi_H$, and $\sfN Z_H$ is isomorphic to $\sct_{Z_H\times\{0\}}
((Z_H)_{y'}\times[0,\epsilon)_x)$. If $M=\XXb$, $\sfN Z_{\lf}=\sfN\bX$
is simply $\sct_{\bX} X$, and similarly for $\rf$.
The induced map $\sfT_{H} M\to \sfN Z_H$ will be denoted by
$\phit_H$; its fibres are isomorphic to $\sct F$. In terms of the
coordinates $(x'',y',y'',\tau,\eta,\mu',\mu'')$ on $\sfT_{H}M$,
the fibers are given by $y'$, $\tau$, $\mu'$ being constant.

It will be very often convenient to realize $\phi_H^*(\sfN Z_H)$
as a subbundle of
$\sfT_{H} M$ by choosing a subbundle complementary to $\sct (F;H)$.
Such a subbundle arises naturally if we have a smooth fiber metric on
$\sfT_{H}M$ which allows us to choose the orthocomplement of $\sct
(F;H)$. It can also be realized very simply if the total space of the
fibration ($H$ in our case) has a product structure. In particular,
for $M=\XXb$,
the bundle decomposition \eqref{eq:str-b-2} gives $\sct (F;\lf)
=(\pbR)^*_{\lf}\sct X$ and $\phi_{\lf}^*\sfN\bX=(\pbL)^*_{\lf}\sct X$.

We wish
to define Legendre distributions and more generally Legendre distributions
with conic points, associated to our structure. In the interior of the main face,
we have a scattering structure, so we will consider Legendre submanifolds of 
$\sfT_{\mf}M$ satisfying certain conditions at the boundary of $\mf$. Before we
do so, we need to consider the behaviour of the contact form on $\sfT_{\mf}M$ at
the boundary. Using the coordinates \eqref{eq:corner-9-a}, we find that the
contact form is
$$
\chi = i_{x^2 \pa_x} (\omega) = d\tau+\eta\, d\xt+\mu'\cdot dy'+\xt\mu''\cdot dy''.
$$
This degenerates at the boundary. In fact it is linear on fibres and
vanishes identically on $\sct (F;H)$
at the corner, and therefore induces a one-form on the quotient bundle 
$\sfT_{H}M/\sct (F;H)$ over the corner. 
In local coordinates, this is $d\tau + \mu' \cdot dy'$, and is the pullback
of a non-degenerate one-form 
on $\sfN Z_H$, thus giving $\sfN Z_H$ a natural
contact structure.

The choice of local coordinates $(x'',\xt,y',y'')$ and the corresponding
basis of $\sfT M$ given in \eqref{eq:corner-9-b} induces the one-form
\begin{equation}
d\tau'' +\mu''\cdot {dy''}
\end{equation}
on the fibers of the restriction of $\phit_H$ to
$\sfT_{H\cap\mf}M$.
A simple computation shows that under a smooth change of coordinates
$$
{\xt} \to a \xt, \quad {x''} \to b x'', \quad {y'} \to Y'(y') + \xt Y'_1, 
\quad y'' \to Y'',
$$
such that ${x} = {\xt} {x''}$ changes according to ${x} \to x (g(y') + \xt g_1)$,
this is
a well-defined one-form on each fibre of $\phit_H$
up to multiplication by the smooth positive factor $a$.
Thus, there is a natural contact structure on
the fibres of $\phit_H$ as well. All three
contact structures play a role in the following definition.

\begin{Def}\label{M-Leg-submfld}
A  Legendre submanifold $G$ of $\sfT M$ is a Legendre
submanifold of $\sfT_{\mf}M$ which is transversal to $\sfT_{H\cap
\mf}M$ for each $H\in M_1(M)\setminus\{\mf\}$,
for which the map $\phit_H$ induces a fibration
$\phit'_H:G\cap\sfT_{H\cap\mf}M\to G_1$ where $G_1$ is a Legendre
submanifold of $\sfN Z_H$ whose fibers are Legendre submanifolds
of $\sct_{\partial F} F$.
\end{Def}

We will say that a function $\phi(\td x, y', y'', v, w)/x$, $v \in \RR^k$, 
$w \in \RR^l$, is a non-degenerate parametrization of $G$ locally
near $q\in\sfT_{\mf\cap H}M\cap G$
which is given in local coordinates as $(\xt_0,y'_0,y''_0,\tau_0,\eta_0,
\mu'_0,\mu''_0)$ with $\xt_0=0$, if $\phi$ has the form 
\begin{equation}\label{eq:fib-ph-9}
\phi(\xt,y',y'',v, w)=\phi_1(y',v)+\xt\phi_2(\xt,y',y'',v, w)
\end{equation}
such that $\phi_1$, $\phi_2$ are smooth on neighborhoods of $(y'_0,v_0)$ 
(in $\bX\times\Real^k$) and
$q'=(0,y'_0,y''_0,v_0,w_0)$ (in $\ffb\times\Real^{k+k'}$)
respectively with
\begin{equation}
(0, y'_0, y''_0, d_{(x,\xt,y',y'')}(\phi/x)(q'))=q,\ d_v\phi(q')=0,
\ d_w\phi_2(q')=0,
\end{equation}
$\phi$ is non-degenerate in the sense that
\begin{equation}\label{eq:corner-20}
d_{(y',v)}\frac{\pa \phi_1}{\partial {v_j}},\ d_{(y'',w)} \frac{\pa \phi_2}
{\partial {w_{j'}}}
\ j=1,\ldots,
k,\ j'=1,\ldots,k'
\end{equation}
are independent at $(y'_0,v_0)$ and $q'$ respectively,
and locally near $q$, $G$ is given by
\begin{equation}
G=\{(\xt,y',y'',d_{(x,\xt,y',y'')}(\phi/x)):\ (\xt,y',y'',v,w)\in
C_\phi\}
\end{equation}
where
\begin{equation}\label{eq:corner-22}
C_\phi=\{(\xt,y',y'',v,w):\ d_v\phi=0,\ d_w\phi_2=0\}.
\end{equation}
Note that
these assumptions indeed imply that in $\sfT_{\interior(\mf)}M$, but in the
region where $C_\phi$ parametrizes $G$, $\phi/x$ is a non-degenerate
phase function (in the usual sense) for the Legendre submanifold $G$.
Also, $\phi_1$ parametrizes $G_1$ and for fixed $(y', v) \in C_{\phi_1}$,
$\phi_2(0, y', \cdot, v', \cdot)$ 
parametrizes the corresponding fiber of $\phit'_H$.

\begin{prop} Let $G$ be a Legendre submanifold. Then for any point
$q \in \pa G$ there is a non-degenerate parametrization of $G$ in some neighbourhood
of $q$. 
\end{prop}

\begin{proof}
The existence of such a phase
function can be established using \cite[Proposition~5]{RBMZw} and the
transversality condition ($G$ to $\sfT_{H\cap\mf}M$). Indeed, this
proposition shows that there exists a splitting $(y'_\flat,y'_\sharp)$
and a corresponding splitting $(\mu'_\flat,\mu'_\sharp)$ such that
$(y'_\sharp,\mu'_\flat)$ give coordinates on $G_1$. It also shows that
there exists a splitting $(y''_\flat,y''_\sharp)$ such that
$(y''_\sharp,\mu''_\flat)$ give coordinates on the fibers of
$G\cap\sfT_{\lf\cap\ffb}M\to G_1$. Taking into account the transversality
condition, we deduce that $(\xt,y'_\sharp,y''_\sharp,\mu'_\flat,
\mu''_\flat)$ give coordinates on $G$. We write
\begin{equation}\begin{split}
&\tau=T_1(y'_\sharp,\mu'_\flat)+
\xt\,T_2(\xt,y'_\sharp,y''_\sharp,\mu'_\flat,\mu''_\flat),\\
&y'_\flat=Y'_{\flat,1}(y'_\sharp,\mu'_\flat)+
\xt\,
Y'_{\flat,2}(\xt,y'_\sharp,y''_\sharp,\mu'_\flat,\mu''_\flat),\\
&y''_\flat=Y''_\flat(\xt,y'_\sharp,y''_\sharp,\mu'_\flat,\mu''_\flat).
\end{split}\end{equation}
Here the special form of $\tau$ and $y'_\flat$ at $\xt=0$ follows
from the fibration requirement. We also write $T=T_1+\xt\,T_2$,
and similarly for $Y'_\flat$. Then the argument
of \cite[Proposition~5]{RBMZw} shows that
\begin{equation}\label{eq:fib-28}
\phi/x=(y'_\flat\cdot\mu'_\flat
-T-Y'_\flat\cdot\mu'_\flat
+\xt\,(y''_\flat\cdot\mu''_\flat-Y''_\flat\cdot\mu''_\flat))/x
\end{equation}
gives a desired parametrization. Indeed, the proof of that proposition
guarantees that $\phi/x$ parametrizes $G$ in the interior of $\sfT_{\mf}M$
but near $q$, and the fact that $G$ is transversal to $\sfT_{H\cap\mf}M$,
hence it is the closure of its intersection with $\sfT_{\interior(\mf)}M$
proves that $\phi/x$ gives a desired parametrization of the 
Legendre submanifold $G$ of $\sfT M$. In particular,
\begin{equation}
\phi_1(y',\mu'_\flat)=y'_\flat\cdot\mu'_\flat-T_1(y'_\sharp,\mu'_\flat)
-Y'_{\flat_1}(y'_\sharp,\mu'_\flat)\cdot\mu'_\flat.
\end{equation}
\end{proof}

Next, we define pairs of  Legendre submanifolds with
conic points. Again it will be an intersecting pair with conic points, in the
sense of Melrose-Zworski, in the interior of mf 
with certain specified behaviour at the boundary of mf. 

\begin{Def}\label{M-leg-submfld-conicpts}
A  Legendre pair with conic points, $(G,\Gs)$,
in $\sfT M$ consists of two  Legendre submanifolds $G$ and $\Gs$
of $\sfT M$ which form an intersecting pair with conic points in $\sct_{\td M} \td M$ 
such that for each $H\in M_1(M)\setminus\{\mf\}$ the
fibrations $\phit'_H$, $\phit^{\sharp}_H$ have the
same Legendre submanifold $G_1$ of $\sfN Z_H$ as base and for
which the fibres are intersecting pairs of 
Legendre submanifolds with conic points
of $\sct_{\partial F} F$. 
\end{Def}

Let us unravel some of the complexities in this definition. First of all, $\Gs$ is
by definition a global section, and consequently, $G_1$ is too. Thus, $\Gs$ is
parametrized by a phase function of the form
$(\phi_1(y') + \xt \phi_2(\xt, y', y''))/x$ near $\xt = 0$. Then,
$G_1$ is parametrized by $\phi_1(y')/x$. Second, the fibres of $\phit^{\sharp}_H$
are pairs with conic points $(G(y'), \Gs(y'))$, which we may parametrize by $y'$
since $G_1$ is projectable. It is not hard to see then that
$\Gs(y')$ is parametrized by $\phi_2(0, y', y'')/x''$. 

As before, we can simplify the coordinate expressions by changing to a total boundary
defining function $\overline x$ so that $\Gs$ is parametrized by $\l/\overline{x}$.
Let us temporarily
redefine $x$ to be $\overline{x}$. Then $G_1$ is parametrized by $\l/x$, 
and $\Gs(y')$ by zero. Since $(G(y'), \Gs(y'))$
form a conic pair in the fibre labelled by $y'$, then $\cl G(y')$ is the image under
blowdown of a smooth submanifold $\Gh(y')$, when $\{ \tau'' = 0, \mu'' = 0 \}$
is blown up in the fibre (cf. remark~\ref{altblowup}). Since $|d\mu''| \neq 0$ on
the boundary of $\Gh(y')$, 
it is also the case that each $\cl G(y')$ is the image
under blowdown of a smooth submanifold when we blow up
\begin{equation}
[ \sfT_{\mf}(M); \{ \mu'' = 0, \tau'' = 0, \mu' = 0 \} ],
\label{sPhi-blowup}
\end{equation}
i.e., the span of $d(1/x)$ over mf. In the interior of mf, this is precisely the
blowup that desingularizes $G$, so we see that the conic geometry of $(G, \Gs)$
inside $\sfT_{\mf} M$ holds uniformly at the boundary of mf.

Let us return to our original total boundary defining function $x$. Then $\Gs$ is
parametrized by a function on mf of the form $\phi_0/x = (\phi_1(y') + 
\xt \phi_2(\xt, y', y''))/x$ near the boundary with $H$. To desingularize $G$, we
blow up the span of $\Gs$ inside $\sfT_{\mf} M$. Coordinates on
the blown-up space are 
$$(\xt, y', y'', \tau, \frac{(\tau'' + \phi_2)}{r}, \frac{(\mu' - d_{y'}
\phi_0)}{r}, r = |\mu''|, \nu = \hat \mu'')
$$ 
near the boundary of $\Gh$, the
lift of $G$ to this blown-up
space which is a smooth manifold with boundary.
By a non-degenerate parametrization of $(G, \Gs)$, as above, near a point 
$q\in \partial\Gh \cap \sfT_{\mf \cap H}M$
is meant a function $\phi/x$, where
$$
\phi(\td x, y', y'', s, w) = \phi_1(y') + \xt \phi_2(\xt, y', y'') + 
s\xt\psi(\xt, y', y'', s, w), \quad w \in \RR^k,
$$
is a smooth function 
defined on a neighborhood of $q'=(0,y'_0,y''_0,0,w_0)$ 
with
\begin{equation}
q = (0, y'_0, y''_0, -\phi_1, -\frac{\psi}{|d_{y''}\psi|}, 
\frac{d_{y'}\psi}{|d_{y''}\psi|}, 
0, \widehat{d_{y''} \psi})(q'), 
\quad d_s\phi(q')=0, \quad 
\ d_w\psi(q')=0,
\end{equation}
$\phi$ non-degenerate in the sense that
\begin{equation}
d_{(y'',w)}\frac{\pa \psi}{\partial {w_j}},\ d_{(y'',w)} \psi, 
\quad j=1,\ldots, k,\ 
\end{equation}
are independent at $q'$
and locally near $q$, $\Gh$ is given by
\begin{gather}
\Gh=\{(\xt,y',y'',-\phi_1, -\frac{\psi}{|d_{y''}\psi|}, \frac{d_{y'}\psi}{|d_{y''}\psi|}, 
s|d_{y''}\psi|, 
\widehat{d_{y''} \psi}) \mid \  (\xt,y',y'',s,w)\in
C_\phi\}, \\
C_\phi=\{(\xt,y',y'',s,w):\ d_s\phi \equiv \psi + s\psi_s = 0,\ d_w\psi=0, \ s \geq 0\}.
\label{eq:corner-22p}\end{gather}
Then $\Gh$ is locally diffeomorphic to $C_\phi$. 

\begin{prop} Let $(G, \Gs)$ be as above. Then there is a
non-degenerate para\discretionary{met-}{riz}{metriz}ation of $(G, \Gs)$ near
any point $q \in \partial\Gh \cap \sfT_{\mf \cap H}M$ .
\end{prop}

\begin{proof} It is easiest to parametrize the pair using the total boundary defining
function $x$ for which $\Gs$ is parametrized by $\l/x$, as above, and then change
back to the original coordinates. 
Using the parametrization argument of \cite[Proposition~6]{RBMZw} for
the conic pairs of Legendre submanifolds $(G(y'_0), \Gs(y'_0))$, we see
that there are coordinates $y''$ and
a splitting $(y''_\sharp,y''_\flat)$, which
we can take to be $y''_\sharp=(y''_1,
\ldots,y''_k)$, $y''_\flat=(y''_{k+1},\ldots,y''_m)$, of these
such that after the span of $d(1/x)$ is blown up inside
$\sfT_{\mf} M$ the image of the base point
$q$ is $y''=0$,
$\tau''/|\mu''_m| = 0$, $\muh''_1=\ldots
=\muh''_{m-1}=0$, $\muh''_m=1$ where $\muh''_j=\mu''_j/\mu''_m$. Then 
$(y''_\sharp,\muh''_{k+1},\ldots,\muh''_{m-1},\mu''_m)$ are
coordinates on $\Gh(y'_0)$, and therefore, on $\Gh(y')$ for $y'$ in a neighbourhood of
$y'_0$. Using the transversality, we
conclude that $(\xt,y',y''_\sharp,\muh''_{k+1},\ldots,\muh''_{m-1},
\mu''_m)$ give coordinates on $\Gh$. Then, following
\cite[Equation~(7.5)]{RBMZw},
\begin{equation}
\phi/x=(-T+\xt(s y''_\flat\cdot w-sY''_\flat(\xt,y',y''_\sharp,s,w)\cdot w+
sy''_{m}-sY''_{m}(\xt,y',y''_\sharp,s,w)))/x,
\end{equation}
where on $\Gh$
\begin{equation}\begin{split}
&\tau=T(\xt,y',y''_\sharp,\mu''_m,\muh''_{k+1},\ldots,\muh''_{m-1})\\
&\qquad\qquad
=\lambda+\mu''_m
T_2(\xt,y',y''_\sharp,\mu''_m,\muh''_{k+1},\ldots,\muh''_{m-1}),\\
&y''_\flat=Y''_\flat(\xt,y',
y''_\sharp,\mu''_m,\muh''_{k+1},\ldots,\muh''_{m-1})
\end{split}\end{equation}
(so $s=\mu''_m$, $w=(\muh''_{k+1},\ldots,\muh''_{m-1})$),
gives a parametrization in our sense. Notice that, taking into account
$sw=(\mu_{k+1},\ldots,\mu_{m-1})$, this is of the same form as
\eqref{eq:fib-28} for $s>0$ (with no $y'_\flat$ variables).
\end{proof}

\section{Legendre distributions on manifolds with codimension 2 corners}
\label{sec:corner-Leg-dist}

We now define a class of distributions on a manifold $M$ with stratification
$\Phi$ as 
discussed above. These will be modelled
on oscillatory functions whose phase is compatible with $\Phi$. 
Thus, an example of a Legendre distribution is
an oscillatory function
\begin{equation}\label{eq:fib-osc}
u=x^q e^{i\phi/x}a,\ a,\phi\in\Cinf(M),\ \phi\ \text{constant on the fibers
of}\ \Phi.
\end{equation}
As mentioned in the previous section, the differentials $d(\phi/x)$ for
such functions $\phi$ span $\sfT M$, so the definition of $\sfT M$ is
precisely designed so that 
\begin{multline}\label{eq:stability}
P\in\Diffsf(M), \ u = x^q e^{i\phi/x} a \text{ with } \phi, \, a \text{ as
  in } (\ref{eq:fib-osc}) \\
\implies \ Pu = x^q e^{i\phi/x} \tilde a, \quad \tilde a \in C^\infty (M).
\end{multline}

In general, we need to consider superpositions of such oscillatory functions.
Thus, near the corner $H\cap\mf$, $H\in M_1(M)\setminus\{\mf\}$,
we consider distributions of the form
\begin{equation}\begin{split}\label{eq:dist-1}
u(x'',\xt,y',y'')
=\int e^{i\phi(\xt,y',y'',v,w)/x} &a(x'',\xt,y',y'',v,w)\\
&(x'')^{m-(k+k')/2+N/4}\xt^{r-k/2+N/4-f/2}\,dv\,dw
\end{split}\end{equation}
with $N=\dim M$, $a\in\Cinf_c([0,\epsilon)\times
U\times\Real^{k+k'})$, $U$ open in $\mf$, $f$ the dimension of the fibres of $H$ and
$\phi$ a phase function parametrizing a Legendrian $G$ on
$U$.
On the other hand, in the interior of $H$ we consider distributions
of the form
\begin{equation}\label{eq:dist-2}
u'(x,y',z'')
=\int e^{i\phi_1(y',w)/x} a(x,y',z'',w)
x^{r-k/2+N/4-f/2}\,dw
\end{equation}
with $N=\dim M$, $a\in\Cinf_c([0,\epsilon)\times
U\times\Real^{k})$, $U$ open in $H$,
$\phi_1$ is a phase function parametrizing the Legendrian $G_1$
(associated to $G$ as in Definition~\ref{M-Leg-submfld}) on
$\phi_H(U)$. Note that distributions of the form \eqref{eq:dist-1}
are automatically of the form \eqref{eq:dist-2} in the interior of $H$
since $\phi = \phi_1 + \xt \phi_2$ as in \eqref{eq:fib-ph-9}, so $e^{i\phi/x}$
is equal to $e^{i\phi_1/x}$ times a smooth function in this region. 

We let
$\Osfh M$ be the half-density bundle induced by
$\sfT M$; near the interior of $\mf$ this is just $\Osch \Mt$. For $M=\XXb$,
\eqref{eq:str-b-2} implies that
\begin{equation}
\Osfh \XXb=(\pbL)^*\Osch\otimes(\pbR)^*\Osch,
\end{equation}
so its smooth sections are spanned by $(\pbL)^*|dg|^\half
\otimes(\pbR)^*|dg|^\half$ over
$\Cinf(\XXb)$. Notice that in the blow-up situation $[Y;\cup_j C_j]$ of
the previous section, we have $\Osfh[Y;\cup_j C_j]=\beta^*\Osch Y$.

To simplify the notation and allow various orders
at the various hypersurfaces $H$, we introduce an `order family',
$\calK$, which assigns a real number to each hypersurface. Thus,
$\calK$ is a function
\begin{equation}\label{eq:ord-fam}
\calK:M_1(M)\to\Real.
\end{equation}
Here we only consider one-step polyhomogeneous (`classical') symbols $a$,
but the discussion could be easily generalized to arbitrary
polyhomogeneous symbols. 

\begin{Def}\label{Def:fib-dist}
A Legendrian distribution $u\in\Isf^{\calK}(M,G;\Osfh)$
associated to a 
Legendre submanifold $G$ of $\sfT M$ is a distribution $u\in
\dist(M;\Osfh)$ with the property
that for any $H\in M_1(M)\setminus\{\mf\}$ and for any
$\psi\in\Cinf_c(M)$ supported away from $\bigcup\{H'\in M_1(M):
\ H'\neq H, H'\neq\mf\}$,
$\psi u$ is of the form
\begin{equation}
\psi u=u_0+u_1+\sum_{j=1}^J w_j\cdot\nu_j+\sum_{j=1}^{J'} w'_j\cdot\nu'_j,
\end{equation}
where $u_0\in\dCinf(M;\Osfh)$, $u_1\in\ I_{\scl,\text{c}}^m(\Mt,G;\Osfh)$
($\text{c}$ meaning compact support in $\Mt$), $m=\calK(\mf)$,
$\nu_j,\nu'_j\in\Cinf(M;\Osfh)$, and $w_j,
w'_j\in
\dist(M)$ are given by oscillatory integrals as in \eqref{eq:dist-1}
and \eqref{eq:dist-2} respectively with $m=\calK(\mf)$, $r=\calK(H)$.
\end{Def}

\begin{rem}
Sometimes, when convenient, we use the notation $\Isf^{m,r}(M,G,\Osfh)$
when we localize a fibred Legendre distribution near a hypersurface $H$
so that only the order $m$ at $\mf$ and $r$ at $H$ are relevant.
\end{rem}

\begin{rem}
The reason for the appearance of the fiber-dimension $f$ in our
order convention, i.e.\ in \eqref{eq:dist-1} -- \eqref{eq:dist-2},
can be understood from the blow-up example $[Y;C]$ of the
previous section. Namely, suppose that $C$ is a closed embedded
submanifold of $\partial Y$, $Y$ being a manifold with boundary.
Let $G$ be 
the scattering conormal bundle $\sci N^* C$ of
$C$; this is simply the annihilator of $\Vsc(Y;C)=x\Vb(Y;C)$, $\Vb(Y;C)$
denoting the set of vector fields in $\Vb(Y)$ tangent to $C$ in the ordinary
sense.
If $C$ is given by
$x=0$, $y'=0$, in local coordinates $(x,y)$, $y=(y',y'')$,
and the dual coordinates
are $(\tau,\mu',\mu'')$, then $G$ is given by $x=0$, $y'=0$, $\tau=0$,
$\mu''=0$.
Suppose that
$u\in\Isc^m(Y,G,\Osch)$.
Thus, $u$ is a sum of oscillatory integrals of the form
$$
x^{m-\frac{k}{2}+\frac{n}{4}}\int e^{iy'\cdot w/x}a(x,y,w)\,dw,
$$
$n=\dim Y$, $k$ the codimension of $C$ in $\partial
Y$, and $a$ smooth and compactly
supported.
Pulling back
$u$ to $[Y;C]$ by the blow-down map $\beta$ gives a function of the form
$$
\beta^*u=x^{m-\frac{k}{2}+\frac{n}{4}}b,\quad b\in\Cinf([Y;C]),
$$
that vanishes to infinite order on $\mf$. Indeed $Y'=y'/x$
is a coordinate in the interior of $\ff$, so above we simply took
the inverse Fourier
transform of $a$ in $w$. Thus, taking $G'$ to be the zero section
of $\sfT_{\mf} [Y;C]$, we see that $\beta^* u\in \Isf([Y;C],G',\Osfh)$
with order $\infty$ on $\mf$ and order $m$ on $\ff$. The agreement
of the Legendrian order in these two ways of looking at $u$ justifies
the order convention in the definition.
\end{rem}

For this class of distributions to be useful, we need to know that it
is independent of the choice of parametrizations used. This follows from

\begin{prop} Let $L$ be a Legendrian as in Definition~\ref{M-Leg-submfld},
let $q \in L$, and let $\phi$, $\phit$ be two parametrizations of a
neighbourhood $U \subset L$ containing $q$. Then if $A$ is smooth and
has support sufficiently close to $q$, the function
$$
u = x^q \xt^r
\int e^{i\phi(y',y'',\xt, v,w)/x} a(y',y'',\xt, x'',v,w) \, dv \, dw
$$
can be written
$$
u = u_0 + x^{\tilde q} \xt^{\tilde r}
\int e^{i\tilde\phi(y',y'',\xt, \tilde v,\tilde w)/x} \tilde
a(y',y'',\xt, x'',\tilde v,\tilde w) \, d\tilde v \, d\tilde w 
$$
for some $u_0 \in \dot \Cinf(X)$, $\tilde q$, $\tilde r$, and 
smooth $\tilde a$.
\end{prop}

\begin{proof} The heart of the proof is an adaptation of H\"ormander's
proof of equivalence of phase functions (\cite{FIO1}, Theorem 3.1.6)
to this setting. If $q$ lies above an interior point of mf, then the
result is proved in \cite{RBMZw}, so we only need to consider $q$ lying above a
point on the boundary of mf. Using coordinates as in \eqref{eq:fib-ph-9},
we will prove the following Lemma.

\begin{lemma} Assume that the phase functions
$$
\phi(\xt,y',y'',v, w)=\phi_1(y',v)+\xt\phi_2(\xt,y',y'',v, w)
$$
and
$$
\phit(\xt,y',y'',\tilde v, \tilde w)=\tilde \phi_1(y',\tilde v)+\xt \tilde 
\phi_2(\xt,y',y'', \tilde v, \tilde w)
$$
parametrize a neighbourhood $U$ of $L$ containing
\begin{multline}
q = (0, y'_0, y''_0, d_{(x,\xt,y',y'')}(\tilde \phi/x)(q'))
= (0, y'_0, y''_0, d_{(x,\xt,y',y'')}(\tilde \phi/x)(\tilde q')), \\
q' = (0, y'_0, y''_0, v_0, w_0) \in C_\phi, \quad \tilde q' = 
(0, y'_0, y''_0, \tilde v_0, \tilde w_0) \in C_{\tilde \phi}.
\end{multline}
If $\dim v = \dim \tilde v$, $\dim w = \dim \tilde w$, and the signatures
of $d^2_{vv} \phi_1$ and $d^2_{ww} \phi_2$ at $(y'_0, v_0)$ and 
$(0, y'_0, y''_0, v_0, w_0)$ are equal to the corresponding signatures for
$\phit$, then the two phase functions are equivalent near $q$. That is,
there is a coordinate change $\tilde v(\xt, y', y'', v, w)$ and
$\tilde w(\xt, y', y'', v, w)$ mapping $q'$ to $\tilde q'$ such that
$$
\phit(\xt, y', y'', \tilde v(\xt, y', y'', v, w), \tilde w(\xt, y', y'', v, w))
= \phi(\xt, y', y'', v, w)
$$
in a neighbourhood of $q'$.
\end{lemma}

\begin{proof}
To prove this we follow the structure of H\"ormander's proof. The first
step is to show that $\phit$ is equivalent to a phase function 
$\psi(\xt, y', y'', v, w)$ that agrees with $\phi$ at the set $C_\phi$,
given by \eqref{eq:corner-22} to second order. This is accomplished by
pulling $\phit$ back with a diffeomorphism of the form
$$
(\xt, y', y'', v, w) \mapsto \big( \xt, y', y'', V(\xt, y', y'', v, w),
W(\xt, y', y'', v, w) \big)
$$
that restricts to the canonical diffeomorphism from $C_\phi$ to
$C_{\phit}$. This diffeomorphism is given by the identifications
of $C_\phi$ and $C_{\phit}$ with $L$. The mapping $F$ may be constructed
as in the first part of H\"ormander's proof.

In the second step, we first note that by the equivalence result for
Legendrians on manifolds with boundary (\cite{RBMZw}, Proposition 5), 
there is a change of variables
$v' = v'(y', v)$ near $(y'_0, v_0)$ such that
$$
\psi(\xt, y', y'', v', w) = \phi_1(y', v) + O(\xt).
$$
Therefore, we may assume that $\phi_1 = \psi_1$. If that is so, then
$\phi - \psi$ is $O(\xt)$, in addition to vanishing to second order
at $C_\phi = \{ \pa_v \phi = 0, \pa_w \phi_2 = 0 \}$. Therefore,
one can express
\begin{equation}
\psi - \phi = \frac{\xt}{2} \Big( \sum \frac{\pa \phi}{\pa v_j}
\frac{\pa \phi}{\pa v_k} b^{11}_{jk} + 2\sum \frac{\pa \phi}{\pa v_j}
\frac{\pa \phi_2}{\pa w_k} b^{12}_{jk} + \sum \frac{\pa \phi_2}{\pa w_j} 
\frac{\pa \phi_2}{\pa w_k} b^{22}_{jk} \Big)
\label{e1}
\end{equation}
for some smooth functions $b^{ab}_{jk}$ of all variables. 

We look for a change of variables of the form 
\begin{equation}\begin{aligned}
(\tilde v - v)_j &= \xt \sum \frac{\pa \phi}{\pa v_k} a^{11}_{jk} +
\xt \sum \frac{\pa \phi_2}{\pa w_k} a^{12}_{jk} \\
(\tilde w - w)_j &= \sum \frac{\pa \phi}{\pa v_k} a^{21}_{jk} +
\sum \frac{\pa \phi_2}{\pa w_k} a^{22}_{jk} .
\end{aligned}
\label{e2}\end{equation}
Using Taylor's theorem, we may write
\begin{multline}
\phi(\xt, y', y'', \tilde v, \tilde w) - \phi(\xt, y', y'', v, w) 
= \sum (\tilde v - v)_j \frac{\pa \phi}{\pa v_j} +
\xt \sum (\tilde w - w)_j \frac{\pa \phi_2}{\pa w_j} \\
+ \sum (\tilde v - v)_j (\tilde v - v)_k c^{11}_{jk}
+ \xt\sum (\tilde v - v)_j (\tilde w - w)_k c^{12}_{jk}
+ \xt\sum (\tilde w - w)_j (\tilde w - w)_k c^{22}_{jk}
\label{e3}
\end{multline}
for some smooth functions $c^{ab}_{jk}$. We wish to equate $\psi$ and
$\phi(\cdot, \tilde v, \tilde w)$. Substituting \eqref{e2} into
\eqref{e3} and equating with \eqref{e1} gives an expression which is
quadratic in the $\pa_v \phi$ and $\pa_w \phi_2$. Matching coefficients
gives the matrix equation
$$
A = B + Q(A, C, \xt)
$$
to be solved, where
$$
A = \begin{pmatrix} a^{11} & a^{12} \\ a^{21} & a^{22} \end{pmatrix},
$$
and similarly for the others. The expression $Q(A, C, x)$ is a polynomial
in the entries which is homogeneous of degree two in $A$. 

By the inverse function theorem,
this has a solution $A = A(B, C, \xt)$
if $B$ is sufficiently small. The proof is completed
by showing that there is a path $\psi_t$ of phase functions parametrizing
$L$ that interpolate between $\phi$ and $\psi$. We define $\psi_t$ by
considering a path of $b^{ab}_{jk}$'s and using \eqref{e1}. 

Recall that the nondegeneracy hypothesis is that
\begin{equation}\begin{aligned}
d_{y',v} \frac{\pa \psi}{\pa v_i} \text{ are independent, } \\
d_{y'',w} \frac{\pa \psi_2}{\pa w_i} \text{ are independent } 
\end{aligned}
\quad \text{ at } q'.
\end{equation}
At this point, by \eqref{e1},
$$
d^2_{vv} \psi = d^2_{vv} \phi \quad \text{and} \quad 
d^2_{y'v} \psi = d^2_{y'v} \phi 
$$
regardless of the value of $B$, so the first condition is automatically
satisfied. The other derivatives are
\begin{equation}\begin{aligned}
d^2_{ww} \psi &= d^2_{ww} \phi + (d^2_{ww} \phi) (b^{22}) (d^2_{ww} \phi) \\
d^2_{y''w} \psi &= d^2_{y''w} \phi + (d^2_{y''w} \phi) (b^{22}) (d^2_{ww} \phi) 
\end{aligned}\end{equation}
where $b^{22}$ is evaluated at $q'$. Thus the second nondegeneracy condition
at $q'$ is that $\Id + (b^{22}) (d^2_{ww}) \phi$ is nondegenerate.
As shown in the last part of H\"ormander's proof, the condition that
the signatures of $d^2_{ww} \phi$ and $d^2_{ww} \psi$ are equal is
precisely the condition that interpolation between $\phi$ and $\psi$
is possible by nondegenerate phase functions. This completes the proof
of the lemma.
\end{proof}

Returning to the proof of the proposition, if the dimension and
signature hypotheses of the lemma are satisfied, then the phases are
equivalent in some neighbourhood of $q$ and thus, by a change of
variables, the function $u$ can be written with respect to $\phit$.
In the general case, one modifies $\phi$ and $\phit$ by adding a nondegenerate
quadratic form $\langle Cv', v' \rangle + \xt \langle Dw', w' \rangle$
in extra variables to $\phi$, and a similar expression to $\phit$, 
so as to satisfy the hypotheses of the
Lemma. This requires a compatibility mod 2 between the
difference of the dimension of $v$ and $w$ and the signature of the
corresponding Hessian, but this is always satisfied by Theorem 3.1.4 of
\cite{FIO1}. The function $u$ can be written with respect to the modified
$\phi$ since the effect of the quadratic term is just to multiply by a 
constant times powers of $x$ and $\xt$. By the Lemma, $u$ can be written with
respect to the modified $\phit$, and therefore, with respect to $\phit$ itself.
\end{proof}

Legendrian distributions with conic points are defined analogously
by oscillatory integrals of the form
\begin{equation}\begin{split}\label{eq:dist-20}
w(x'',\xt,y',y'')=\int_{[0,\infty)\times\Real^{k'}}
&e^{i\phi(\xt,y',y'',s,w)/x} a(x'',\xt,y',y'',x''/s,s,w)\\
&\quad
\left(\frac{x''}{s}\right)^{m-(k'+1)/2+N/4}s^{p-1+N/4}\xt^{r+N/4-f/2}\,ds\,dw
\end{split}\end{equation}
with $N=\dim M$, $a\in\Cinf_c([0,\epsilon)\times
U\times[0,\infty)\times[0,\infty)\times\Real^{k'})$, $U$ open in $\mf$,
$\phi$ is a phase function parametrizing a Legendrian with
conic points $(G,\Gs)$ on $U$. Note that $k'+1$ appears instead of $k'$
since there are $k'+1$ parameters: $w$ and $s$. The integral converges
absolutely if $x''>0$ and $\xt>0$ (so $x>0$ and we are in the
interior of $M$) giving a smooth function in the interior of
$M$. Indeed, if $x''>c>0$ and the support of $a$ is in $x''/s<c'$, then
on the support of $a$, $s>c/c'$, so there is no problem with the
$s$ integral at $0$. Also, if $a$ is
supported in $s>0$, then this gives a parametrization completely
analogous to \eqref{eq:dist-1} (without the $v$ variables, since we
are assuming that we do not need them), i.e.\ in this case we simply
have a distribution in $\Isf(M,\Gs;\Osfh)$ with order $m$ on $\mf$ and
$r$ on $H$. On the other hand, if
$a$ is supported in $x''/s>0$, then we can change variables $\rho=s/x''$
(in place of $s$) allows us to write the integral as one analogous
to \eqref{eq:dist-1} with parameters $(\rho,v)$, so the result is in
$\Isf(M,G,\Osfh)$ with order $p$ on $\mf$ and order $r$ on $H$.
In addition, if we just keep $x''$ away from
$0$ (but possibly let $\xt$ to $0$), i.e.\ in the interior of $H$,
the previous argument showing $s>0$ on the support of $a$ works
whereby we conclude that the oscillatory integral in \eqref{eq:dist-20}
reduces to an oscillatory section of the form \eqref{eq:dist-2} without
any $v$ variables to integrate out. In particular, localized in the
interior of $H$, the integral gives an element of $\Isf(M,G;\Osfh)$
of order $r$ on $H$ (and infinite order vanishing at $\mf$)
(or of $\Isf(M,\Gs;\Osfh)$ with order $r$ at $H$).

Since the order on $H\neq\mf$ is forced to be the same for both $\Gs$
and $G$, the order of a distribution associated to a Legendrian pair 
with conic points can
be captured by a single order family $\calK$, which we take to give the
orders on $G$, and an additional real number $p$ to give the order on
$\Gs$ at $\mf$.

\begin{Def}
A Legendrian distribution $u\in\Isf^{\calK,p}(M,\Gt;\Osfh)$
associated to a 
Legendre pair with conic points
$\Gt=(G,\Gs)$ of $\sfT M$ (with $\calK$ as in \eqref{eq:ord-fam})
is a distribution $u\in
\dist(M;\Osfh)$ with the property
that for any $H\in M_1(M)\setminus\{\mf\}$ and for any
$\psi\in\Cinf_c(M)$ supported away from $\bigcup\{H'\in M_1(M):\ H'\neq\mf,
\ H'\neq H\}$,
$\psi u$ is of the form
\begin{equation}
\psi u=u_0+u_1+u_2+\sum_{j=1}^J w_j\cdot\nu_j
\end{equation}
where $u_0\in\Isf^{\calK'}(M,\Gs;\Osfh)$ with $\calK'(H)=\calK(H)$ if
$H\neq\mf$, $\calK'(\mf)=p$, $u_1\in\Isf^{\calK}(M,G;\Osfh)$,
$u_2\in I_{\scl,\text{c}}^{\calK(\mf),p}(\Mt,\Gt;\Osfh)$
($\text{c}$ meaning compact support in $\Mt$),
$\nu_j\in\Cinf(M;\Osfh)$, and $w_j\in
\dist(M)$ are given by oscillatory integrals as in \eqref{eq:dist-20}.
\end{Def}

\begin{rem}
Again, for the sake of convenience, we occasionally use the notation
$\Isf^{m,p,r}(M,\Gt,\Osfh)$
when we localize a Legendre distribution, associated to an
intersecting pair $\Gt$ with conic points, near a hypersurface $H$
so that only the order $m$ on $G$ and $p$ on $\Gs$ at $\mf$, and the
(common) order $r$ at $H$ are relevant.
\end{rem}

Again, for this class of distributions to be useful, we need to know that it
is independent of the choice of parametrizations. This follows in the
same way as above once we have the following result about equivalence
of phase functions parametrizing Legendrian pairs with conic points.

\begin{lemma} Let $(G, \Gs)$ be a Legendrian pair with conic points
and 
$$
\phi(\xt,y',y'',s, w)=\phi_1(y')+\xt(\phi_2(\xt, y',y'') + s\psi(y', s, w))
$$
and
$$
\phit(\xt,y',y'',\tilde s, \tilde w)=\phi_1(y')+\xt (
\phi_2(\xt,y',y'') + \tilde s\tilde \psi(\xt, y',y'',\tilde s,\tilde w))
$$
two phase functions parametrizing a neighbourhood $U$ of $L$ containing
\begin{gather}
q = (0, y'_0, y''_0, -\phi_1, -\frac{\psi}{|d_{y''}\psi|}, 
\frac{d_{y'}\psi}{|d_{y''}\psi|}, 
0, \widehat{d_{y''} \psi})(q') \\
= (0, y'_0, y''_0, -\phi_1, -\frac{\tilde \psi}{|d_{y''}\tilde \psi|}, 
\frac{d_{y'}\tilde \psi}{|d_{y''}\tilde \psi|}, 
0, \widehat{d_{y''} \tilde \psi})(\tilde q'), \\
q' = (0, y'_0, y''_0, 0, w_0) \in C_\phi , \quad \tilde q' =
(0, y'_0, y''_0, 0, \tilde w_0) \in C_{\tilde \phi}.
\end{gather}
If $\dim w = \dim \tilde w$, and the signatures
of $d^2_{ww} \psi$ and $d^2_{\tilde w \tilde w} \tilde \psi$ 
at $(0, y'_0, y''_0, 0, w_0)$ are equal, 
then the two phase functions are locally equivalent.
\end{lemma}

\begin{proof}
Again we follow the structure of H\"ormander's proof. First,
we may assume that $\phit$ agrees with $\phi$ at the set $C_\phi$,
given by \eqref{eq:corner-22p} to second order. As before, this is possible
if we can find a diffeomorphism of the form
$$
(\xt, y', y'', s, w) \mapsto \big( \xt, y', y'', S(\xt, y', y'', s, w),
W(\xt, y', y'', s, w) \big)
$$
that restricts to the canonical diffeomorphism from $C_\phi$ to
$C_{\phit}$. This is simply a version of the result from
\cite{RBMZw}, Proposition 7, with $y'$ an extra parameter. 

In the second step, we fix $y'$ and $\xt = 0$ and use the result of
Melrose and Zworski for Legendrian pairs with conic points, which says
that there is a coordinate change mapping $\psi(0,y',y'',s,w)$ to
$\tilde \psi(0,y',y'',\tilde s, \tilde w)$. We can deduce the same result,
with the same proof, letting $y'$ vary parametrically. Thus we can assume
that $\psi = \tilde \psi$ at $\xt = 0$. Consequently, 
$\phit - \phi = O(\xt^2 s)$. Thus, we may write 
\begin{equation}
\phit - \phi = \frac{\xt^2 s}{2} \Big( (\psi + s \psi_s)^2 b_{00}
+ 2 \sum (\psi + s \psi_s)(\psi_{w_j}) b_{0j} + \sum (\psi_{w_j})(\psi_{w_k})
b_{jk} \Big).
\label{f1}
\end{equation}
We look for a change of variables of the form 
\begin{equation}\begin{aligned}
(\tilde s - s) &= s \big( (\psi + s \psi_s) a_{00} + \sum \psi_{w_j} a_{0j}
\big) \\
(\tilde w - w)_j &= \sum (\psi + s \psi_s) a_{j0} + \sum
\psi_{w_k} a_{jk} 
\end{aligned}
\label{f2}\end{equation}
Using Taylor's theorem, we may write
\begin{multline}
\phi(\xt, y', y'', \tilde s, \tilde w) - \phi(\xt, y', y'', s, w) 
= \xt \Big( (\tilde s - s) (\psi + s \psi_s) +  
s \sum (\tilde w - w)_j \psi_{w_j} \\
+ \sum (\tilde s - s)^2 c_{00}
+ 2s\sum (\tilde s - s) (\tilde w - w)_k c_{0j}
+ s\sum (\tilde w - w)_j (\tilde w - w)_k c_{jk} \Big)
\label{f3}
\end{multline}
for some smooth functions $c_{\alpha k}$. Equating $\phit$ and
$\phi(\cdot, \tilde s, \tilde w)$ gives the matrix equation
$$
A = \xt B + Q(A, C, s)
$$
where $Q$ is homogeneous of degree two in $A$. 
This always has a solution for small $\xt$ by the inverse function theorem.
This completes the proof of the Lemma. 
\end{proof}

\section{Poisson operators and spectral projections}\label{sec:Poisson}
In this section we will prove a statement asserted in the introduction,
namely, that the spectral projection at energy $\l^2 > 0$ can be expressed in
terms of the Poisson operators. 

We briefly recall the structure of the Poisson operator from the paper
of Melrose and Zworski~\cite{RBMZw}. We also remark that a somewhat
different approach has been presented in \cite{Vasy:Geometric}; in this
approach the form \eqref{eq:rev-con-16} of Legendre distributions
associated to intersecting pairs with conic points arises very explicitly.
Following \cite{RBMZw}, we often make $P(\pm
\lambda)$ a map from half-densities to half-densities to simplify
some of the notation. The correspondence between the smooth functions
and half-densities is given by the trivialization of the half-density
bundles by the Riemannian densities; so for example we
have
\begin{equation}
P(\pm\lambda)(a|dh|^{1/2})=(P(\pm\lambda)a)\,|dg|^{1/2}
,\quad a\in\Cinf(\bX).
\end{equation}
With this normalization the kernel of $P(\pm\lambda)$ will be a section of
the kernel density bundle
\begin{equation}
\KDsc=\pi_L^*\Osch X\otimes\pi_R^*\Oh\bX.
\end{equation}
where $\pi_L:X\times\bX\to X$ and
$\pi_R:X\times\bX\to\bX$ are the projections. We remark that
smooth sections of $\KDsc$ are of the form $a|dg|^\half|dh|^\half$,
$a\in\Cinf(X\times\bX)$, while
smooth sections of $\Osch(X\times\bX)$ are of the form $a|dg|^\half
|dh|^\half x^{-(n-1)/2}$. Melrose and Zworski constructed
the kernel of $P(\lambda)$ in \cite{RBMZw} as a 
section of $\Osch(X\times\bX)$, essentially by identifying $\KDsc$ with
$\Osch(X\times\bX)$ (via the mapping $a|dg|^\half|dh|^\half\mapsto
a|dg|^\half
|dh|^\half x^{-(n-1)/2}$), given a
choice of boundary defining function $x$.
    
We introduce the following Legendre submanifolds of
$\sct_{\bX\times\bX}(X\times\bX)$:
\begin{equation}
\Gs(\l)=\{(y,y',\tau,\mu,\mu'):\ \mu=0,\ \mu'=0,\ \tau=-\l\},
\label{gs}
\end{equation}
and for $\lambda > 0$ let
\begin{equation}\begin{split}
\label{eq:def-Gp}
G(\lambda)=&\{(y,y',\tau,\mu,\mu'):\ (y,\muh)=\exp(sH_{\half h})(y',\muh'),
\ \tau=\lambda\cos s, \\
&\ \mu=\lambda(\sin s)\muh,\ \mu'=-\lambda(\sin s)\muh',
\ s\in[0,\pi),\ (y',\muh')\in S^*\bX\},
\end{split}\end{equation}
\begin{equation}\begin{split}\label{eq:def-Gm}
G(-\lambda)=&\{(y,y',\tau,\mu,\mu'):\ (y,\muh)=
\exp((s-\pi)H_{\half h})(y',\muh'),
\ \tau=\lambda\cos s,\\
&\ \mu=\lambda(\sin s)\muh,\ \mu'=-\lambda(\sin s)\muh',
\ s\in(0,\pi],\ (y',\muh')\in S^*\bX\}.
\end{split}\end{equation}
We recall that $G(\lambda)$ is actually smooth at $s=0$,
and $G(-\l)$ is smooth at $s=\pi$. 
We introduce notation for the Legendre pair
\begin{equation}
\Gt(\lambda)=(G(\lambda),\Gs(\lambda)), \quad \Gt(-\lambda)=(G(-\lambda),
\Gs(-\lambda)).
\end{equation}
For later use, we also define $G_{y_0} = G_{y_0}(\l)$ to be the Legendrian in 
$\sct_{\pa X} X$
\begin{equation}
G_{y_0}(\l) = \{ (y', \tau, \mu') \mid (y_0, y', \tau, \mu, \mu') \in G(\l) 
\text{ for some } \mu \}.
\label{Gy0}
\end{equation}
Thus, Melrose and Zworski construct the kernel of $P(\pm\lambda)$, $\l > 0$,
as a Legendre distribution associated to a Legendre pair with conic points:
\begin{equation}\label{eq:param-4}
P(\pm \l) \in \Isc^{-(2n-1)/4,-1/4}(X\times\bX,\Gt(\pm\lambda),\KDsc).
\end{equation}

To express the spectral projection in terms of the Poisson operator we first need to
note the following technical result. 

\begin{lemma}\label{continuity} 
The resolvent $R(\sigma)$ extends to a continuous family of operators
from $x^{1/2 + \ep} L^2$ to $x^{-1/2 - \ep} L^2$ on each of the closed sets
$$ \{ \sigma \in \CC \mid \re \sigma \geq \ep \text{ and }
\pm \im \sigma \geq 0 \}
$$ 
for any $\ep > 0$. 
\end{lemma} 

\begin{proof} (sketch) This result follows readily from \cite{RBMSpec}
(see also \cite[Proposition~2.5]{Vasy:Geometric} for more details) by
noting that the estimates of Propositions 8 and 9, which are uniform scattering
wavefront set estimates on $R(\sigma + i\ep)$ as $\ep \to 0$, $\sigma \in \RR$,
can also be shown
to be uniform on a small $\sigma$-interval not containing zero. Then the argument
in the proof of Proposition 14 extends to show joint continuity of $R(\sigma + i\ep)$
in both $\sigma$ and $\ep$.
\end{proof}

This allows us to obtain the desired formula:

\begin{lemma}\label{Poisson=sp} 
The generalized spectral projection at energy $\l$ is 
$$
\sp(\l) = (4\pi \l)^{-1} P(\l) P(\l)^* = (4\pi \l)^{-1} P(-\l) P(-\l)^* .
$$
\end{lemma}

\begin{proof} Regarding the spectral projection $E_{(a^2, b^2)}$ as a map
$\dCinf(X) \to \dist(X)$, and using Lemma~\ref{continuity} 
we may take the pointwise limit as $\ep \to 0$ in the
integral \eqref{Stone} to obtain
$$
E_{(a^2, b^2)} = \frac1{2\pi i} \int_a^b \Big( R(\l^2 + i0) - R(\l^2 - i0) \Big) 
2\l d\l. 
$$
Thus it remains to show that 
\begin{equation}\label{eq:intro-5}
R(\lambda^2+i0)-R(\lambda^2-i0)=\frac{i}{2\lambda}P(\lambda)P(\lambda)^*
=\frac{i}{2\lambda}P(-\lambda)P(-\lambda)^*
\end{equation}
as operators $\dCinf(X)\to\dist(X)$. 
To derive \eqref{eq:intro-5}, first recall that $R(\lambda^2\pm i0)f$,
$f\in\dCinf(X)$, has an asymptotic expansion similar to
\eqref{eq:smooth-exp-1}. Namely,
\begin{equation}\label{eq:intro-6}
R(\lambda^2\pm i0)f=\pm
e^{\pm i\lambda/x}x^{(n-1)/2}v_\pm,\ v_\pm\in\Cinf(X).
\end{equation}
Thus,
\begin{equation}\label{eq:intro-7}
R(\lambda^2+i0)f-R(\lambda^2-i0)f=P(\lambda)(v_-|_{\bX}),
\end{equation}
since both sides are generalized eigenfunctions of $H$,
have an expansion as in \eqref{eq:smooth-exp-1} and the
leading coefficients of the $e^{-i\lambda/x}$ term agree, namely they are
$v_-|_{\bX}$. Thus, we only have to determine $v_-|_{\bX}$
in \eqref{eq:intro-6}.

For this purpose, recall the boundary pairing formula from Melrose's paper
\cite[Proposition~13]{RBMSpec}. Thus, suppose that
\begin{equation}\label{eq:intro-8}
u_j=e^{-i\lambda/x}x^{(n-1)/2}v_{j,-}+e^{i\lambda/x}x^{
(n-1)/2}v_{j,+},\ v_{j,\pm}\in\Cinf(X),
\end{equation}
and $f_j=(H-\lambda^2)u_j\in\dCinf(X)$. Let $a_{j,\pm}=v_{j,\pm}|_{\bX}$.
Then
\begin{equation}\label{eq:intro-9}
2i\lambda\int_{\bX}(a_{1,+}\,\overline{a_{2,+}}-a_{1,-}\,\overline{a_{2,-}})
\,dh=\int_X(u_1\,\overline{f_2}-f_1\,\overline{u_2})\,dg.
\end{equation}
We apply this result with $u_1=-R(\lambda^2-i0)f$, $u_2=P(\lambda)
a$, $a\in\Cinf(\bX)$. Thus, $f_1=-f$, $f_2=0$, $a_{2,-}=a$, $a_{1,+}=0$
and $a_{1,-}=v_-|_{\bX}$ with the notation of \eqref{eq:intro-6}.
We conclude that
\begin{equation}\label{eq:intro-10}
2i\lambda\int_{\bX}v_-|_{\bX}\,\overline{a}
\,dh=-\int_X f\,\overline{P(\lambda)a}\,dg.
\end{equation}
Thus, $v_-|_{\bX}=\frac{i}{2\lambda}P(\lambda)^*f$.
Combining this with \eqref{eq:intro-7}, we deduce the first equality in
\eqref{eq:intro-5}. The second equality can be proved similarly.
\end{proof}


\section{An example}\label{sec:example}

Before we discuss the general case, let us consider the case of the resolvent
of the Laplacian on Euclidean space. Using the Fourier transform, it is not hard
to get expressions for the kernels of the spectral projections and resolvent.

The spectral projection at energy $\l^2$ is
$$
\frac{\l^{(n-2)}}{2 (2\pi)^n} \int_{S^{n-1}} e^{i(z'-z'')\cdot \l \omega} \ d\omega
\ |dz'dz''|^\half
$$
Since $z' - z''$ is a homogeneous function of degree one on
$\RR^{2n}$, this is a Legendrian distribution on the radial compactification of
$\RR^{2n}$; denote the Legendre submanifold, which will be described in
detail below, by $L$. 
Let us see what sort of object this is on $\XXb$, where here $X$ is
the radial compactification of $\RR^n$. In fact, as noted in 
\cite{Hassell-Vasy:Symbolic},
$$
\XXb = [\ol{\RR^{2n}}; S', S'']
$$
where $S'$ is the submanifold of the boundary of $\ol{\RR^{2n}}$ which is the
intersection of the plane $z' = 0$ with the boundary, and similarly for $S''$. 
Let $S$ be the intersection of the plane $z'=z''$ with the boundary. Away from
$S$, $|z'-z''|$ is a large parameter, and stationary phase in $\omega$ shows
that the Legendrian $L$ is the union of two projectable Legendrians $L_\pm$
parametrized by $\pm \l |z'-z''|$. At
$S$, stationary phase is no longer valid, and $L$ is no longer projectable. In fact
the part of $L$ lying above $S$ is
\begin{equation}\label{L-above-S}
\{ \hat z' = \hat z'', |z'|/|z''| = 1, \mu' = -\mu'', \tau = 0, \eta^2 + 2|\mu'|^2 
= 2\l^2 \}
\end{equation}
in the coordinates \eqref{eq:sp-a-10}, which is codimension one in $L$. It is
the intersection of $\sci N^*\diag_\bl$ with $\Sigma(\Lap - 2\l^2)$, where
$\sci N^*\diag_\bl$ is the boundary of the closure of $N^* \diag_\bl$ in 
$\sct \ol{\RR^{2n}}$,
$$
\sci N^*\diag_\bl = \{ y' = y'', \sigma = 1, \mu' = -\mu'', \tau = 0 \},
$$
 and $\Lap = \Lap_L + \Lap_R$ is the Laplacian in $2n$ variables. 
Notice that $(\Lap - 2\l^2) \sp(\l) = 0$, so we know that 
$L \subset \Sigma(\Lap - 2\l^2)$ and the bicharacteristic flow of $\Lap - 2\l^2$ is
tangent to $L$. Under bicharacteristic flow,
$\dot \tau = -2(|\mu'|^2 + |\mu''|^2 + \eta^2)$, which is strictly negative on
\eqref{L-above-S}, so $L_\pm$ are in fact the flowouts from $L \cap \sci N^*\diag_\bl$
by backward/forward bicharacteristic flow. 

Next consider the structure of $\sp(\l)$ near lb (the case of rb is similar). 
After blowup of
say $S''$, creating the boundary hypersurface lb, we have smooth coordinates
$x'' = 1/|z''|$, $\sigma = |z''|/|z'|$, $\hat z'$ and $\hat z''$. Let
$x' = x'' \sigma$. Then near the corner $x'' = \sigma = 0$, the phase function
takes the form 
$$
\l|z' - z''| = \l(|z' - z''|^2)^{1/2} 
= \l(1 - \sigma \hat z' \cdot \hat z'' + O(\sigma^2))/x',
$$
which parametrizes a Legendrian distribution on $\XXb$ near $\lb \cap \bfc$.
In fact we can read off that the Legendre submanifold $G_1$ is $\Gs(\l)$, and
the Legendre submanifold lying in the fibre above $\hat z'$ is $G_{\hat z'}(\l)$.

If there were a potential $V$ present, $H = \Lap + V$, say $V \in \dCinf(\RR^n)$,
then the lifts of $V$ from the left and right factor to $\ol{\RR^{2n}}$ would
be singular, and in fact Legendrian, at $S'$ and $S''$ respectively. Then the
analysis of \cite{Hassell:Plane} suggests that there would be another
Legendrian propagating out from $S'$ and $S''$. This Legendrian may be considered
to be the origin of the second Legendre submanifold, denoted $\Ls$, that we shall
see is present in the spectral projection in the general case.

Next consider the resolvent $R(\l^2 + i0)$. This has kernel
$$
\int e^{i(z' - z'') \cdot \zeta} (|\zeta|^2 - \l^2 - i0)^{-1} \, d\zeta \ |dz'dz''|^\half.
$$
To understand this kernel, we need to cut it into pieces. First we introduce
a cutoff in $\zeta$, $1 = \psi(\zeta) + (1 - \psi(\zeta))$ where $\psi$ vanishes
in a neighbourhood of $|\zeta| = \l$ and $1 - \psi$ has compact support. Then
$$
\int e^{i(z' - z'') \cdot \zeta} (|\zeta|^2 - \l^2 - i0)^{-1} \psi(\zeta)
\, d\zeta \ |dz'dz''|^\half
$$
is (the kernel of)
a constant coefficient pseudodifferential operator on $\RR^n$, and hence lies
in the scattering calculus. To deal with the other piece, we write
\begin{equation}\label{eq:plus-i0}
(|\zeta|^2 - \l^2 - i0)^{-1} = (|\zeta| + \l)^{-1} (|\zeta| - \l - i0)^{-1}
= \frac{i}{|\zeta| + \l} \int_0^\infty e^{-i(|\zeta| - \l)\ol t} d\ol t,
\end{equation}
change variable to $t = \ol t |z|$ and substitute into the previous expression. The result
is 
\begin{equation}\label{eq:ex-res}
\int \int_0^\infty dt \ e^{i(z'-z'')\cdot \zeta - (|\zeta| - \l)t|z|} (1- \psi(\zeta)) 
\frac{i|z|}{|\zeta| + \l} \, |\zeta|^{n-1} d|\zeta| d\hat\zeta  \ |dz'dz''|^\half.
\end{equation}
The phase $\Phi = (z' - z'')\cdot \zeta - (|\zeta| - \l)t |z|$
now parametrizes the pair of intersecting Legendrian submanifolds 
$(L_+, \sci N^*\diag_\bl)$ (see section~\ref{sec:sc-calculus}). 
Notice that away from $S$, ie, away from $\sci N^*\diag_\bl$,
$\Phi$ is an exact nondegenerate quadratic in $(|\zeta|,t)$ and 
the rest of the integral is independent
of $t$. Therefore, the top order part of the stationary
phase approximation is accurate to all orders, and we find that microlocally
in this region, the resolvent is $2\pi i$ times the spectral projection. 
We have ignored some analytical details here: namely, the fact that the integrand
in \eqref{eq:ex-res} is not compactly supported in $t$, and we should retain a factor
$e^{-\ep \ol t}$ in the integral \eqref{eq:plus-i0} where $\ep$ is later sent to zero. 
These details will be taken care of when we treat the general situation.

As with the spectral projection, we expect to get an extra Legendrian whenever
we perturb the Laplacian by a potential. 

In the rest of this paper we show that when we look at these kernels on $\RR^{2n}$
from the point of scattering microlocal analysis they become a rather good guide
to the general situation.

\section{Legendre structure of the spectral projections}\label{sec:sp}
In this section we will show that the spectral projection $\sp(\l)$
is a Legendre distribution
associated to an intersecting pair $(L = L(\l), \Ls = \Ls(\l) \cup \Ls(-\l))$ of Legendre
submanifolds with conic points, on $M = \XXb$. 

We now introduce the Legendre submanifolds $L$ and $\Ls$. 
First note that the another coordinate system near $\ffb$,
different from that discussed at the beginning of
Section~\ref{sec:corner-geo}, which is
sometimes more convenient, is given by replacing $\sigma$ by $\theta
=\arctan\sigma$.
Thus, \eqref{eq:fib-b-2} shows that
we can write $\alpha\in\sfT_p M$, $p=(x'',\theta,y',y'')$ near
$\lf\cap\ffb$, as
\begin{equation}\label{eq:sp-a-15}
\alpha=\tau\frac{dx}{x^2}+\eta\frac{d\theta}{x}+\mut'\cdot\frac{dy'}{x}
+\mut''\cdot\frac{dy''}{x''}.
\end{equation}
Note that we have $x=x''\,\sin\theta$, so we should take $\xt=\sin\theta$
in \eqref{eq:corner-9-a}. Comparison with \eqref{eq:corner-9-a} also
shows that $\eta$ is defined slightly differently here, but this should
not cause any trouble. Also, we use the notation $\mut'$, $\mut''$ in place
of $\mu'$, $\mu''$, since these will be reserved from the coordinates
coming from the bundle decomposition \eqref{eq:str-b-7}.

We note that
an attractive feature of the coordinates $(x,\theta,y',y'',\tau,\eta,\mut',
\mut'')$ induced by \eqref{eq:sp-a-15} is that
in these coordinates the product metric
\begin{equation}
\gt=(\pbL)^* g+(\pbR)^* g
\end{equation}
becomes
\begin{equation}\label{eq:sp-a-17}
\gt=\frac{dx^2}{x^4}+\frac{d\theta^2+\cos^2\theta \,(\pbL)^*h+\sin^2
\theta\,(\pbR)^*h}{x^2}
\end{equation}
which is of the form
\begin{equation}
\gt=\frac{dx^2}{x^4}+\frac{\htil}{x^2},
\end{equation}
so $\gt$ is a scattering metric on $\Mt$.

Corresponding to the direct sum decomposition of $\sfT M$ given
by the pull-backs
of $\sct X$ by the stretched projections,
\eqref{eq:str-b-7}, we can also use coordinates
$(x'',\sigma,y',y'',\tau',\tau'',\mu',\mu'')$ on $\sct\Mt$ where we write
$q\in\sct_p \Mt$, $p=(x'',\sigma,y',y'')$, as
\begin{equation}\label{eq:sp-a-10}
q=\tau'\frac{dx'}{(x')^2}+\tau''\frac{dx''}{(x'')^2}
+\mu'\cdot\frac{dy'}{x'}+\mu''\cdot\frac{dy''}{x''};
\end{equation}
we may of course replace $\sigma$ by $\theta=\arctan\sigma$.
A reason why \eqref{eq:sp-a-10} is often convenient is that it matches
the structure of the composition $P(\lambda)P(\lambda)^*$.
The relationship between \eqref{eq:sp-a-10} and \eqref{eq:sp-a-15} is
given by
\begin{equation}\label{relnship}
\tau=\tau'\cos\theta+\tau''\sin\theta,\ \eta=\tau'\sin\theta-\tau''\cos
\theta,\ \mut'=(\cos\theta)\mu',\ \mut''=(\sin\theta)\mu''.
\end{equation}

Now, the Laplacian of the product metric $\gt$, dicussed in
\eqref{eq:sp-a-17}, is
\begin{equation}\label{Lap-gt}
\Delta_{\gt}=\Delta_{L}+\Delta_{R},
\end{equation}
$\Delta_L$ and $\Delta_R$ being the Laplacians of $g$ lifted from either
factor of $X$ using the product $X\times X$, and lifting further to the
blown up space $\XXb$. Similarly lifting the perturbation $V$, and
using that $(\Delta+V-\lambda^2)P(\lambda)=0$, we see that
\begin{equation}\label{eq:sp-13}
(\Delta_{\gt}+V_L+V_R-2\lambda^2) \sp(\l)=0,
\end{equation}
the density factors being
trivialized by the Riemannian density.
Thus, the kernel of $\sp(\l)$, also denoted by $\sp(\l)$, is a generalized
eigenfunction of
\begin{equation}
\Ht=\Delta_{\gt}+V_L+V_R=\Delta_L+\Delta_R+V_L+V_R\in\Diffsf^2(M).
\end{equation}
Away from $\lf\cup\rf$, i.e.\ on $\Mt$, $\Ht$ is a scattering
differential operator (just as every other operator in $\Diffsf(M)$),
$\Ht\in\Diffsc^2(\Mt)$,
so based on the paper of Melrose and Zworski \cite{RBMZw},
we can expect the appearance of Legendrian submanifolds
of $\sct \Mt$ so that $\sp(\l)$ becomes a (singular) Legendre distribution.
It is thus natural to expect that the full analysis, including the behavior
of the kernel up to $\lf\cup\rf$ would involve Legendre
distributions on the manifold $M$ with fibred boundary faces.

We first define the smooth Legendre submanifold $L$ of
$\sfT_{\ffb}M$
by using the parametrization that will appear in the composition $P(\lambda)
P(\lambda)^*$. Thus, we let
\begin{equation}\begin{split}\label{eq:sp-1}
L = L(\l) 
=&\{(\theta,y',y'',\tau',\tau'',\mu',\mu''):\ \exists(y,\muh)\in S^*\bX,
\ s,s'\in[0,\pi],\\
&\quad (\sin s)^2+(\sin s')^2>0,\ \text{s.t.}\\
&\quad \sigma = \tan \theta = \frac{\sin s'}{\sin s},
\ \tau'=\lambda\cos s',
\ \tau''=-\lambda\cos s,\\
\quad(y',\mu')&=\lambda\sin s'\exp(s'H_{\half h})(y,\muh),
(y'',\mu'')=-\lambda\sin s \exp(s H_{\half h})(y,\muh)\}\\
&\cup\{(\theta,y,y,\lambda,-\lambda,0,0):\ \theta\in[0,\pi/2],\ y\in\bX\}\\
&\cup\{(\theta,y,y,-\lambda,\lambda,0,0):\ \theta\in[0,\pi/2],\ y\in\bX\}.
\end{split}\end{equation}
(Note that the requirement $(\sin s)^2+(\sin s')^2>0$ just
means that $s$ and $s'$ cannot take values in $\{0,\pi\}$ at the same time.)
Although it
is not immediately obvious, with this definition $L$ is a smooth submanifold of
$\sfT_{\ffb}M$. 

Let us show this. 
Note that to get points away from $\rf$
and away from the last two sets in \eqref{eq:sp-1}, we must
have $\sin s>0$, so we can restrict the parameter $s$ to $(0,\pi)$ in this
region. Thus, $\tau''$, $\sigma$ and $(y'',\mu''/|\mu''|)$ (note that
$\mu''\neq 0$ as $\sin s\neq 0$) give coordinates on $L$ away from $\rf$
and the last two sets. A similar argument works near $\rf$ if we replace
$\sigma=\tan\theta$ by $\cot \theta$.
To see the smoothness near the last two
sets in the definition of $L$ above, note that near $\tau'=\lambda$,
$\tau''=-\lambda$, $\sigma\in[0,C)$ where $C>1$, $L$ is
given by
\begin{equation}\begin{split}\label{eq:sp-1-a}
\{(\sigma&,y',y'',\tau',\tau'',\mu',\mu''):\ \exists(y,\mu)\in T^*\bX,
\ |\mu|<C^{-1},\ \sigma\in[0,C)\ \text{s.t.}\\
&\tau'=\lambda(1-|\sigma\mu|^2)^{1/2},
\ \tau''=-\lambda(1-|\mu|^2)^{1/2},\\
&(y',\mu')=\lambda\exp(f(\sigma\mu)H_{\half h})(y,\sigma\mu),
\ (y'',\mu'')=-\lambda\exp(f(\mu)H_{\half h})(y,\mu)\}
\end{split}\end{equation}
where
\begin{equation}
f(\mu)=\frac{\arcsin(|\mu|)}{|\mu|}
\end{equation}
is smooth and $f(0)=1$. Thus, the differential of the map
\begin{equation}
(y,\mu)\mapsto-\lambda\exp(f(\mu)H_{\half h})(y,\mu)=(y'',\mu'')
\end{equation}
is invertible near $\mu=0$, so it gives a diffeomorphism near $|\mu|=0$.
Hence, $\sigma$ and $(y'',\mu'')$ give coordinates on $L$ in this
region, so $L$ is smooth here. Near $\rf$, but still
with $\tau'$ near $\lambda$, $\tau''=-\lambda$, we replace $\sigma=\tan\theta$
by $\cot\theta$ and deduce that $\cot\theta$ and $(y',\mu')$ give
coordinates on $L$ here.
Near $\tau'=-\lambda$,
$\tau''=\lambda$ a similar argument works if we replace $s$ by $\pi-s$,
$s'$ by $\pi-s'$, $(y,\muh)$ by $\exp(\pi H_{\half h})(y,\muh)$ in the
definition of $L$ above. Notice that the in the first set on the
right hand side of \eqref{eq:sp-1} (and also in the other pieces)
\begin{equation}\label{eq:sp-16-b}
(y'',\muh'')=\exp((s-s') H_{\half h})(y',\muh')
\end{equation}
for some appropriate $\muh'$, $\muh''$, so the distance of $y'$ and $y''$ is
at most $\pi$. Also notice that with the notation of
\eqref{eq:sp-a-15}, $\tau$ is given by
\begin{equation}\label{eq:tau-on-G}
\tau=\lambda(\cos s'\,\cos\theta-\cos s\,\sin\theta)
\end{equation}
so that $\tau\in(-\sqrt{2}\lambda,\sqrt{2}\lambda)$ on $L$,
i.e.\ $L$ is
disjoint from the `radial set' in $\sct_{\partial\Mt} \Mt$ induced by
$\Delta_{\gt}-2\lambda^2$.

\

Now we discuss the `ends' of $L$ with the goal of understanding
the composition
$P(\lambda)P(\lambda)^*$ in mind.
The closure of $L$ in $\sfT_{\ffb} M$ is
\begin{equation}
\cl L =L\cup F_{\lambda}\cup F_{-\lambda}
\end{equation}
where
\begin{equation}
F_{\l}=\{(\theta,y',y'',-\l,-\l,0,0):\ \exists\muh\in S^*_{y'}\bX
\ \text{s.t.}\ \exp(\pi H_{\half h})(y',\muh)\in S^*_{y''}\bX\}.
\end{equation}
The condition in the definition of $F_{\l}$ just amounts to saying that
there is a geodesic of length $\pi$ between $y'$ and $y''$; it arises if
we let $s\to 0$, $s'\to\pi$ (or vice versa) in \eqref{eq:sp-1}.
We also define
\begin{equation}
\Ls(\l)=\{(\theta,y',y'',-\l,-\l,0,0):\ y',y''\in\bX,\ \theta
\in[0,\pi/2]\},
\label{eq:L-sharp}\end{equation}
so $\Ls(\l)$ is a Legendrian submanifold of $\sct_{\ffb}\Mt$
for each $\l$ and
\begin{equation}
\cl L\cap\Ls(\pm \l)=F_{\pm\lambda}.
\end{equation}

\begin{prop} The pair 
\begin{equation}
\Lt = \Lt (\l) = (L(\l), \Ls(\l) \cup \Ls(-\l))
\end{equation}
is a pair of intersecting Legendre manifolds with conic points.
\end{prop}

\begin{proof}
First we show near $\bfc \cap \lb$ that $L$ and $\Ls(\pm \l)$ satisfy the conditions
of Definition~\ref{M-Leg-submfld}. This is clear for $\Ls(\pm \l)$, which 
projects under $\phit_{\lb}$ to 
$\Gs(\pm \l)$ (see \eqref{gs})
at the corner with fibres $\Gs(\pm \l)$. To see this for $L$, notice
that at the corner, i.e.\ when $\sigma = 0$, then by \eqref{eq:sp-1} we have
\begin{equation}\begin{gathered}
L \cap \sfT_{\lb \cap \bfc}\XXb = \{ (0, y', y', \pm \l, \mp \l, 0, 0) \} \\
\cup \{ (0, y', y'', \pm \l, -\l \cos s, 0, \mu'') 
\mid (y'', \mu'') = \mp \l \sin s \exp(H_{\half h})(y', \hat \mu') \} 
\end{gathered}\end{equation}
which also projects to $\Gs(\l) \cup \Gs(-\l)$, with fibre above $(y', \pm \l, 0)$ equal
to $G_{y'}(\pm \l)$. It is shown in \cite{RBMZw}
that $(G_{y'}(\pm \l), \Gs(\l) \cup \Gs(-\l))$ form
an intersecting pair with conic points. Since $\Lt$ is invariant under the
transformation $(z', z'') \to (z'', z')$, $\tau' \to -\tau''$, $\tau'' \to -\tau'$, 
$\mu' \to -\mu''$, $\mu'' \to -\mu'$, this is also the case near $\rb$. It remains
to show that in the interior of $\bfc$ that they form an intersecting pair
with conic points.

To see this,
notice that using the description \eqref{eq:sp-1} of $L$, near $\Ls(\l)$,
$s'$ is near
$\pi$ and $s$ is near $0$ there, so we may write, using \eqref{eq:sp-1},
$|\mu''|=\lambda\sin s$, $|\mu'|=\lambda\sin s'=\sigma|\mu''|$,
$s'=\pi-\arcsin(\sigma|\mu''|/\lambda)$ ($\arcsin$ takes values in
$[-\pi/2,\pi/2]$ as usual),
$(y,\muh)=\exp(-sH_{\half h})(y'',-\muh'')$,
and then near $\Ls(-\lambda)$,
$L$ is given by
\begin{equation}\begin{split}\label{eq:sp-1-b}
&\{(\sigma,y',y'',\tau',\tau'',\mu',\mu''):\ \mu''=|\mu''|\muh'',\\
&\qquad \tau'=-(\lambda^2-\sigma^2|\mu''|^2)^{1/2},
\ \tau''=-(\lambda^2-|\mu''|^2)^{1/2},\\
&\qquad (y',\mu')=\sigma|\mu''|\exp((\pi-\arcsin(\sigma|\mu''|/\lambda)
-\arcsin(|\mu''|/\lambda))H_{\half h})(y'',-\muh'')\},
\end{split}\end{equation}
i.e.\ the $N-1$ variables $\sigma,y'',|\mu''|,\muh''$ give coordinates
on $L$ near $\Ls(-\lambda)$. This means that $L$
is the image of a
smooth manifold in $\mu''$-polar coordinates. However, notice also that
$|\mu'|=\sigma|\mu''|$, i.e.\ away from $\rf$ it is a bounded multiple
of $|\mu''|$, and similarly $\tau'-\tau''$ is a bounded multiple of
$|\mu''|$, so $L$ is the image of a smooth manifold in $(\tau'-\tau'',
\mu',\mu'')$-polar coordinates. A similar argument works near $\rf$
with the role of $\mu'$ and $\mu''$ interchanged, so we can
conclude that $L$ is the image of a
closed embedded submanifold $\hat L$ of
\begin{equation}
[\sfT_{\ffb}M;\Span\{\frac{dx'}{(x')^2}+\frac{dx''}{(x'')^2}\}]
=[\sfT_{\ffb}M;\Span\{d((\sigma+1)/x')\}],
\end{equation}
under the blow-down map.
Also $\hat L$
is transversal to the front face of the blow-up (since away from $\rf$
resp.\ $\lf$, $|\mu''|$ resp.\ $|\mu'|$
is a coordinate on $\Gh_\lambda$).
These, together with corresponding arguments
near $\Ls(-\lambda)$, establish the result.
\end{proof}

\

Finally we come to one of the main results of this paper.

\begin{thm}\label{thm:sp-proj}
The kernel of the spectral projection at energy $\l^2$ is
a Legendre distribution associated to the intersecting pair with conic
points, $\Lt = \Lt (\l) = (L(\l), \Ls(\l) \cup \Ls(-\l))$. More
precisely, we have
\begin{equation}
\sp(\l) = (4\pi \l)^{-1} P(\lambda)P(\lambda)^* \in\Isf^{\calK,\frac{n-2}{2}}
(\XXb, \Lt, \Osfh)
\end{equation}
with $\calK(\lb)=\calK(\rb) =\frac{n-1}{2}$, $\calK(\mf) = -\half$.
\end{thm}

\begin{proof}
We prove the result by microlocalizing in each factor of the composition
$P(\l) P(\l)^* = P(-\l) P(-\l)^*$ and showing quite explicitly that we get
a Legendre conic pair associated to $(L, \Ls)$. Thus, let $Q_i$ and
$R_j$ be pseudodifferential operators on $X$, $\sum_i Q_i = \sum_j R_j = \Id$,
associated to a microlocal partition of unity,
such that each $Q_i P(\l)$, $Q_i P(-\l)$, $R_j P(\l)$ and $R_j P(-\l)$ (which
are all Legendrian conic pairs associated to $(G(\l), \Gs(\l))$ or 
$(G(-\l), \Gs(-\l))$) can be parametrized by a single phase function. Then
we write 
$$
\sp(\l) = \sum_{i,j} Q_i P(\l) P(\l)^* R_j^*
$$ 
and analyze each term separately.  There
are three qualitatively different types of parametrizations that can occur: 
(1) if $Q_i$ is microlocalized near $\Gs(-\l)$ then $Q_i P(\l)$ has a phase function
of the form $-\cos d(y', y)/x'$; (2) if $Q_i$ is microlocalized where $-\l < \tau
< \l$ then $Q_i P(\l)$ has phase function $\phi(y', y, v)/x'$ parametrizing
$G(\l)$; and (3) if $Q_i$ is localized near $\Gs(\l)$ then $Q_i P(\l)$ has phase
function $(1 + s\psi(y', y, s, w))/x'$ parametrizing $(G(\l), \Gs(\l))$. Similarly
for $Q_i P(-\l)$, $(R_j P(\l))^*$, etc. A subcase of the last case is $\psi \equiv 0$,
i.e.\ the phase function parametrizes only $\Gs(\l)$. 

Let us denote by $(a,b)$ the case that $Q_i$ falls into case $a$ and $R_j$
into case $b$. Then there are nine possibilities in all, but we can eliminate
several of them. First, writing $Q_i P(\l) P(\l)^* R_j^*$ as
$Q_i P(-\l) P(-\l)^* R_j^*$ this has the effect of `switching $G(\l)$ and
$G(-\l)$' so that case $(3,a)$ is the same as case $(1,4-a)$ with the sign of
$\l$ reversed. Thus we can
neglect cases $(3,1), (3,2), (3,3)$. Similarly, case $(2,3)$ is the same as 
$(2,1)$. 
Thus we have only to consider cases $(1,1)$, $(1,2)$, $(2,1)$, $(2,2)$ and $(1,3)$. 
Also, since $\sp(\l)$ is formally self-adjoint in the sense that
$$
\sp(\l)(z',z'') = \ol{ \sp(\l)(z'',z') },
$$
we may restrict attention to a neighbourhood of lb and the interior of bf. 

First consider case $(2,2)$. In this case, the composition is an integral of
the form
\begin{equation}\begin{gathered}
\int e^{i\phi/x} {x'}^{- k/2} a'(y', y, v, x') {x''}^{- k'/2}
a''(y, y'', w, x'') \, dv \, dw \, dy , \\
\phi(y', y'', \theta, v, w, y) = 
\cos \theta \,\phi'(y', y, v) - \sin \theta \,\phi''(y, y'', w) 
\end{gathered}
\label{22}
\end{equation}
where $v \in \RR^k, w \in \RR^{k'}$ and $\phi'$ and $\phi''$
 both parametrize $G(\l)$ in region two. (We will omit the $s\Phi$-half-density
factors for notational simplicity.)
The integral is rapidly decreasing in $x$ unless $d_{y}\phi = 0$. This implies
that 
$$
\sigma = \tan \theta = \frac{|d_y \phi' |}{|d_y \phi'' |},
$$
but in case $(2,2)$, $d_y \phi_i = \mu_i$ is never zero, so the integral is
rapidly decreasing at lb and rb. Thus we only need to analyze \eqref{22} in the
interior of bf.

Let us show that $\phi$ is a nondegenerate phase function in the interior
of bf, with phase variables $(v, w, y)$. Consider the phase function
$\phi'$. With $y$ held fixed, $\phi'$ parametrizes $G_y$, the Legendrian
defined in \eqref{Gy0}. (It is a non-degenerate parametrization since
the coordinate functions $y_i$ are linearly independent on $G(\l)$, hence 
on $C_{\phi}$, $d(y_i)$ and $d_{(y',y,v)}\pa_{v_j}\phi'$ are independent. That
is equivalent to independence of 
$d_{(y',v)}\pa_{v_j} \phi'$ when $d_{v}\phi'
= 0$, which is to say that $\phi$ with $y$ held fixed parametrizes $G_y$
nondegenerately.) Note that on $G_y$, the rank of the $(n-1)$ differentials
$d(\mu_i)$ is everywhere at least $(n-2)$, by \eqref{eq:def-Gp}. Therefore,
when $d_{v}\phi' = 0$, the rank of the $k+n-1$ differentials
$d_{y',v} \pa_{y_i} \phi'$ and $d_{y',v} \pa_{v_j} \phi'$ is at least
$k+n-2$. Next, note that 
$$ 
d_{\theta} \pa_y \phi  = -\sin \theta \frac{ \pa \phi'}{\pa y} - \cos \theta
\frac{\pa \phi''}{\pa y},
$$
and when $d_y \phi = 0$, we have $\cos \theta \, d_{y} \phi' = \sin \theta \,
d_y \phi''$, so this is
$$
d_\theta \phi = -\frac1{\cos \theta} (\sin^2 \theta + \cos^2 \theta) d_y \phi''
\neq 0.
$$
Finally, by the same reasoning as above, the differentials $d_{(y'',w)} \pa_{w_k}
\phi'' $ are independent when $d_w \phi'' = 0$. Putting these facts together
we have shown that $\phi$ is a nondegenerate phase function. 

It is not hard to see that the Legendrian parametrized by $\phi$ is $L$. Setting
$d_v\phi = d_w \phi = d_y \phi = 0$ we see that $(\sigma, y',y'',\tau',
\tau'', \mu', \mu'')$ is on the Legendrian only if there are points
$(y,y',\tau',\mu_1,\mu')$ and $(y,y'',\tau'',\mu_2,\mu'')$ in $G(\l)$ with 
$y',y''$ and $y$ on the same geodesic, with $y$ not the middle point, with
$\sigma = |\mu_1|/|\mu_2|$. With $|\mu_1| = \sin s'$ and $|\mu_2| = \sin s$
this is the same as \eqref{eq:sp-1}. Thus \eqref{22} is in $\Isf^m (\XXb, L)$
with $m = -\half$ by \eqref{eq:dist-1}. 

Next consider case $(1,1)$. In fact we can treat $(1,1)$ and $(1,2)$ simultaneously,
by considering the integral of the form 
\begin{equation}
x^{-k'} \int e^{i \big( -\cos \theta \cos d(y',y) - \sin\theta \phi''(y,y'',w) \big)/x}
a'(y',y',x') a''(y,y'',x'',w) \ dy \ dw
\label{12}
\end{equation}
where $w \in \RR^{k'}$ and
$a'$ is supported where $y'$ is close to $y$.
In the interior of bf the phase function $\phi$ is non-degenerate since the differentials
$d_{y'} \pa_{y_i} \phi$ are
linearly independent when $y'$ is close to $y$, and $d_{(y'',v)} \pa_{v_j} \phi''$ are
linearly independent as above. The same reasoning as above shows
that $\phi$ parametrizes $L$ in the interior of bf. 
Near lb the phase is of the form $-\cos d(y',y) - \theta \phi'' + O(\theta^2)$, which
has the form \eqref{eq:fib-ph-9} with $\phi_1 = -\cos d(y',y)$ and $\phi_2 = -\phi''$.
It is a non-degenerate parametrization of $L$ near the corner
since $d_{y'} \pa_{y_i} \phi_1$ and $d_{(y'',v)} 
\phi_2$ are independent. Since $\phi_1$ parametrizes $\Gs(-\l)$ and
for fixed $y = y'$, $\phi_2$ parametrizes $G_{y'}(-\l)$, we see that \eqref{12} is
in $\Isf^{-\half, \frac{n-1}{2}} (\XXb, \Lt)$. 

Near the interior of lb, we have an integral expression
\begin{equation}\label{eq:26a}
\int e^{-i\l \cos d(y',y)/x'} a'(y',y,x') a''(y,z'') dy.
\end{equation}
The second factor is a smooth function in this region since $z''$ lies in the
interior of $X$. 
The phase is nondegenerate at $y$, so performing stationary phase we get
\begin{equation}\label{eq:26b}
e^{-i\l /x'} {x'}^{(n-1)/2} c(x',y',z'') 
\end{equation}
which is an expression of the form \eqref{eq:dist-2} with $r = (n-1)/2$.

The case $(2,1)$ is easily treated, because in this case we have an integral
\begin{equation}
x^{-k} \int e^{i \big( \cos \theta \phi'(y',y,v) +  \sin\theta \cos d(y,y'') \big)/x}
a'(y',y',\theta,x') a''(y,y'',\theta,x'') \ dy
\end{equation}
where $|d_y \phi'| \neq 0$ whenever $d_v \phi' = 0$. In this case, we
get nontrivial contributions when $d_y \phi = d_v \phi = 0$, which implies that
$$
\sigma = \frac{|\sin d(y',y)|}{|d_y \phi'|} \text{ and } d_v \phi' = 0
$$
there. But this means that $\sigma \geq c$, so we only get contributions at the
interior of bf (since we are neglecting a neighbourhood of rb), and in this region
the argument is similar to that of case $(2,2)$. 

Finally we have case $(1,3)$. There the composite operator takes the form
\begin{equation}\begin{split}
\int_0^\infty ds \int &e^{i\lambda \big(-\cos\theta \cos d(y',y) - \sin\theta (1 + s \psi(y,y'',s,w))
\big) /x} \\
& a'(y',y,x') a''(y,y'',w,s,x''/s) 
\big( \frac{x''}{s} \big)^{-(k'+1)/2} s^{(n-3)/2}\,dw\,dy
\end{split}\end{equation}
where $a''$ is supported where $s$ is small. 
When $s=0$, the Hessian of the phase function in $y$ is nondegenerate, with a
critical point at $y=y'$, so for small $s$ we get a critical point
$y = f(y',y'',s, w)$. When $\theta = 0$ the critical point is also $y' = y$ so
the critical value takes the form $y = y' + s \theta g(y',y'',\theta,s,w)$. 
Using the stationary phase lemma
we may rewrite the integral as
\begin{equation}\begin{split}
\int_0^\infty ds 
\int &e^{i\lambda( -\cos\theta - \sin\theta + s \theta \td \psi(y',y'',\theta,s, w))/x}
a(y',y'',\theta,s,x/s,w) \\
&\big(\frac{x''}{s}\big)^{(n-1)/2-(k'+1)/2} s^{(2n-3)/2} \theta^{(n-1)/2} \,dy\,dw
\end{split}\end{equation}
where $\td \psi(y',y'',0,0,w) = \psi(y',y'',0,0,w)$. This is a nondegenerate
parametrization of $(L, \Ls(\l))$ both near the interior of bf and near the
corner $\lb \cap \bfc$ with orders $-\half$ at $L$, $(n-2)/2$ at $\Ls$ and
$(n-1)/2$ at lb, so the integral is in $\Isf^{-\half, \frac{n-2}{2}, \frac{n-1}{2}}
(\XXb, \Lt)$. The argument in the interior of lb is the same as in case $(1,2)$. 
This completes the proof of
the theorem. 
\end{proof}

\section{The resolvent}\label{sec:resolvent}

In this section we identify the structure of the resolvent $R(\l^2 \pm i0)$ on the
real axis. The resolvent is
a little more complicated than the spectral projection $\sp(\l)$ since it
has interior singularities along the diagonal. However we can cut the resolvent
into three pieces, each of which has a structure which has been described. One
piece lies in the scattering calculus; this takes care of the interior
singularities. A second piece is an intersecting Legendrian distribution
microsupported in the interior of bf. The third
piece is microlocally identical to the spectral projection in the region
$\tau < 0$ (for $R(\l^2 + i0)$) or $\tau > 0$ (for $R(\l^2 - i0)$). See 
section~\ref{sec:example} for a discussion of the resolvent on $\RR^n$. 

In preparation
for the description of the second piece, we describe the intersection properties
of the two Legendrians $L$ and $\sci N^*\diag_\bl$, the boundary of the closure
of $N^* \diag_\bl$ in $\sct \Mt$:
$$
\sci N^*\diag_\bl = \{ \sigma = 1, y'=y'', \mu' = -\mu'', \tau = 0 \}.
$$
The following proposition shows that $L$ and $\sci N^*\diag_\bl$ have the
correct geometry for an intersecting Legendrian distribution. 

\begin{prop} The intersection $L \cap \sci N^*\diag_\bl$ is codimension one in
$L$. In fact, the function $\tau$, restricted to $L$, has non-vanishing
differential at $\tau = 0$, and $L \cap \{ \tau = 0 \}$ is precisely
equal to $L \cap \sci N^*\diag_\bl$. Thus $L$ is divided by $\sci N^*\diag_\bl$
into two pieces, $L_+ = \{ \tau < 0 \}$ and $L_- = \{ \tau > 0 \}$. 
\end{prop}

\begin{proof} Clearly $L \cap \sci N^*\diag_\bl \subset L \cap \{ \tau = 0 \}$, so
we prove the reverse inclusion. Setting $\tau = 0$ in \eqref{eq:sp-1} and
using \eqref{relnship} we see
that 
$$
\sigma = \frac{\cos s'}{\cos s} = \frac{\sin s'}{\sin s},
$$
so $s=s'$ on the first piece of $L \cap \{ \tau = 0 \}$ in \eqref{eq:sp-1} and then
we have $\sigma = 1$, $y=y'$, $\mu' = -\mu''$. It is easy to check on the second and
third pieces using \eqref{relnship}. 

Since $L$ is contained in the characteristic variety of $\Delta_{\gt} - 2\l^2$ 
(see \eqref{Lap-gt})
which is $\{ \tau^2 + |\mut'|^2 + |\mut''|^2 + \eta^2 = 2\l^2 \}$, it is tangent to
the bicharacteristic flow of $\Delta_{\gt}$. Under this flow $\tau$ evolves by
$\dot \tau = - 2(|\mut'|^2 + |\mut''|^2 + \eta^2)$, so on $\Sigma(\Delta_{\gt} - 2\l^2)
\cap \{ \tau = 0 \}$, $\dot \tau = -4\l^2$ and hence $d\tau \neq 0$.
\end{proof}

\begin{thm}\label{thm:resolvent}
The boundary value $R(\lambda^2+ i0)$ of the resolvent of $H$ on the
positive real axis is of the form $R_1+R_2 + R_3$, 
where $R_1$ is a scattering pseudodifferential operator of order $(-2,0)$,
$R_2$ is an intersecting Legendrian distribution of order $-\frac1{2}$ associated to 
$(\sci N^*\diag_\bl, L_+)$ 
and $R_3$ is a Legendrian distribution associated to the conic pair
$(L_+, \Ls(\l))$, which is microlocally identical to $\sp(\l)$ for $\tau < 0$,
and in particular has the same orders as $\sp(\l)$ (see
Theorem~\ref{thm:sp-proj}). 
\end{thm}

\begin{proof}
We restrict our attention to $R(\lambda^2+i0)$ for the sake of definiteness. Similar
results hold for $R(\l^2 - i0)$ with some changes of signs; in particular, $L_+$
changes to $L_-$ in the statement above. 

We start by defining $R_1$.
Let $\psi\in\Cinf_c(\Real)$ be identically $1$ near $\lambda^2$ and
supported in a small neighborhood of $\lambda^2$, and let
\begin{equation}
R_1=(\Id-\psi(H))R(\lambda^2+i0),\ R' = \psi(H)R(\lambda^2+i0).
\end{equation}
Then $R_1$ is the distributional limit of $(1-\psi(H)) R(\l^2 + i\ep) =\psit_\ep(H)$ 
where $\psit_\ep (t)=(1-\psi(t))(t-\lambda^2 - i\ep)^{-1}$, so
$\psit_\ep\in S^{-1}(\Real)$ is a `classical' symbol for all $\ep \geq 0$,
as
$\lambda^2\nin\supp(1-\psi)$. Thus, by the symbolic functional calculus
of \cite{Hassell-Vasy:Symbolic} $\psit_\ep(H)\in\Psisc^{-2,0}(X)$ for
all $\ep \geq 0$; this means that its
kernel is conormal to $\diag_\scl$ on $\XXsc$, polyhomogeneous
at $\ffsc$, and vanishes to infinite order
on all other boundary hypersurfaces. Since $\sup |\psit_\ep - \psit_0 | \to 0$,
we have $\psit_\ep(H) \to \psit_0(H)$ as bounded operators on $L^2(X)$, hence
certainly as distributions on $\XXb$. Hence, $R_1 = \psit_0(H)$ is a scattering
pseudodifferential operator. Note that if $\psi_1$
is another cutoff function with the same properties and $\psit_1$ is defined
similarly to $\psit$ then $\psit-\psit_1\in\Cinf_c(\Real)$, so
$\psit(H)-\psit_1(H)\in\Psisc^{-\infty,0}(X)$ by the symbolic
functional calculus, i.e.\ its kernel is polyhomogeneous on $\XXsc$ (has
no singularity at $\diag_\scl$) and vanishes to infinite order an all other
boundary hypersurfaces. Thus, the splitting $R_1+R'$ is in fact well
defined modulo such terms, and it separates the diagonal singularity (which is
really a high energy phenomenon) represented by $R_1$
and the singularity due to the existence
of generalized eigenfunctions of energy $\lambda^2$,
represented by $R'$.

To analyze $R'$, note that the kernel of $R'$ is the distributional limit of
$\psi(H)(H - \l^2 - i\ep)^{-1}$ which is a function of $H$ supported on the positive
real axis; hence we can use our knowledge of the spectral projections for $\l > 0$
to compute this term as
$$
\int_0^\infty \psi(\rho) (\rho^2 - \l^2 - i\ep)^{-1} \sp(\rho) 2\rho \, d\rho.
$$
The term $(\rho^2 - \l^2 - i\ep)^{-1}$ can be expressed as
\begin{equation}
(\rho^2 - \l^2 - i\ep)^{-1} = \frac{i}{\rho + \l} \int_0^\infty 
e^{-i(\rho - \l)\ol t} e^{-\ep \ol t} d\ol t
\end{equation}
(presently we will change variable to $t = \ol t/x$). We have shown in the proof
of Theorem~\ref{thm:sp-proj} that the kernel of $\sp(\l)$ is a sum of expressions
of the form
\begin{equation}
x^q \int e^{i\l\phi/x} a
\label{a}
\end{equation}
where $\phi$ parametrizes some part of $(L(1), \Ls(1))$ and $a$ is smooth on $\XXb$
and in the phase variables. We need to prove smoothness of $a$ in $\l$ as
well; for ease of exposition this is left until the end of the section. We may suppose
that we have broken up $\sp(\l)$ in such a way that on each local expression,
$\phi$ is either strictly less than zero, strictly greater than zero or takes
values in $(-\delta, \delta)$ for some small $\delta$. Thus to understand $R'$ we
must understand the expressions
\begin{equation}
\int \int_0^\infty dt \ \frac{ 2i\rho}{x(\rho + \l)}  \psi(\rho) 
e^{-i(\rho - \l)t/x} e^{-\ep t/x}
x^q e^{i\rho \phi(y',y'',\theta, v)/x} a(y',y'',\theta, x',v,\rho) \,d\rho \, dv
\label{R'}
\end{equation}
as $\ep \to 0$. (Note that $\phi$ will also depend on an $s$ variable if it is
parametrizing $(L, \Ls)$ near at a conic point; this is not indicated in notation
for simplicity.) We define the phase function $\Phi = -(\rho - \l)t + 
\rho\phi$. 

It is inconvenient that \eqref{R'} is not compactly supported in $t$, so first
we insert $\chi(t) + (1 - \chi(t))$ into the integral, where $\chi$ is identically
equal to one on $[0, 2]$ and is compactly supported. The phase function $\Phi$
is stationary in $\sigma$ when $t = \phi$ which is always less
than $2$. We may also produce any negative power $t^{-k}$ in the integrand by
writing $e^{-i(\rho - \l)t/x} = -ixt^{-1} \pa_\rho e^{-i(\rho - \l)t/x}$ and repeatedly
integrating the $\rho$ derivative by parts. Thus, by stationary phase, when
$(1 - \chi(t))$ is inserted the integral is smooth on $\XXb$ and
rapidly decreasing as $x \to 0$, uniformly
as $\ep \to 0$. Hence we may ignore this term and insert the factor $\chi(t)$ into
\eqref{R'}. 

Let us consider \eqref{R'} at $\ep = 0$ with the term $\chi(t)$ inserted, that is,
\begin{equation}
\int \int_0^\infty \frac{ 2i\rho}{x(\rho + \l)}
dt \ \psi(\rho) e^{-i(\rho - \l)t/x} \chi(t)
x^q e^{i\rho \phi(y',y'',\theta, v)/x} a(y',y'',\theta, x',v,\rho) d\rho dv.
\label{Phi}
\end{equation}
when the phase $\phi$ is negative. In this
case there are no critical points of $\Phi$ since $d_\rho \Phi = 0$ requires
$t = \phi$ and $t$ is nonnegative. Thus in this case, again, the
expression is smooth and rapidly decreasing as $x \to 0$ uniformly in $\ep$. 

When the phase is positive, we get critical points when $\phi$ parametrizes
$(L, \Ls)$, $\rho = \l$ and $t = \phi$. We can ignore the fact
that the range of integration of $t$ has a boundary since $d_\rho \Phi \neq 0$
there. The phase function $\Phi$ is nondegenerate in $(t, \rho)$ and exactly
quadratic, and moreover the rest of the expression is independent of $t$. It
follows that the first term of the stationary phase expansion is accurate to
all orders, so that up to a smooth function which is rapidly decreasing in $x$,
\eqref{Phi} is equal to 
$$
(2\pi i) \int x^q e^{i\l\phi(y',y'',\theta, v)/x} a(y',y'',\theta, x',v,\l) dv
$$
which is $2\pi i$ times the corresponding piece of $\sp(\l)$. We take
$R_3$ to be the sum of the contributions from all terms with $\phi > 0$ in our
decomposition. Since $\tau < 0$
when $\phi > 0$, we see that $R_3 \in \Isf^{\calK,\frac{n-2}{2}}
(\XXb, (L_+(\l), \Ls(\l)), \Osfh)$ 
with $\calK(\lb)$ as in Theorem~\ref{thm:sp-proj}, and is (up to a numerical factor)
microlocally identical to $\sp(\l)$ there. 

When the phase changes sign, then we have an intersecting Legendrian distribution
associated to $(L_+, \sci N^*\diag_\bl)$. To see this, note that when $t > 0$, 
$\pa_t \Phi = 0$ implies $\rho = \l$, $\pa_\rho \Phi = 0$ implies $t = \phi$
and $d_v \Phi = 0$ implies that $\phi$ parametrizes $L(1)$. Thus we get the
Legendrian $L(\l)$ for $t>0$, subject to the restriction $\tau = -\phi < 0$ which
gives us $L_+$. The phase is non-degenerate since $\phi$ is nondegenerate and the
Hessian of $\Phi$ in $(t, \rho)$ is nondegenerate. In addition, $dt \neq 0$ at
$t = 0$ on $C_\Phi$. If $\Phi_0 = \rho \phi$ denotes $\Phi$ restricted to 
$t=0$, then $d_\rho 
\Phi_0 = 0$ implies that $\phi = 0$ and $d_v \Phi_0 = 0$ implies that $\phi$ parametrizes
$L(1)$. The Legendrian parametrized is then
$$
\{ (y',y'',\theta,\tau,\mu',\mu'',\eta) \mid \tau = 0, (y',y'',\theta,0,
\rho\mu',\rho\mu'',\rho\eta) \in L(1) \}
$$
which is $\sci N^*\diag_\bl$. Thus in this case \eqref{R'} is an intersecting Legendrian
distribution when $\ep = 0$. We define $R_2$ to the sum of all the contributions when
$\phi$ changes sign. It is easy to see that as $\ep \to 0$ we get
convergence in distributions to $R_2$, since by multiplying by a suitable power of $x$, we
have convergence in $L^1$ by the Dominated Convergence Theorem; therefore the
expression converges in some weighted $L^2$ space. 

We have now written $R(\l^2 + i0)$ as a sum $R_1 + R_2 + R_3 + \td R$, where $R_i$ have
the required properties and $\td R$ is a smooth kernel vanishing with all derivatives
at the boundary, hence trivially fits any of the descriptions of the $R_i$. 
Thus, together with the proof of smoothness of \eqref{a} in $\l$, the 
proof of the theorem is complete.
\end{proof}
 
Finally we show smoothness of \eqref{a} in $\l$. To do this it is enough to show
smoothness of the symbol of $P(\l)$ in $\l$. The Poisson operator $P(\l)$ is
constructed in \cite{RBMZw} from a parametrix $K(\l)$, which is explicit and
easily seen to be smooth in $\l$, plus an correction term,
$U(\l) = R(\l^2 + i0) (\Lap - \l^2 ) K$. Here $K$ is constructed so that
$(\Lap - \l^2) K(\l)$ is in $\dCinf (X \times \pa X)$. Therefore, the following
lemma completes the proof. 

\begin{lemma} The resolvent $R(\l^2 + i0)$ acting on $f \in \dCinf(X)$ is of
the form
$$
x^{(n-1)/2} e^{i\l/x} a(x, y, \l),
$$
where $a$ is smooth in $x, y$ and $\l$.
\end{lemma}

\begin{proof} Smoothness in $x$ and $y$ is proved in \cite{RBMSpec}, so it
remains to prove smoothness in $\l$. Let $\sigma = \l^2$. 
One can show, as in \cite{RBMSpec}, proof
of Proposition 14, that for $|\re \sigma| \geq \epsilon$, $R(\sigma)^2 f$ has a
limit in $x^{-3/2-\epsilon} L^2$ as $\sigma$ approaches the real axis from
$\im \sigma > 0$, and the
limit is continuous in $\sigma$. (Indeed, this is completely
analogous to Jensen's proof of the corresponding statement
in the Euclidean setting, see \cite{Jensen:Propagation}.)
Moreover, it is the unique solution to
$(\Lap - \sigma)^2 u = f$ satisfying the radiation condition $\tau < 0$
on $\WFsc(u)$. Let us denote the limit by $R(\sigma + i0)^k f$. 
Since $\pa_\sigma  R(\sigma) = -R(\sigma)^2$, we see that both 
$R(\sigma + i\epsilon)f$ and $\pa_\sigma R(\sigma + i\ep)
f$ have limits in $x^{-3/2 - \ep} L^2$ as $\ep \to 0$. Hence, we have
$$
\pa_\sigma R(\sigma + i0) f = -R(\sigma + i0)^2 f.
$$
Since $a = x^{-(n-1)/2} e^{-i\l/x} R(\l^2 + i0) f$, we see that the $\l$ derivative
of $a$ exists in $x^{-3/2 - \ep} L^2$. 

If we calculate 
$$
\frac{d}{d\l} (\Lap - \l^2) (e^{i\l/x} a) = 
(\Lap - \l^2) \big( e^{i\l/x} \frac{\pa a}{\pa \l} \big)
+ (\Lap - \l^2) \big( \frac{i}{x} e^{i\l/x} a \big) - 2\l (e^{-i\l/x} a) = 0,
$$
we see that
\begin{equation}\begin{gathered}
(\Lap - \l^2) \big( e^{i\l/x} \frac{\pa a}{\pa \l} \big) = 
(\Lap - \l^2) \big( \frac{i}{x} e^{i\l/x} a \big) - 2\l (e^{-i\l/x} a) \\
 = \frac{i}{x} f - 2i(i\l) e^{i\l/x} a + 2i e^{i\l/x} (x^2 \pa_x) a - (n-1) ix e^{i\l/x}
a - 2\l (e^{-i\l/x} a)\\
= \frac{i}{x} f + 2ix e^{i\l/x} (x \pa_x - \frac{n-1}{2}) a \\
= e^{i\l/x} x^{(n+3)/2} b(x, y, \l) + \frac{i}{x} f, \quad b \in \Cinf(X).
\end{gathered}\end{equation}
We can find a solution to 
$$
(\Lap - \l^2) v_1 = e^{i\l/x} x^{(n+3)/2} b \quad \text{mod } \dCinf
$$
of the form $v_1 = e^{i\l/x} 
a_1$, with $a_1 \in x^{(n+1)/2} \Cinf(X)$ and varying continuously with $\l$
(since $b$ depends on $a$ only through
taking $x$ and $y$ derivatives, not $\l$ derivatives). 
Therefore, 
\begin{equation}
(\Lap - \l^2)  \big( e^{i\l/x} (\frac{\pa a}{\pa \l} + a_1) \big) \in \dCinf(X),
\label{a-deriv} \end{equation}
with the right hand side continuous in $\l$.
It follows from Proposition 12 of \cite{RBMSpec}, and Lemma~\ref{continuity}
that $\pa_\l a$ is a smooth function on $X$ which is continuous in $\l$. 

This implies that actually the right hand side of \eqref{a-deriv} is $\calC^1$ in $\l$.
Thus, we can repeat the argument with $e^{i\l/x} (\pa_\l a + a_1)$ in place of
$e^{i\l/x} a$. Inductively we see that $a$ is $\calC^k$ in $\l$ for any $k$, proving
the lemma.
\end{proof}

\begin{rem}
In this paper we analyzed the kernel of the resolvent at the real axis
by using the
Poisson operator. Conversely, the resolvent
determines, in a simple way, 
the other analytic operators associated to $H$. We have
already seen that the generalized spectral projection is
$\sp(\l)=(R(\l^2+i0)-R(\l^2-i0))/(2\pi i)$, i.e.\ essentially the
difference of the limits of the resolvent taken from either side of the
real axis. Next, the Poisson operator $P(\pm \l)$
is a multiple of the principal symbol
of the kernel of $R(\l^2 \pm i0)$ 
at either $\lb$ or $\rb$. (We did not define the
principal symbol for Legendre distributions on manifolds with corners
because we did not need this notion,
but for oscillatory functions as in \eqref{eq:fib-osc},
it is, up to a constant multiple,
the restriction of the amplitude $a$ to $x=0$. Note that the kernel
of $R(\l^2 \pm i0)$ 
is indeed such an oscillatory function at lb and rb since $L^\sharp(\pm\l)$
are global sections of $\sfT_{\ffb} \XXb$.) Finally, $S(\l)$ is
a multiple of the principal
symbol of $R(\l^2 \pm i0)$ along $L^\sharp(\pm \l)$ in $\ffb$.

These statements follow from the proof of Theorem~\ref{thm:sp-proj}, together
with the fact that $R(\l^2 \pm i0)$ is microlocally identical to $\sp(\l)$ for
$\tau < 0$ and $\tau > 0$, respectively.
For example, in \eqref{eq:26a}, $a'(y,y,0)$ is a constant (it is the
initial data for the transport equation for the principal symbol of
$P(\l)$, see \cite{RBMZw}),
and $a''(y,z'')$ is a multiple of the kernel of $P(\l)^*$,
so the stationary phase argument yields \eqref{eq:26b} with $c(0,y',z'')$
being a multiple of the kernel of $P(\l)^*$. The statement about $S(\l)$
follows similarly, but it can also be deduced from the fact that its kernel
is the principal symbol of $P(\l)$ along $\Gs(\l)$.
\end{rem}

\section{Distorted Fourier Transform}\label{sec:distorted-FT}
Here we interpret some of our results in terms of a generalized `Distorted Fourier
transform'. Namely, we show that the Poisson operators determine two isometries $P_\pm$
from the absolutely continuous spectral subspace of $H$ to $L^2(\pa X \times \RR^+)$
(cf. \cite{Herbst-Skibsted:Free-channel} in the more complicated $N$-body case).
The coordinate on $\RR^+$ should be thought of as the square root of the energy, $\l$.
We will show that they 
determine the scattering matrix through equation \eqref{Sequation}, which is
analogous to 
a standard formula in scattering theory (see \cite{Taylor}, chapter 9, section 3,
or \cite{RS3}, section XI.4).

\begin{prop} The maps $P^*_\pm$, defined on $\Cinf_c(X) \cap H_{ac}(H)$ by 
\begin{equation}
(P^*_\pm u)(y,\l) = (2\pi)^{-1/2} (P(\pm \l)^* u)(y),
\end{equation}
map into $L^2(\pa X \times \RR^+)$ and extend to unitary maps from $H_{ac}(H)$
to $L^2(\pa X \times \RR^+)$, with adjoints 
\begin{equation}
(P_\pm f)(z) = (2\pi)^{-1/2} \int_0^\infty P(\pm\l)(f(\cdot, \l)) \, d\lambda .
\end{equation}
They intertwine the Hamiltonian with multiplication by $\l^2$:
$$
P^*_\pm H = \l^2 P^*_\pm .
$$
\end{prop}

\begin{proof}
By Lemma~\ref{continuity}, for $0 < \ep < a < b < \infty$, 
and a dense set of $f \in L^2$, namely,
$f \in x^{1/2 + \ep} L^2$, the function $\ang{f}{R(\l) f}$ is continuous on 
$|\re \l| \geq \ep$ and $\pm \im \l \geq 0$. It follows that for such $f$, and
$\l \in \RR$,
\begin{equation}\begin{gathered}
\lim_{\ep \to 0} \int_a^b \ang{f}{R(\l^2 + i\ep) f} - 
\ang{f}{R(\l^2 - i\ep) f} \, 2\l d\l \\
= \int_a^b i \ang{f}{P(\l) P(\l)^* f} \, d\l \quad \text{ by equation } 
\eqref{eq:intro-5} \\
= \int_a^b i \big| P(\l)^*f \big|^2 \, d\l.
\end{gathered}\end{equation}
Sending $a \to 0$ and $b \to \infty$, we see that
$$
\ang{E_{(0, \infty)}f}{f} = \frac1{2\pi} \int_0^\infty \big| P(\l)^*f \big|^2 \, d\l
$$
for a dense set of $f$. Since $H_{ac}(H) = E_{(0, \infty)}$, this shows that
$P^*_+$ extends to an isometry from $H_{ac}(H)$ to
$L^2(\pa X \times \RR^+)$. The same reasoning shows that $P^*_-$ also extends as an
isometry. It is easy to see that the adjoint of $P^*_\pm$ is $P_\pm$. 

The intertwining property follows from $(H - \l^2) P(\pm \l) = 0$ by multiplying
on the left and right by $P^*_\pm$. 
\end{proof}

\begin{prop}
If the operator $S$ is defined by 
$$
(Sf)(y,\l) = (S(\l) f(\cdot, \l)(y)),
$$
then we have
\begin{equation}
S = P_-^* P_+ .
\label{Sequation}
\end{equation}
\end{prop}

\begin{proof}
Multiplying both sides of the equation on the left by $P_-$, it suffices to
establish $P_- S = P_+$. Since each operator acts fibrewise in $\l$, it is
sufficient to establish $P(-\l) S(\l) = P(\l)$. Recall that for $v \in \Cinf(\pa X)$,
$P(\l) v$ is the solution to $(H - \l^2)u = 0$ with incoming boundary data $v$. Denote
by $w$ the outgoing boundary data of $u$, ie, $w = S(\l) v$. Then $P(-\l)w$ is
also equal to $u$. But this proves $P(-\l) S(\l) = P(\l)$, so the proof is complete.
\end{proof}

\begin{rem} For $H = \Lap + V$ on $\RR^n$, our normalization for the two
distorted Fourier transforms is different from the standard normalization,
which would have the kernels of the two being $c
\cdot e^{-i\l z \cdot \theta}$ plus
incoming/outgoing correction term. This makes the scattering matrix for the free
Laplacian equal to the identity operator. Here, the two distorted Fourier transforms
have the form $c \cdot
e^{\pm i\l z \cdot \theta}$ plus incoming/outgoing correction, and the
scattering matrix for the free Laplacian is $i^{-n+1} R$, where $R$ is the antipodal
map on the $(n-1)$-sphere. In \cite{Hassell:S-matrices} 
these two operators are called the analytic
and geometric scattering matrices, respectively. Notice that only the geometric
scattering matrix has an analogue for a general scattering metric.
\end{rem}

\section{Wavefront set relation}\label{sec:wavefront-reln}

Here we will use the microlocal information about the resolvent from 
Theorem~\ref{thm:resolvent} to find a
bound on the scattering wavefront set of $R(\l + i0)f$ in terms of the scattering
wavefront set of $f$. The result is not new, since Melrose proved it in \cite{RBMSpec}
using positive commutator estimates. However, it is quite illuminating to see how
the result arises from the structure of the resolvent kernel itself. The fundamental
tool is the pushforward theorem of \cite{Hassell:S-matrices}. 

\begin{thm} If $f \in \Cinf(X^o) \cap \dist(X)$, and $\WFsc(f)$ does not
intersect $\Gs(-\l)$, then $\WFsc(R(\l^2 + i0)f)$
is contained in the union of $\WFsc(f)$, $\Gs(\l)$ and bicharacteristic rays of
$H - \l^2$ that start at $\WFsc(f) \cap \Sigma(H - \l^2)$ and propagate in the
direction of decreasing $\tau$. That is,
\begin{equation}\begin{gathered}
\WFsc( R(\l^2 + i0)f ) \subset \WFsc(f) \cup \Gs(\l) \cup \\
\{ (y, \tau, \mu) \mid \exists \ (y', \tau', \mu') \in \WFsc(f) \text{ with }
\tau' = \l \cos s, \ \mu' = \l \sin s \hat \mu', \ |\hat \mu'| = 1, \\
\text{ and } (y, \mu) = \l \sin(s+t) \exp(tH_{\half h})(y', \hat \mu'), 
\tau = \l \cos(s+t) \text{ for some } s \text{ and } t > 0. \}
\end{gathered}\label{eq:scwf-res-f}\end{equation}
\end{thm}

\begin{proof} We write $R(\l^2 + i0) = R_1 + R'$ as in the proof of 
Theorem~\ref{thm:resolvent}. Then $R_1$ is a scattering pseudodifferential operator,
so $\WFsc(R_1 f) \subset \WFsc(f)$. Thus we may restrict attention to $R' f$. 
The action of $R'$ on $f$ may be represented
by pulling $f$ up to $\XXb$ from the left factor, multiplying by the kernel of
$R'$ and then pushing forward to the right factor of $X$ (see figure). 
We will apply the pushforward
theorem to this pushforward to obtain scattering wavefront set bounds. 

Recall that the pushforward theorem from \cite{Hassell:S-matrices} 
applies in the situation when we have
a map $\pi: \td X \to X$ between manifolds with boundary, which is locally of the form
$(x, y, z) \mapsto (x, y)$ in local coordinates where $x$ is a boundary defining
function for $X$, and $\pi^* x$ is one for $\td X$. Corresponding coordinates on the
two scattering cotangent bundles over the boundary are then $(y, z, \tau, \mu, \zeta)$
and $(y, \tau, \mu)$. Then the pushforward theorem states that for densities $u \in
\Cinf(\td X^o) \cap \dist(\td X)$,
\begin{equation}
\WFsc(\pi_* u) \subset \{ (y, \tau, \mu) \mid \exists \ (y, z, \tau, \mu, 0) \in
\WFsc(u) \}.
\label{pushforward}\end{equation}
Unfortunately there is no pushforward theorem currently
available on manifolds with corners.
Since $\XXb$ has codimension 2 corners, we cannot apply the theorem directly, but
must cut off the resolvent kernel away from the corners and treat the part
near the corners by a different argument. It turns out that nothing special happens
at the corners and the bound on the wavefront set is just what one would expect
from a consideration of the boundary contributions alone. 

To localize the kernel of $R'$, choose a partition of unity
$1 = \chi_1 + \chi_2 + \chi_3 + \chi_4$, where $\chi_1$ is zero in a neighbourhood
of lb and rb, $\chi_2$ is supported near $\lf \cap \bfc$, $\chi_3$ is zero on
bf and rb and $\chi_4$ is supported close to rb. Correspondingly decompose
$$
R' = R'_1 + R'_2 + R'_3 + R'_4.
$$
We investigate each piece separately. \begin{figure}\centering
\epsfig{file=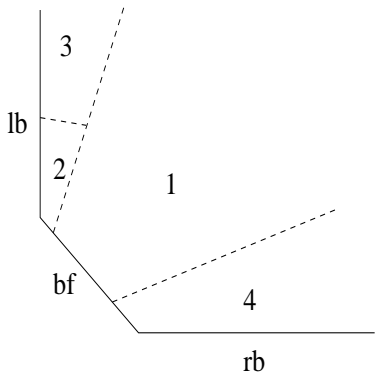,width=5cm,height=5cm}
\caption{The approximate supports of $R'_i$ on $\XXb$}
\end{figure}

The easiest piece to deal with is $R'_3 f$. The product of the two distributions
on $\XXb$ is localized near lf, and $f$ is smooth there. We know that $R'_3$ is
a Legendrian associated to $\{ \tau' = -1, \mu' = 0 \} = \Gs(\l)$ and is smooth in
$z''$ (where $z''$ is a set of coordinates in the interior of the right factor of $X$). 
Thus if $\zeta''$ is a rescaled dual variable to $z''$, the wavefront set of
the product is contained in $\{ \tau' = -1, \mu' = 0, \zeta'' = 0 \}$ and by
\eqref{pushforward}, the result has scattering wavefront set in $\Gs(\l)$.

Next we treat $R'_1 f$. This is localized near the interior of bf. We will see
that the contribution
here is the part that propagates inside the characteristic
variety of $H - \l^2$. To apply \eqref{pushforward}, we must first choose
compatible coordinates on $\XXb$ and $X$. Let us take $x'$, $y'$, $y''$ and
$\sigma = x'/x''$ near the interior of bf, and let $\ol \tau'$, 
$\ol \mu'$, $\ol \mu''$ and $\ol \nu$
be the corresponding rescaled cotangent coordinates. They
are related to coordinates \eqref{eq:sp-a-10} by 
\begin{equation}
\ol \mu' = \mu', \quad \ol \mu'' = \sigma \mu'', \quad \ol \tau' = \tau' + \tau'' \sigma,
\quad \ol \nu = -\tau''.
\label{coordchange}\end{equation}
The scattering wavefront
set of the lift of $f$ to $\XXb$ is
$$
\{ (y', y''_0, \sigma, \ol \tau' = \sigma \tau''_0, 
\ol \mu'=0, \ol \mu'' = \sigma \mu''_0, \ol \nu = -\tau''_0) \mid (y''_0, \tau''_0, \mu''_0)
\in \WFsc(f) \}.
$$
The wavefront set of $R'_1$ is contained in the union of $\Nsc^*(\diag_\bl)$
and $L_+$. In our present coordinates this is the union of 
\begin{equation}\label{eq:Nsc-diag}
\{ (y'' = y', \sigma = 1, \tau' = 0, \ol \mu' = -\ol \mu'', \ol \nu = 0) \},
\end{equation}
corresponding to $\Nsc^*(\diag_\bl)$, and
\begin{equation}\begin{split}
&\{(y',y'',\sigma,\ol \tau',\ol\mu',\ol\mu'', \ol\nu):\ \exists(y,\muh)\in S^*\bX,
\ s,s'\in[0,\pi],
\ (\sin s)^2+(\sin s')^2>0,\\ 
&\quad s+ s' > 0 \ \text{ such that } \ \sigma=\frac{\sin s'}{\sin s},\quad
\ol\tau'=\lambda\cos s' + \sigma \l \cos s, \quad \ol\nu = -\l \cos s,\\
&\quad(y',\ol\mu')=\lambda\sin s'\exp(s'H_{\half h})(y,\muh), \quad
(y'',\ol\mu'')=-\sigma\lambda\sin s \exp(s H_{\half h})(y,\muh)\}\\
&\cup\{(y,y,\sigma, \l - \sigma \l ,0,0, \l):\ \sigma\in[1, \infty),\ y\in\bX\}\\
&\cup\{(y,y,\sigma, -\l +\sigma \l ,0,0, -\l):\ \sigma\in (0, 1],\ y\in\bX\}
\end{split}\end{equation}
by \eqref{eq:sp-1} and \eqref{coordchange} (the condition $s+s' > 0$, and the 
restrictions on $\sigma$ in the last two lines, gives the 
correct `half' of the Legendrian for $R(\l^2 + i0)$ since this is precisely the
part of \eqref{eq:sp-1} where $\tau < 0$). 
It is not hard to see that \eqref{eq:Nsc-diag}, associated to $\Nsc^*(\diag_\bl)$,
gives us $\WFsc(f)$ back again, so consider the second part.  The first piece 
can be re-expressed as
\begin{equation}\begin{split}
&\{(y',y'',\sigma,\ol \tau',\ol\mu',\ol\mu'', \ol\nu): \exists \ s,s'\in[0,\pi],
\ (\sin s)^2+(\sin s')^2>0,\ s+ s' > 0 \\
&\text{ such that } \ \sigma=\frac{\sin s'}{\sin s},
\quad \ol\tau'=\lambda\cos s' + \sigma \l \cos s, \quad \ol\nu = -\l \cos s,\\
&\quad \ol\mu' = \l \sin s' \hat{\eta'}, \ |\hat \eta'| = 1, \ 
\quad \ol\mu'' = \l \sigma \sin s \hat{\eta''}, \ |\hat \eta''| = 1, \\
&\quad (y', \eta') = \exp((s'-s)H_{\half h})(y'', \eta'') \}.
\end{split}\end{equation}
The scattering wavefront set of the product is in each fibre the setwise sum of the 
scattering wavefront sets of the factors, ie, 
\begin{equation}\begin{split}
&\{(y',y'',\sigma,\ol \tau',\ol\mu',\ol\mu'', \ol\nu): \exists \ s,s'\in[0,\pi],
\ (\sin s)^2+(\sin s')^2>0,\ s+ s' > 0 \ \text{ with }\\
&\quad\sigma=\frac{\sin s'}{\sin s}, \quad
\ol\tau'=\lambda\cos s' + \sigma \l \cos s + \sigma \tau''_0, \quad \ol\nu = -\l \cos s - 
\tau''_0 \\
&\quad \ol\mu' = \l \sin s' \hat{\eta'}, \ |\hat \eta'| = 1, \ 
\quad \ol\mu'' = \l \sigma \sin s \hat{\eta''} + \sigma \mu''_0, \ |\hat \eta''| = 1, \\
&\quad (y', \eta') = \exp((s'-s)H_{\half h})(y'', \eta'') \mid (y''_0, \tau''_0, \mu''_0)
\in \WFsc(f) \} \\
&\cup\{(y,y,\sigma, \l - \sigma \l + \sigma \tau''_0,\sigma \mu''_0, \l - \tau''_0):\ 
\sigma\in[1, \infty),\ y\in\bX\}\\
&\cup\{(y,y,\sigma, -\l +\sigma \l + \sigma \tau''_0,\sigma \mu''_0, -\l - \tau''_0):\ 
\sigma\in (0, 1],\ y\in\bX\}
\end{split}\end{equation}
By \eqref{pushforward}, we 
only get a contribution to the scattering wavefront set of the pushforward when 
$\ol\nu = \ol\mu'' = 0$, ie, when 
$$
\tau''_0 = -\l \cos s, \quad \mu''_0 = -\l \sin s \hat \eta'',
$$
ie when the corresponding point of $\WFsc(f)$ is in the characteristic variety, and then
we get points
\begin{equation}\begin{split}
&\{ (y', \tau', \mu') \mid \tau' = \l \cos s', (y', \mu') = \l \sin s' 
\exp((s'-s)H_{\half h})(y''_0, \hat \eta''), \ s' > s \} \\
& \quad \cup \{ (y, \l, 0) \mid (y, \l, 0) \in \WFsc(f)\} .
\end{split}\end{equation}
The first set is precisely the bicharacteristic ray from $(y''_0, \tau''_0, \mu''_0)$
in the direction in which $\tau$ decreases, and the second set is the intersection
of the scattering wavefront set of $f$ with $\Gs(\l)$. 
This shows that the contribution from the
interior of bf is contained in the set specified in the theorem. 

\

Next we consider the corner contributions. We use a rather crude method here to
bound the scattering wavefront set contributions near the corner; it is not optimal,
but as the support of $\chi_2$ shrinks to the corner, we recover the optimal result. 
First we treat $R'_2$. Notice that $\ol \tau' = -1$ at the corner, so, by taking
the support of $\chi_2$ sufficiently small, we may suppose that the phase function,
$\phi$, of the kernel of $R'_2$ takes values arbitrarily close to $1$. 

If we localize in the coordinate $y'$, then via a diffeomorphism that takes a
neighbourhood of the boundary of $X$ to a neighbourhood to a boundary point of
$\ol{\RR^n}$ (the radial compactification of $\RR^n$) we may suppose that
$R'_2f$ is a function on $\RR^n$. Then we may exploit the following characterization
of the scattering wavefront set on $\RR^n$ (see \cite{Hassell:S-matrices}):
\begin{equation}\begin{gathered}
\WFsc(u)^\complement = \{ (\hz_0, \zeta_0) \mid 
\exists \Cinf \text{ function } \phi(x, \hz) \\
\text{ with } \phi(0, \hz_0) \neq 0, \widehat{\phi u}(\zeta) \text{ smooth near }
\zeta = \zeta_0 \}.
\label{alt}
\end{gathered}\end{equation}
Here, $\zeta$ is related to the $(\tau, \mu)$ coordinates by $\tau = -\zeta \cdot \hz$
and $\mu = \zeta +(\zeta \cdot \hz)\hz$, ie, $-\tau$ and $\mu$ are the components of
$\zeta$ with respect to $\hz$. Thus, we need to investigate where 
\begin{equation}\begin{gathered}
\int e^{-i z \cdot \zeta} dz \int \frac{dx'' \, dy''}{{x''}^{n+1}}
\int dv \, dw \,
e^{i\l\phi(\hz', v, y'', \sigma, w)/x'} f(x'', y'') |z'|^{-k} \\
\chi_2(z, y'', \sigma) \td \chi(z) a(z, y'', \sigma, v, w) 
\end{gathered}
\end{equation}
is smooth in $\zeta$. Here $\td \chi$ is a function that localizes the integral
in $z$. We have multiplied by a large negative power of $|z'|$ to make the integral
convergent; this does not affect the scattering wavefront set. 

Note that $f \cdot \chi_2 \cdot \td \chi \cdot a$ has stable regularity under 
repeated application
of the vector fields $z_i \pa_{z_j}$, since $f$ is independent of $z$ and
$\chi_2, \td \chi, a$ are symbolic in $z$. Integrating
by parts, and expressing $z_i$ as $i\pa_{\zeta_i}$, we find that the integral 
has stable regularity under the application of
$$
(\zeta_j + \hz'_j \l \phi + \l \hz'_j \sigma \pa_{\sigma} \phi ) \pa_{\zeta_i}.
$$
But $\phi$ is close to $1$, and $\sigma$ is close to zero, so for some small
$\delta > 0$, if $|\zeta + \lambda \hz'| > \delta$ then the factor multiplying the vector
field is always nonzero, showing smoothness in $\zeta$. Thus,
$$
\WFsc(R'_2(f)) \subset \{ |\zeta + \l \hz'| \leq \delta \}.
$$
As the support of $\chi_2$
shrinks to the corner, $\delta$ shrinks to zero, so we see that the contribution
of the corner to the scattering wavefront set is contained in
$\{ \zeta' = -\l \hz \}$, which corresponds under \eqref{alt} to $\Gs(\l)$. 
This is just what one would expect by taking the result for the left boundary
or the b-face and assuming that it were valid uniformly to the corner. 

Finally we come to $R'_4$. To deal with this, we use a duality argument and exploit
the fact that the adjoint of $R'_4$ is $R'_2$ for $R(\l^2 - i0)$. We use the
characterization

\begin{lemma} Let $q \in K$ and $f \in \Cinf(X^o) \cap \dist(X)$. Then
$q \notin \WFsc(f)$ if there exists $B \in \Psisc^{0,0}(X)$, elliptic at $q$,
and a sequence $(f_n) \in \dCinf(X)$ converging to $f$ in $\dist(X)$ such that
$\langle u, Bf_n \rangle $ converges for all $u \in \dist(X)$. \end{lemma}

To prove the lemma, note that the condition on the sequence $B f_n$ implies that
$B f_n$ converges weakly in every weighted Sobolev space, hence by compact embedding
theorems, strongly in every weighted Sobolev space, and thus in $\dCinf(X)$. Also,
we know that $Bf_n \to Bf$ in $\dist(X)$, so this implies that $Bf \in \dCinf(X)$. 
Since $B$ is elliptic at $q$ this shows that $q \notin \WFsc(f)$. 

Now given $q$ that is not in the set \eqref{eq:scwf-res-f}, choose a self-adjoint
pseudodifferential
operator $B$ elliptic at $q$ and with $\WFsc'(B)$ disjoint from the set 
\eqref{eq:scwf-res-f}. Since \eqref{eq:scwf-res-f} is invariant under the forward
bicharacteristic flow of $H - \l^2$, one may assume that $\WFsc'(B)$ is invariant under the backward
bicharacteristic flow. Also, since by assumption $\WFsc(f)$ does not meet $\Gs(-\l)$,
we may assume that $B$ is elliptic on $\Gs(-\l)$. 
Choose a zeroth order self-adjoint pseudodifferential operator
$A$ with is microlocally trivial on \eqref{eq:scwf-res-f} and microlocally equal to
the identity on $\WFsc'(B)$. 

Since $Af \in \dCinf(X)$, it is possible to choose $f_n \in \dCinf(X)$ converging to 
$f$ in $\dist(X)$, such that $Af_n \to Af$ in $\dCinf(X)$. Consider the distributional
pairing
$$
\langle B R'_4 f_n, u \rangle
$$
for $ u \in \dist(X)$. This is equal to 
\begin{equation}\begin{gathered}
\langle B R'_4 A f_n, u \rangle + \langle B R'_4 (\Id - A) f_n, u \rangle \\
= \langle A f_n, {R'_4}^* B u \rangle + \langle f_n, (\Id - A) {R'_4}^* B u \rangle 
\end{gathered}\end{equation}
Note that ${R'_4}^*$ is defined on any distribution since the kernel is smooth in the
interior and vanishes at rb (hence there is no integral down to any boundary hypersurface
which would lead to convergence problems). Since $Af_n$ converges in $\dCinf(X)$, the
first term converges. As for the second term, the scattering wavefront set of $Bu$ is
contained in $\WFsc'(B)$ which is invariant under the {\em backward} 
bicharacteristic flow and
contains $\Gs(-\l)$, so taking into account that ${R'_4}^*$
is essentially the term $R'_2$ for the incoming resolvent
$R(\l^2 - i0)$, the above results show that the scattering
wavefront set of ${R'_4}^*(Bu)$ is contained in $\WFsc'(B)$ and is disjoint from the
essential support of $\Id - A$. Thus $(\Id - A) {R'_4}^* Bu \in \dCinf(X)$, and so the
second term converges. This establishes that $q \notin \WFsc(f)$ by the Lemma above, 
and completes the proof of the proposition. 
\end{proof}

\begin{rem}
Although the proof is quite long, the important point is that the main, propagative
part of the scattering wavefront set comes from the contribution from bf and 
is obtained by a routine computation based on the Pushforward theorem.
\end{rem}

\begin{rem} The condition that the scattering wavefront set of $f$ not meet $\Gs(-\l)$
cannot be removed. For example, on $\RR^n$, the resolvent $R(\l^2 + i0)$ is conjugated
via Fourier transform to the operation of multiplication by 
$(|\zeta|^2 - \l^2 - i0)^{-1}$.
This is undefined on certain functions, for example $(|\zeta|^2 - \l^2 + i0)^{-1}$, 
whose inverse Fourier transforms $f$ have scattering wavefront set intersecting
$\Gs(-\l)$. 
\end{rem}


\bibliographystyle{plain}
\bibliography{sm}
\end{document}